\documentclass[11pt]{article}
 \usepackage{epsfig}
 \usepackage{amssymb,amsmath,amsthm, amscd, amsfonts}
\usepackage{latexsym}
\usepackage[normal]{subfigure}

\renewcommand{\eqref}[1]{(\ref{#1})}
\usepackage[numbers,sort&compress]{natbib}

\newcommand{\bsub}{\begin{subequations}}
\newcommand{\esub}{\end{subequations}$\!$}

\newcommand{\ds}[0]{\displaystyle}

\newcommand{\bey}{\begin{eqnarray}}
\newcommand{\eey}{\end{eqnarray}}

\newcommand{\beq}{\begin{equation}}
\newcommand{\eeq}{\end{equation}}
\theoremstyle{plain}

\newtheorem{pro}{\hspace{6mm}Proposition}[section]

\theoremstyle{definition}

\theoremstyle{remark}

\title{Permanent charge effects on ionic flow: a numerical study of flux ratios and their bifurcation}
\date{}

\author{
Weizhang Huang\thanks{Department~of Mathematics, University of Kansas, Lawrence, Kansas 66045, USA {\tt whuang@ku.edu}} \and
Weishi Liu\thanks{Department~of Mathematics, University of Kansas, Lawrence, Kansas 66045, USA  {\tt wsliu@ku.edu}} \and
Yufei Yu\thanks{Department~of Mathematics, University of Kansas, Lawrence, Kansas 66045, USA
{\tt yh910612@gmail.com}.
}
}
\begin{document}

\maketitle

\begin{abstract}
Ionic flow carries electrical signals for cells to communicate with each other. The permanent charge of an ion channel is a crucial protein structure for flow properties while boundary conditions play a role of the driving force. Their effects on flow properties have been analyzed via a quasi-one-dimensional Poisson-Nernst-Planck model for small and relatively large permanent charges. The analytical studies have led to the introduction of flux ratios that reflect permanent charge effects and have a universal property. The studies also show that the flux ratios have different behaviors for small and large permanent charges. However, the existing analytical techniques can reveal neither behaviors of flux ratios nor transitions between small and large permanent charges. In this work we present a numerical investigation on flux ratios to bridge between small and large permanent charges.  Numerical results verify the analytical predictions for the two extremal regions. More significantly, emergence of non-trivial behaviors is detected as the permanent charge varies from small to large. In particular, saddle-node bifurcations of flux ratios are revealed, showing rich phenomena of permanent charge effects by the power of combining analytical and numerical techniques. An adaptive moving mesh finite element method is used in the numerical studies.
\end{abstract}



\section{Introduction}
\label{SEC:Introduction}
\setcounter{equation}{0}
Ion channels are large proteins embedded in membranes of cells. They serve as a major way for cells to communicate and interact with each other and with the outside world. Ion channels may open and close depending on transmembrane voltage, pressure, light, etc.
The movement of ions through channels
produces electrical signals that control many biological functions.  
Two key structures of an ion channel are its shape and permanent charge. 
An ion channel has a varying cross-section area   along its longitudinal axis and,  typically, has a relatively short and narrow neck where the permanent charge is distributed.  While the structures are crucial, it is the properties of ionic flows that are the main interest on ion channels. An equally important but sometimes overlooked factor
is the boundary conditions -- ion concentrations and electric potential on the sides of the ion channel. Their interactions with the channel structures determine specifics of ionic movement. The ultimate goal of ion channel studies is to understand the correspondence from the channel structures and boundary conditions to ionic flow properties. This is a challenging task due to the multi-scale and multi-parameter nature of the problem as well as the fact that present experimental techniques are unable to measure or observe internal (within the channel) behavior of ionic flows.
Although it has been known experimentally that ionic flows exhibit extremely rich phenomena, there still lacks so far a good set of mathematical characteristics for ionic flow properties, and thus studies of ion channel problems based
on simple models have been playing a unique role in identifying critical characteristics and separatrices among distinct behaviors.

The movement of ion species through membrane channels is affected by multiple physical quantities that interact with each other nonlinearly and  non-locally.  
The basic  models for electrodiffusion are self-consistent Poisson-Nernst-Planck (PNP) type models.  Those models consider open stage of channels and treat the medium implicitly as dielectric continuum.
They are not direct limits of molecular dynamic models as the number of ions approaches infinite.
They miss details of  motions of individual ions but capture thermodynamic quantities of the ionic flow
such as fluxes, pressure, and energy.
PNP systems can be viewed as the Fokker-Planck systems of molecular dynamic models \cite{SNE2001} coupled with the Poisson equation for the drift (electric field) that is a part of unknown state variables, and they can also be  derived from Boltzmann equations \cite{B1992} or energetic variational principles \cite{HFEL2012, HEL2014}.

Rigorous analysis has the advantage to discover important properties of biological interest and provides detailed classifications of distinct behaviors over different physical domains, in limiting or ideal setups. Numerical simulation has the power to extend  the analytical discovery to realistic parameter ranges of physical problems, and often, discover further phenomena along the continuation. This is the methodology  of this work.  We consider here open channels with fixed  shape and permanent charge distribution and combine the advantages of analysis and numerics to  examine the effects of permanent charges on individual fluxes. 
More precisely,  previous analysis based on PNP has revealed a number of interesting, some counterintuitive, phenomena of permanent charge effects for small  and large permanent charges \cite{JLZ2015, ZEL2019}.   
For channels with  permanent charge density that is small relative to the characteristic concentrations, it has been widely known that the current is increasing with respect to the transmembrane electric potential. The saturation effect due to large permanent charge density has also been established recently \cite{ZEL2019}. It is still unclear how the flux of each species is influenced by the electrochemical potential interacting with the permanent charge. It seems intuitive that the permanent charge always promotes the fluxes of counter-ion species (those with opposite charge signs as the permanent charge), and reduces the fluxes of co-ion species (those with the same charge sign as the permanent charge). However, the intuition is incorrect, as shown in \cite{JLZ2015}. Besides, it seems reasonable that if the  transmembrane electrochemical potential is increased in magnitude, the flux through the channel should increase. This has been proven false too; more precisely,  it is shown in \cite{ZEL2019} that the counter-ion flux   can decrease   as the transmembrane electrochemical potential increases.  This declining phenomenon seems to be well-known but its mechanism due to large permanent charge is discovered in \cite{ZEL2019} only recently.  

Due to the nonlinear, multi-scale, and multiple parameter nature of the problem, it is unrealistic to expect that effects of permanent charges interacting with boundary conditions can be fully understood with the analytical approach for all magnitudes of permanent charges. In this work, incorporating with the analytical results we take the advantage of a numerical method (an adaptive moving mesh finite element method \cite{HR11}) to obtain a complete diagram of the effects on ionic flows of the interplay between the permanent charge and the   electric potential.  It is well recognized that the numerical computation of PNP problems itself is challenging too for many reasons. The problem is generally stiff and, more importantly, lacks a set of concrete mathematical/physical characteristics to simulate. Motivated by available analytical results, we numerically study a  particular quantity, {\em the flux ratio}   introduced in \cite{Liu18}, and its dependence on permanent charges and  electric potentials. The analysis in \cite{JLZ2015, ZEL2019, ZL2020} on flux ratios for small and large permanent charge 
for two species reveals detailed specific dependence of flux ratio on boundary conditions and also a universal property. These analytical results provide a crucial  insight and a starting point of our numerical study. The adaptive moving mesh finite element method shows to be very suitable for such a problem. Not only it reproduces the known analytical results accurately, but also provides a global picture of flux ratio behavior including intermediate size permanent charges. In particular, our numerical study allows a rather complete partition of the parameter space into domains according to different characters of flux ratios.


The rest of the paper is organized as follows. In \S\ref{model_lambda}, we describe a three-dimensional PNP model, and a quasi-one-dimensional PNP model and its dimensionless form for ionic flow.    In \S\ref{SEC:recentresults}, we recall  the flux ratio for permanent charge effects introduced in \cite{Liu18} and some analytical results from \cite{JLZ2015,ZEL2019,ZL2020} where the dependence of flux ratio on small and large permanent charges along with the boundary conditions has been fully studied for the flow of two ion species.
In \S\ref{SEC:newresults}, we study the dependence of flux ratio on general values of permanent charge. With fixed boundary concentrations of both ion species, we examine the dependence of the flux ratio on permanent charge and on 
the boundary condition of the electric potential. Combining these results, we give a complete diagram for the effects of permanent charges on the flux of each ion species.
In \S\ref{SEC:numerical}, we describe the adaptive moving mesh finite element method that is used to obtain the numerical results in \S\ref{SEC:newresults} for PNP models. Mesh adaptation is necessary for efficient and accurate computation since the solutions to the PNP models have sharp layers caused by sudden changes in channel shape and permanent charge.

\section{Ionic flow   and Poisson-Nernst-Planck  models}\label{model_lambda}
\setcounter{equation}{0}

 PNP type systems are primitive models for ionic flows that treat the aqueous medium (in which salts are dissolved into free ions and ions are migrating) as dielectric continuum. For an ionic mixture with $n$ ion species, PNP reads as
 \begin{align}\label{model-3D}\begin{cases} \begin{split} 
&- \nabla\cdot\Big(\varepsilon_r({\bf r})\varepsilon_0 \nabla\Phi\Big)=e_0\Big(\sum_{j=1}^nz_j C_j+{\cal Q}({\bf r})\Big),\qquad {\bf r}\in \Omega\\
 &\nabla\cdot \vec{\cal J}_k =0, \quad -\vec{\cal J}_k=\frac{1}{k_BT}{\cal D}_k({\bf r}) C_k\nabla \mu_k, \quad
 k=1,\cdots, n, \qquad {\bf r}\in \Omega
 \end{split}
 \end{cases}
\end{align}
where $\Omega$ is a three-dimensional cylinder-like domain representing the channel, ${\bf r}$ denotes the physical coordinates,  ${\cal Q}({\bf r})$ is the permanent charge density, $\varepsilon_r({\bf r})$ is the  relative dielectric coefficient, $\varepsilon_0$ is the vacuum permittivity, $e_0$ is the elementary charge, $k_B$ is the Boltzmann constant, $T$ is the absolute temperature, $\Phi$ is the electric potential,  and, for  the $k$th ion species, $C_k$ is the
concentration, $z_k$ is  the valence (the number of charges per particle), $\mu_k$ is the electrochemical potential depending on $\Phi$ and $\{C_j\}$,   $\vec{\cal J}_k$ is the flux density vector, and ${\cal D}_k({\bf r})$ is  the diffusion coefficient.  

The Poisson equation  (the first equation in (\ref{model-3D})) for the electric potential $\Phi$ is  the continuum version of Coulomb's Law. The Nernst-Planck equations (the second and third equations in (\ref{model-3D})) describe the steady state of the conservation of mass of the ionic flow. The main modeling component is the electrochemical potential $\mu_k$. It can be  split  into two components, i.e., $\mu_k = \mu_k^{id} + \mu_k^{ex}$, where $\mu_k^{id}$ and $\mu_k^{ex}$ represent ideal and excess components, respectively. The ideal or point-charge component is given by
\[
\mu^{id}_k = z_k e_0 \Phi + k_B T \ln\frac{C_k}{C_0},
\]
where $C_0$ is a characteristic concentration. The first term is the electric potential while the second term
is the ideal gas potential. The excess component $\mu_k^{ex}$ accounts for ion-to-ion interactions and ion size effects, and multiple models for $\mu_k^{ex}$ have been proposed.  A hard-sphere   model is  
\begin{equation}
\begin{split}
\frac{1}{k_BT}\mu_k^{HS}   & =  - \ln \Big( 1 - \frac{4\pi}{3}\Sigma_j  r_j^3 C_j\Big) + \frac{4 \pi r_k \Sigma_j  r_j^2 C_j}{1 - \Sigma_j \frac{4\pi}{3} r_j^3 C_j} + \frac{4 \pi r_k^2 \Sigma_j r_j C_j}{1 - \Sigma_j \frac{4\pi}{3} r_j^3 C_j}
+ \frac{4\pi}{3}  \frac{r_k^3 \Sigma_j C_j}{1 - \Sigma_j\frac{4\pi}{3} r_j^3 C_j} ,
\label{HS}
\end{split}
\end{equation}
where $r_j$ denotes the radius of the $j$th hard sphere ion species \cite{QLCL2016}.
Including a hard-sphere model improves the point-charge with the volume exclusive effect of finite ion sizes. 
For most of this work, we consider only the ideal component $\mu_k^{id}$.
In \S\ref{sec:HS}, we  examine the properties with the inclusion of the hard sphere component $\mu_k^{HS}$.


Since ion channels have narrow cross-sections relative to their lengths, the three-dimensional systems (\ref{model-3D}) can be reduced to quasi-one-dimensional models, which were first proposed in \cite{NE1998}, and a special case of the reduction is justified in \cite{LW2010}. A quasi-one-dimensional PNP boundary value problem reads as
\begin{equation}
\begin{cases}
\begin{split}
& - \frac{1}{A(X)}\frac{d }{d X} \left (\epsilon_r(X) \epsilon_0 A(X)
\frac{d \Phi}{d X}\right ) = e_0\left (\sum_{j=1}^n z_j C_j(X) + \mathcal{Q}(X)\right ),\qquad X \in (a_0, b_0)
\\
& \frac{d \mathcal{J}_k}{d X}  = 0,\quad
-\mathcal{J}_k = \frac{1}{k_B T} \mathcal{D}_k(X) \mathcal{A}(X) C_k(X) \frac{d \mu_k}{d X},\quad k = 1, ..., n,
\qquad X \in (a_0, b_0)
\end{split}
\end{cases}
\label{model1}
\end{equation}
subject to the boundary conditions \cite{EL2007}
\begin{equation}\label{BC1}
\Phi(a_0) = \mathcal{V}, \quad C_k(a_0) = \mathcal{L}_k > 0; \quad \Phi(b_0) = 0, \quad C_k(b_0) = \mathcal{R}_k > 0, 
\end{equation}
where $X$ is the coordinate along the axis of the channel, $A(X)$ is
the area of the cross section of the channel at location $X$, and $\mathcal{J}_k$ is the flux of the $k$th ion species
through the cross section of the channel.

It is convenient to work on a dimensionless form of the  quasi-one-dimensional PNP. To this end, we introduce the dimensionless variables as
\begin{equation}
\begin{cases}
\begin{split}
& \epsilon^2 = \frac{\epsilon_r \epsilon_0 k_B T}{e_0^2 (b_0 - a_0)^2 C_0}, \quad x = \frac{X - a_0}{b_0 - a_0}, \quad h(x) = \frac{A(X)}{(b_0 - a_0)^2}, \quad Q(x) = \frac{\mathcal{Q}(X)}{C_0}, \\ 
& D_k(x)\mathcal{D}_0 = \mathcal{D}_k(X), \quad \phi(x) = \frac{e_0}{k_B T} \mathcal{\phi}(X), \quad c_k = \frac{C_k(X)}{C_0}, \quad J_k = \frac{\mathcal{J}_k}{(b_0 - a_0)C_0\mathcal{D}_0}, \\            & \bar{\mu}_k(x) = \frac{1}{k_B T} \mu_k(X), \quad V = \frac{e_0}{k_B T} \mathcal{V},\quad L_k = \frac{\mathcal{L}_k}{C_0}, \quad R_k = \frac{\mathcal{R}_k}{C_0},
\end{split}
\end{cases}
\end{equation}
where $\mathcal{D}_0$ is a characteristic diffusion coefficient.
The dielectric coefficient $\epsilon$ can be rewritten as $\epsilon = \frac{\lambda_D}{b_0-a_0} > 0$, i.e., the Debye
length divided by the length between the two electrodes. 
Then the dimensionless quasi-one-dimensional PNP reads as
\begin{equation}
\begin{cases}
\begin{split}
& - \frac{1}{h(x)}\frac{d}{d x} \left ( \epsilon^2 h(x) \frac{d \phi}{d x} \right ) = \sum_{j=1}^n z_j c_j + Q(x),
\qquad x \in (0,1)
\\
& \frac{d J_k} {d x} = 0, \quad - J_k = D_k(x) h(x) c_k\frac{d\bar{\mu}_k}{dx},
\quad k = 1, \cdots, n, \quad x \in (0,1)
\end{split}
\end{cases}
\label{dimensionlessPNP}
\end{equation}
with the boundary conditions at $x = 0$ and $x = 1$
\begin{equation}
\phi(0) = V, \quad c_k(0) = L_k; \quad \phi(1) = 0 \quad c_k(1) = R_k.
\label{bound}
\end{equation}
We will assume that $\epsilon$ is constant and small. As a result, one can treat (\ref{dimensionlessPNP})
as a singularly perturbed problem with $\epsilon$ as the singular parameter.

In this work, we consider an ionic mixture of two ion species with $z_1 > 0 > z_2$.
We  assume that
the boundary conditions satisfy the electroneutrality 
\[
z_1 L_1+z_2 L_2 = 0 = z_1 R_1+z_2 R_2 .
\]
The reason for this is that, otherwise, there will be sharp boundary layers which cause significant changes  (large gradients) in the electric potential and concentrations near the boundaries
and thus lead to non-trivial uncertainties in measuring their values.
It is worth pointing out that the geometric singular perturbation framework for PNP type  models developed in \cite{EL2007,  Liu2009} can treat the situation without the electroneutrality assumption.


In the dimensionless form, the electrochemical potential can be written as $\bar{\mu}_k = \bar{\mu}_k^{id} + \bar{\mu}_k^{ex}$,
where $\bar{\mu}_k^{id} = z_k \phi + \ln c_k$ and $\bar{\mu}_k^{ex}$ is taken to be zero
for most cases we consider in this work.
In \S\ref{sec:HS}, a  hard-sphere case (cf. (\ref{HS}) or (\ref{muLHS})) will be examined and the permanent charge effects corresponding to different boundary values will be studied and characterized via numerical approach.
 
In this work, we choose  $D_1 = D_2 = 1$, $z_1 = 1$, $z_2 = -1$,   $L_1 = L_2 = L$,
and $R_1 = R_2 = R$,
where $L$ and $R$ are the parameters.
The  channel geometry is taken to  be 
\begin{equation}
h(x) = 
\begin{cases}
3\left (0.4x + 20(\frac{1}{3} - x)\right ), & \quad  0 \leq x < \frac{1}{3} \\
3\left ( 0.4(x-\frac{1}{3}) + 0.4(\frac{2}{3} - x)\right ), & \quad \frac{1}{3} \leq x < \frac{2}{3} \\
 3\left ( 20(x - \frac{2}{3}) + 0.4( 1 - x)\right ), & \quad \frac{2}{3} \leq x < 1. \\
\end{cases}
\label{shape}
\end{equation}
The variable cross-section area is chosen to reflect the fact that the channel is not uniform and
much narrower in the neck
($\frac{1}{3} < x < \frac{2}{3}$) than other regions \cite{JLZ2015}.
As in \cite{JLZ2015},   the permanent charge is assumed to have the form 
\begin{equation}
Q(x) =
\begin{cases}
0, & \quad 0 < x < \frac{1}{3} \\
2 Q_0, & \quad \frac{1}{3} < x < \frac{2}{3} \\
0, & \quad \frac{2}{3} < x < 1.
\end{cases}
\label{permanentcharge}
\end{equation}
We will denote $q_0 = 2 Q_0$ at various occasions in the rest of the paper.

In computation, the discontinuity of $Q(x)$ and $h(x)$ at $x= \frac{1}{3}$ and $x= \frac{2}{3}$ can cause
divergence in Newton's iteration when applied to the nonlinear system resulting from discretization of (\ref{dimensionlessPNP}). To avoid this difficulty, we apply a regularization on $Q(x)$ and $h(x)$ as 
\begin{equation}
Q_\delta(x) = 2 Q_0\left [ \tanh\left ( \frac{x - \frac{1}{3}}{\delta}\right ) - \tanh\left ( \frac{x - \frac{2}{3}}{\delta}\right )\right ],
\label{Qregular}
\end{equation}
\begin{equation}
h_{\delta_x}(x) = 
\begin{cases}
3[0.4x + 20(\frac{1}{3} - x)], & \quad  0 \leq x < \frac{1}{3} - \delta_x \\
14.7\delta_x^{-1} (x - \frac{1}{3} - \delta_x) + 0.4, & \quad \frac{1}{3} - \delta_x \leq x < \frac{1}{3} + \delta_x \\
3[0.4(x-\frac{1}{3}) + 0.4(\frac{2}{3} - x)], & \quad \frac{1}{3} + \delta_x \leq x < \frac{2}{3} - \delta_x \\
14.7\delta_x^{-1} (x - \frac{2}{3} + \delta_x) + 0.4, & \quad \frac{2}{3} - \delta_x \leq x <\frac{2}{3} + \delta_x \\
 3[20(x - \frac{2}{3}) + 0.4( 1 - x)], & \quad \frac{2}{3} + \delta_x \leq x < 1 \\
\end{cases}
\label{Aregular}
\end{equation}
where $\delta > 0$ and $\delta_x > 0$ are regularization parameters.
It can be verified that $Q_{\delta}(x) \to Q(x)$ as $\delta \to 0$
and $h_{\delta_x}(x) \to h(x)$ as $\delta_x \to 0$ both in $L^2$~norm.
The values $\delta_x = 10^{-7}$ and $\delta = 1/800$ have been used in our computation.

\section{A flux ratio for permanent charge effects and relative analytical results}
\label{SEC:recentresults}
\setcounter{equation}{0}

The major interest for an ion channel is in the fluxes of individual ion species. More precisely, one would like to understand how individual fluxes depend on the channel structure and boundary conditions. On the other hand, only the total current ${\cal I} = \sum_ {k = 1}^n z_k \mathcal{J}_k$ is measurable in most experiments while individual ionic fluxes $J_k$'s are difficult to measure directly. For instance, a typical way to measure the flux of  sodium in the sodium-chloride  (Na$^+$Cl$^-$)  solutions is: a small amount of a radioactive isotope of sodium is added, the flux of this isotope is measured by its radioactivity, and the flux of sodium is then estimated using the flux of the isotope; see, e.g., \cite{HK55, JEL18, Uss49b} for detail. As such, analytical and numerical studies of the PNP systems are crucial for understanding and gaining insights to flux dynamics, channel properties, and effects of permanent charge on each single species.

To study the effects of permanent charge on fluxes, the flux ratios have been recently introduced in \cite{Liu18} as
\begin{equation}
\lambda_k(Q):=\lambda_k(Q;V,L,R, h) = \frac{J_k(Q;V,L,R,h)}{J_k(0;V,L,R,h)},
\label{lambdak}
\end{equation}
where $J_k(Q;V,L,R,h)$ is the flux of the $k$th ion species with permanent charge $Q$
and the boundary conditions (\ref{bound}), and $J_k(0;V,L,R,h)$ is the flux
with the zero permanent charge and the same boundary conditions. 
 
 It has been observed in \cite{ELX2015} that, for any $Q$,  $J_k(Q;V,L,R,h)$ has the same sign as that of $(\bar{\mu}_k(0)-\bar{\mu}_k(1))$ (the difference of the electrochemical potential at $x=0$ and $x=1$), which is determined by the boundary values $(V,L,R)$. Thus $\lambda_k>0$,   independent of the permanent charge $Q$.
Note that the permanent charge enhances the flux of the $k$th ion species (in the sense  $|J_k(Q)|>|J_k(0)|$) when  $\lambda_k (Q)> 1$. As a consequence, the equation $\lambda_k=1$ determines the boundary sets of the parameter space $(Q,V,L,R,h)$ separating the regions where the permanent charge promotes or reduces the flux $J_k$, and it is thus crucial to understand those hypersurfaces determined by $\lambda_k = 1$ in the parameter space.

For $n=2$ with $z_1=1$ and $z_2=-1$,  the effects of permanent charge of the form (\ref{permanentcharge}) on flux ratios have been analyzed in detail for small and large $Q_0$ \cite{JLZ2015, Liu18, ZEL2019,ZL2020}.  
Hereafter, the notation $\lambda_k(Q_0)$ and $\lambda_k(Q)$ will be used interchangeably. 
Moreover, $\lambda_k(Q_0,V)$ will be used when the dependence on $V$ is emphasized.
The following universal property has been established  in \cite{Liu18}.

\begin{pro}
For general $Q(x)>0$ (not necessarily in the form  (\ref{permanentcharge})),  one always has  $\lambda_1(Q) < \lambda_2(Q)$, independent of boundary conditions and channel geometry $(V,L,R,h)$.
 \label{universal-property}
 \end{pro}
 
On the other hand, the effects of permanent charge on fluxes are shown to depend on $(V,L,R,h)$ in a very rich way \cite{JLZ2015, ZEL2019,ZL2020}. 
For example, it seems intuitive that, if $Q(x)>0$, then $\lambda_1(Q) < 1 < \lambda_2(Q)$. But this is not always the case. 
Moreover, there are significant behavioral mismatches between $\lambda_k(Q_0)$'s of small $Q_0$ and those of large $Q_0$ for certain choices of $(V,L,R,h)$, which imply that changes in behavior occur  as $Q_0$ varies from small to large.

We note that results for $Q(x)<0$ can be obtained by symmetry from those with $Q(x)>0$. For example, we have $\lambda_1(Q) > \lambda_2(Q)$
for $Q(x)<0$. For this reason and for notational simplicity, we consider only nonnegative permanent charge $Q(x) \ge 0$ in the rest of the paper.

In the following we describe some analytical results on $\lambda_k(Q_0)$ for small and large $Q_0$
from \cite{JLZ2015,ZEL2019,ZL2020}.

{\bf Flux ratios for small $Q_0$.} In \cite{JLZ2015}, based on a result in  \cite{EL2007} from geometric singular perturbation analysis on  the boundary value problem (BVP) of PNP (\ref{dimensionlessPNP}) and (\ref{bound}),  expansions of $J_k$'s   in small $Q_0$ are obtained as
\begin{equation}
J_k(Q_0) = J_{k0}^0 + J_{k1}^0Q_0 +  \mathcal{O}(Q_0^2),
\label{JKsmallQ}
\end{equation}
where
\begin{align}
\label{smallQs}
\begin{cases}
\begin{split}
J_{10}^0=& \frac{(L-R)(V+\ln L-\ln R)}{H(1)(\ln L-\ln R)},\quad J_{11}^0 =\frac{A((B-1)V +\ln L-\ln R)}{2H(1)(\ln L-\ln R)^2}(V+\ln L-\ln R),\\
J_{20}^0=& \frac{(L-R)(-V+\ln L-\ln R)}{H(1)(\ln L-\ln R)}, \quad
J_{21}^0 =-\frac{A((1-B)V +\ln L-\ln R)}{2H(1)(\ln L-\ln R)^2}(-V+\ln L-\ln R),
\end{split}
\end{cases}
\end{align}
\begin{align}
\label{smallQs+1}
\begin{cases}
\begin{split}
&    H(x) = \int_0^x h^{-1} (s) ds, \quad  \alpha = \frac{H(a)}{H(1)},  \quad  \beta = \frac{H(b)}{H(1)},\\
& A = \frac{(\beta - \alpha)(L - R)^2}{((1-\alpha)L + \alpha R)((1-\beta)L + \beta R)},  \quad B = \frac{(1 - \beta) L + \beta R - \ln((1 - \alpha) L + \alpha R)}{A}.
\end{split}
\end{cases}
\end{align}
Here, we use the superscript $0$ for the case of small $ Q_0$. Note that $J_k(0)=J_{k0}^0$. Hence, for $Q_0$ small, 
\[
\lambda_k(Q_0)=1+\frac{J_{k1}^0}{J_{k0}^0}Q_0+\mathcal{O}(Q_0^2).
\]
As shown in \cite{JLZ2015}, for small $Q_0$, the parameters $(V,L,R,\alpha,\beta)$ can be decomposed  as the union of three types of regions defined by $1 < \lambda_1 < \lambda_2$, $\lambda_1 < 1 < \lambda_2$, and $\lambda_1 < \lambda_2 < 1$. More precisely, let
\begin{equation}
  t = L/R , \quad  \gamma(t) = \frac{t \ln t - t + 1}{(t - 1) \ln t},
\label{preresult}
\end{equation}
\begin{equation}
g(\beta) =  \left ( (1-\alpha) t + \alpha \right ) \left ((1-\beta) t + \beta \right )
\ln t \ln \frac{(1-\beta) t + \beta}{(1-\alpha) t + \alpha}  + (\beta - \alpha ) (t-1)^2,
\label{g-fun}
\end{equation}
\begin{equation}
V_1^0 = V_1^0(L,R, \alpha, \beta) = - \frac{\ln L - \ln R}{z_2(1-B)}, \quad V_2^0 = V_2^0(L,R, \alpha, \beta) = - \frac{\ln L - \ln R}{z_1(1-B)}.
\label{Vq12}
\end{equation}
It is shown in \cite{JLZ2015} that $g(\beta) = 0$ has a unique solution $\beta_1 \in (\alpha, 1)$ when $\alpha < \gamma(t)$.
Then,  for small $Q_0$ with $ t = L/R > 1$, one has the following results.
\begin{itemize}
\item If $\alpha < \gamma(t)$ and $\beta \in (\alpha, \beta_1)$, then 
\[
V_1^0 < 0 < V_2^0  \qquad \text{and} \qquad
\begin{cases}
1 < \lambda_1 < \lambda_2, &\mbox{ for } V < V_1^0 \\
\lambda_1 < 1 < \lambda_2, &\mbox{ for } V_1^0 < V < V_2^0 \\
\lambda_1 < \lambda_2 < 1, &\mbox{ for } V > V_2^0. 
\end{cases}
\]

\item If $\alpha < \gamma (t)$ and $\beta > \beta_1$ or $\alpha \geq \gamma(t)$, then
\[
V_1^0 > 0 > V_2^0  \qquad \text{and} \qquad
\begin{cases}
1 < \lambda_1 < \lambda_2, &\mbox{ for } V > V_1^0 \\
\lambda_1 < 1 < \lambda_2, &\mbox{ for } V_2^0 < V < V_1^0 \\
\lambda_1 < \lambda_2 < 1, &\mbox{ for } V < V_2^0 .
\end{cases}
\]
\end{itemize}

Similar analytical results for $t = L/R < 1$ have been obtained. The interested reader is referred to
\cite{JLZ2015} for detail.

We note that, with $(L, R,\alpha,\beta)$ fixed, along the curves $\lambda_k(Q_0,V) = 1$ in the $Q_0$-$V$ plane, we have $V\to V_k^0$ as $Q_0 \to 0$.

 {\bf Flux ratios for large $Q_0$.}  Flux ratios and other effects have been examined for large $Q_0$ in \cite{ZEL2019}.  In particular, for $\nu=1/{Q_0}\ll 1$,  the following expansions of the flux $J_k$'s with respect to $\nu$ have been rigorously established,
\begin{equation}
J_k(\nu) = J_{k0}^{\infty} + J_{k1}^{\infty} \nu + \mathcal{O}(\nu^2),
\label{JKbigQ}
\end{equation}
where  
\begin{equation}
J_{10}^{\infty}  = 0, \quad J_{11}^{\infty} =\frac{1}{2H(1)(\beta-\alpha)}\left(\frac{(1-\beta)L+\alpha R}{(1-\beta)\sqrt{e^VL}+\alpha\sqrt{R}}\right)^2(e^VL-R), 
\label{zeroJs}
\end{equation}
\begin{align}\label{J11J21}
\begin{cases}
\begin{split}
 J_{20}^{\infty}  =& \frac{2 \sqrt{LR}}{H(1)} \frac{1}{(1 - \beta) \sqrt{L} + \alpha \sqrt{e^{-V}R}} (\sqrt{e^{-V}L} - \sqrt{R}),\\
  J_{21}^{\infty} =&-\frac{(\beta-\alpha)e^VLR\big((1-\beta)L+\alpha R\big)}{H(1)\big((1-\beta)\sqrt{e^VL}+\alpha\sqrt{R}\big)^3}
 (\sqrt{e^{-V}L}-\sqrt{R}) \\
 &+\frac{(e^VL-R)\big(-V+\ln L-\ln R\big)\big((1-\beta)L+\alpha R\big)^3}{4(\beta-\alpha)H(1)(\sqrt{e^{-V}L}-\sqrt{R})\big((1-\beta)\sqrt{e^VL}+\alpha\sqrt{R}\big)^3}\\
&-\frac{e^VL-R}{2(\beta-\alpha)H(1)}\left(\frac{(1-\beta)L+\alpha R}{(1-\beta)\sqrt{e^VL}+\alpha\sqrt{R}}\right)^2.
\end{split}
\end{cases}
 \end{align}
Here, we use the superscript $\infty$ for the case of large $Q_0$.
Recall from (\ref{smallQs}) that $J_k(0)=J_{k0}^0$. For cations, one then has, for $Q_0\gg 1$ or ${\nu}=1/{Q_0}\ll 1$,
\begin{align*}
\lambda_1(Q_0)=\frac{J_1(Q_0)}{J_1(0)}= \frac{J_{11}^{\infty}}{J_{10}^0}\nu+\mathcal{O}(\nu^2)<1.
\end{align*}
Thus, large positive permanent charges inhibit the flux of cations. This contrasts sharply to the situation of small positive permanent charge where, under some boundary conditions, it could enhance the flux of cations.

 For anions, using (\ref{J11J21}) for $J_{20}^{\infty}$ and (\ref{smallQs}) for $J_2(0)=J_{20}^0$,  one has, for $Q_0\gg 1$ or ${\nu}=1/{Q_0}\ll 1$,
 \begin{align}\label{fK}
 \lambda_2(Q_0)=&\frac{J_{2}(Q_0)}{J_2(0)}= \frac{J_{20}^{\infty}}{J_{20}^0} +\mathcal{O}(\nu)
 =\frac{2(t-\sqrt{e^Vt})\ln t}{(t-1)((1-\beta)\sqrt{e^Vt}+\alpha)(\ln t-V)}+\mathcal{O}(\nu),
 \end{align}
 where $t=L/R$. It follows from Lemma 4.6 in \cite{ZL2020} that, there exist $V_k^{\infty}=V_k^{\infty}(L,R,\alpha, \beta)$ with $V_1^{\infty}<V_2^{\infty}$ such that, for large $Q_0$, one has 
$\lambda_2(Q_0)>1$ if $V\in (V_1^{\infty},V_2^{\infty})$, and $\lambda_2(Q_0)\le 1$ otherwise.
 Explicit analytical formulas are unavailable for $V_1^{\infty}$ and $V_2^{\infty}$.
 See \S\ref{sec:Valternate} for their numerical values in case studies.
 
 \vspace{20pt}
 
 From the above results, we have seen that the behaviors of $\lambda_k$'s are very different for small and large $Q_0$.
 In particular, for the situation with $n=2$, $z_1=1$, $z_2=-1$ and fixed $(L, R, h)$, $\lambda_k$'s are functions
of $V$ and $Q_0$. For small $Q_0$, as $V$ increases, the change occurs from $1 < \lambda_1 < \lambda_2$
to $\lambda_1 < 1 < \lambda_2$ and to $\lambda_1 < \lambda_2 < 1$
or from $\lambda_1 < \lambda_2 < 1$ to $\lambda_1 < 1 < \lambda_2$  and to $1 < \lambda_1 < \lambda_2$, depending on the values of $\alpha$ and $\beta$.
On the other hand, for large $Q_0$, we have $\lambda_1 < 1$ independent of $V$ while $\lambda_2 > 1$ for
a bounded interval of $V$ values. In other words, $(\lambda_1, \lambda_2)$ has the region $\lambda_1 < \lambda_2 < 1$ for large
$Q_0$ and large $V$ and region $\lambda_1 < 1 < \lambda_2$ for large $Q_0$ and a bounded interval of $V$.
This indicates that there are changes in the behaviors of $\lambda_k$ between small and large $Q_0$.
 Our main objective of the current work is to study the specifics of these transitions and their mechanisms,
 which seems beyond existing analytical techniques. Our strategy is to allow $V$ to vary as well as $Q_0$:
 we shall study the effects of $Q_0$ and $V$ on $\lambda_1$ and $\lambda_2$  with several choices of fixed boundary conditions $L$, $R$ and channel geometry $h$.    
 For each fixed $(L,R,h)$, we take the advantage of numerical approach to identify  the parameter boundary curves defined by $\lambda_k = 1$, which partitions the $Q_0$-$V$ plane into three different types of regions. From this partition of the Q$_0$-V space, one can read out occurrence
of bifurcations of $\lambda_k(Q_0,V)=1$ among other features. It seems that a typical bifurcation is
the saddle-node bifurcation.  We shall also see that, for a  choice of $(V,L,R,h)$,  even if there is
no mismatch between the effect of small $Q_0$ and that of large $Q_0$, there are still (multiple) transitions.
 

\section {Numerical study of flux ratios for general $Q_0$}
\label{SEC:newresults}
\setcounter{equation}{0}

In this section we conduct numerical studies  on quantitative  and qualitative behaviors of flux ratios $\lambda_k$ for  permanent charge effects on fluxes. Recall that $\lambda_k=\lambda_k(Q_0; V, L, R, h)$ depends on the multiple parameters $(V, L, R, h)$ in additional to $Q_0$.  As mentioned in the introduction that our numerical study is incorporated with the analytical results in \cite{JLZ2015,ZEL2019,ZL2020}. More precisely, the flux ratios for small and large permanent charges are well understood not only qualitatively but also  quantitatively. On the other hand, it is nearly impossible for analytical analysis to provide even qualitative results for moderate size of permanent charge.
In this section  we present numerical studies that bridge the two extremes of very small and large $Q_0$ and demonstrate a complete picture of how the combination of $(Q_0,V)$ affects the flux ratios, with boundary conditions given. We also present multiple numerical results of $\lambda_1$ and $\lambda_2$, with fixed boundary concentrations of species. Throughout this section, unless stated otherwise, we take $\epsilon = 10^{-5}$, $L = 0.008$, and $R = 0.001$, which are scaled values associated with common boundary values of the concentration in practice. 

The PNP system (\ref{dimensionlessPNP}) has been solved using an adaptive moving mesh finite element method.
The method employs piecewise linear finite elements with mesh adaptation
to increase numerical resolution in the regions with discontinuities in permanent charges.
To keep the flow of the discussion, we leave the description of the numerical method to \S\ref{SEC:numerical}.
 
\subsection{Dependence of $\lambda_1$ and $\lambda_2$ on $Q_0$ for fixed $V$}
\label{sec:Qalternate}

\begin{figure}[htbp]

\subfigure[$V=50$]{\includegraphics[width = 0.47\textwidth]{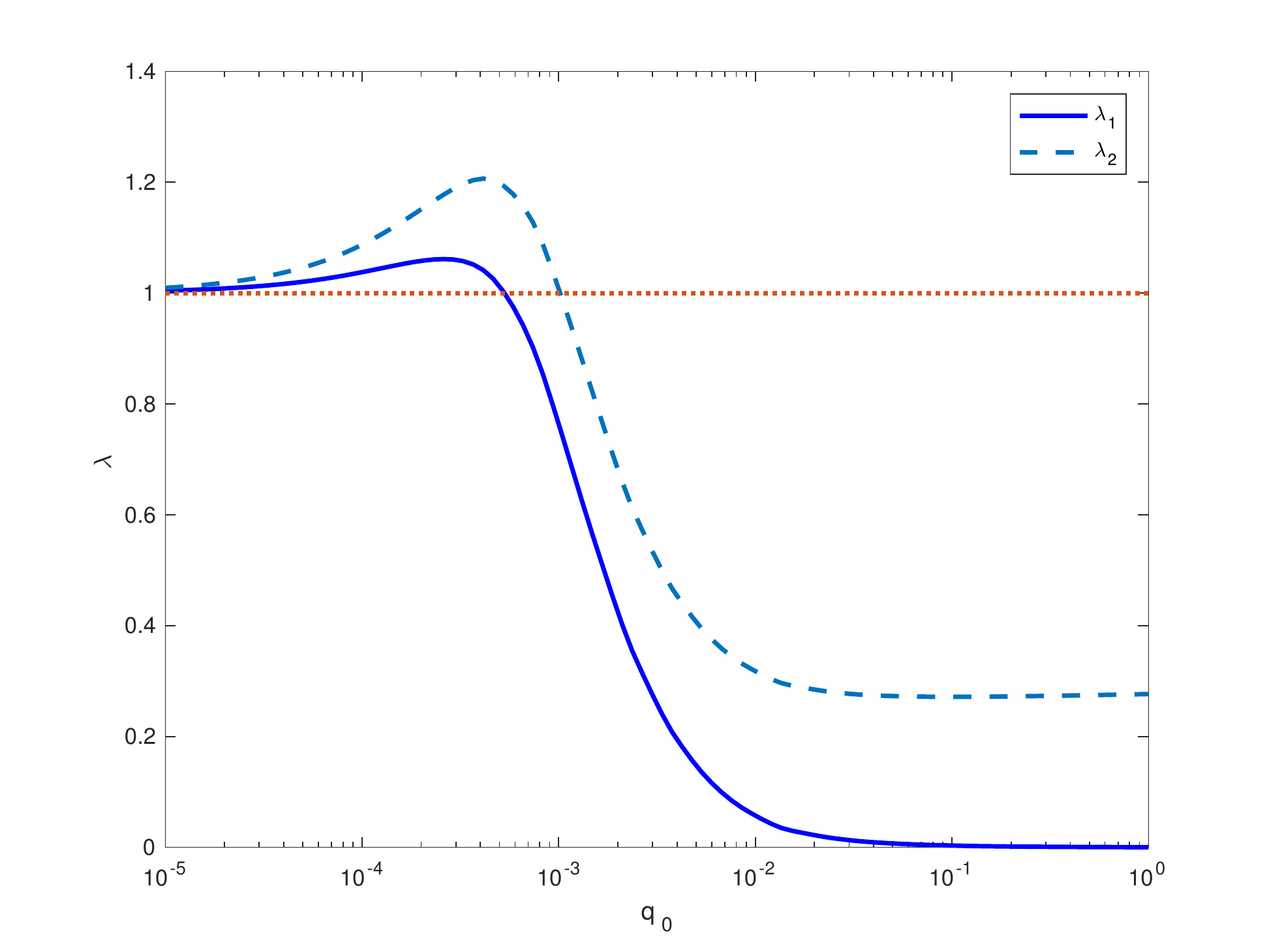} \label{fig:lambda-Q-V50}}\qquad
\subfigure[$V=10$]{\includegraphics[width = 0.47\textwidth]{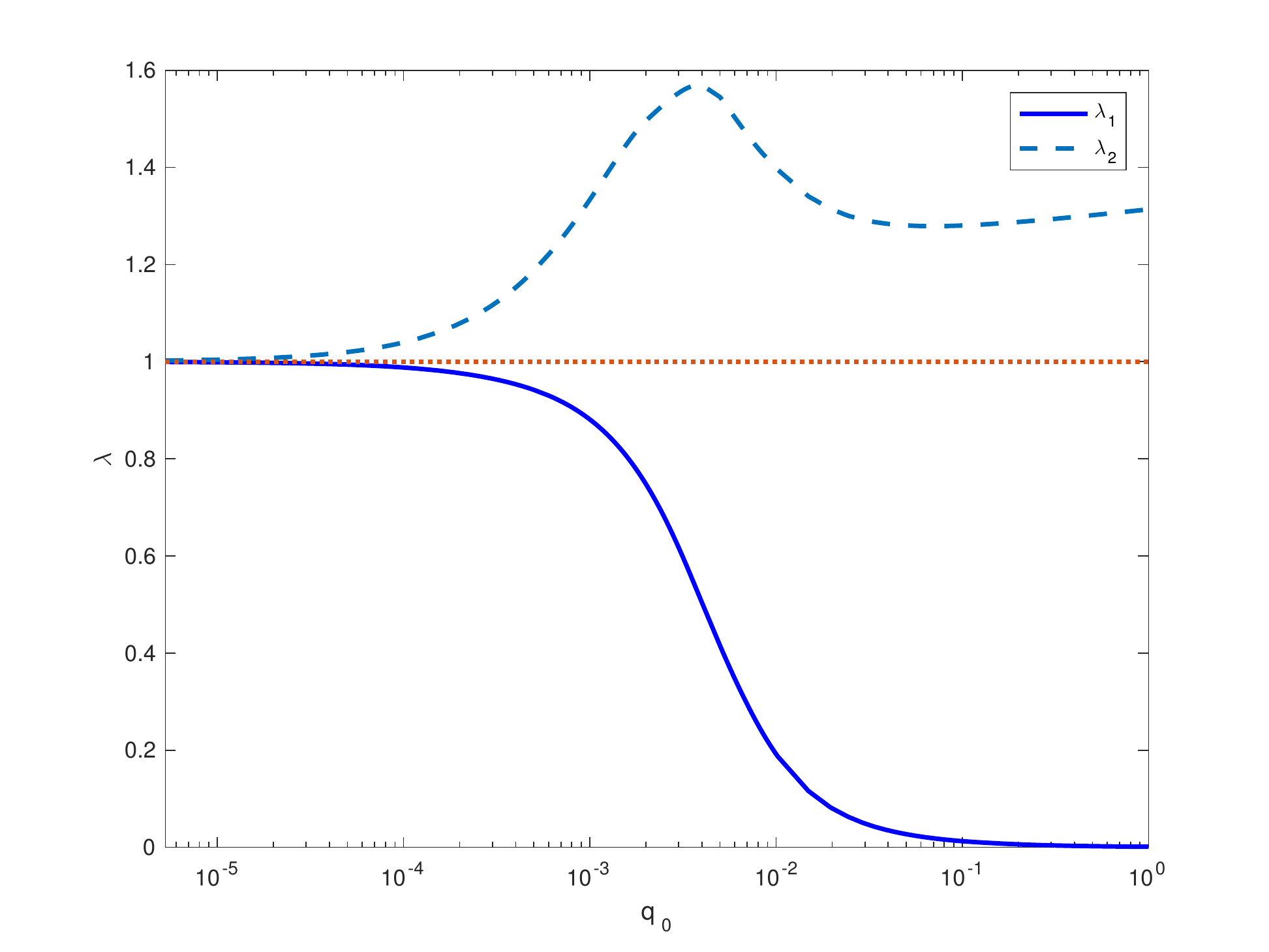} \label{fig:lambda-Q-V10}}\qquad
\subfigure[$V=-60$]{\includegraphics[width = 0.47\textwidth]{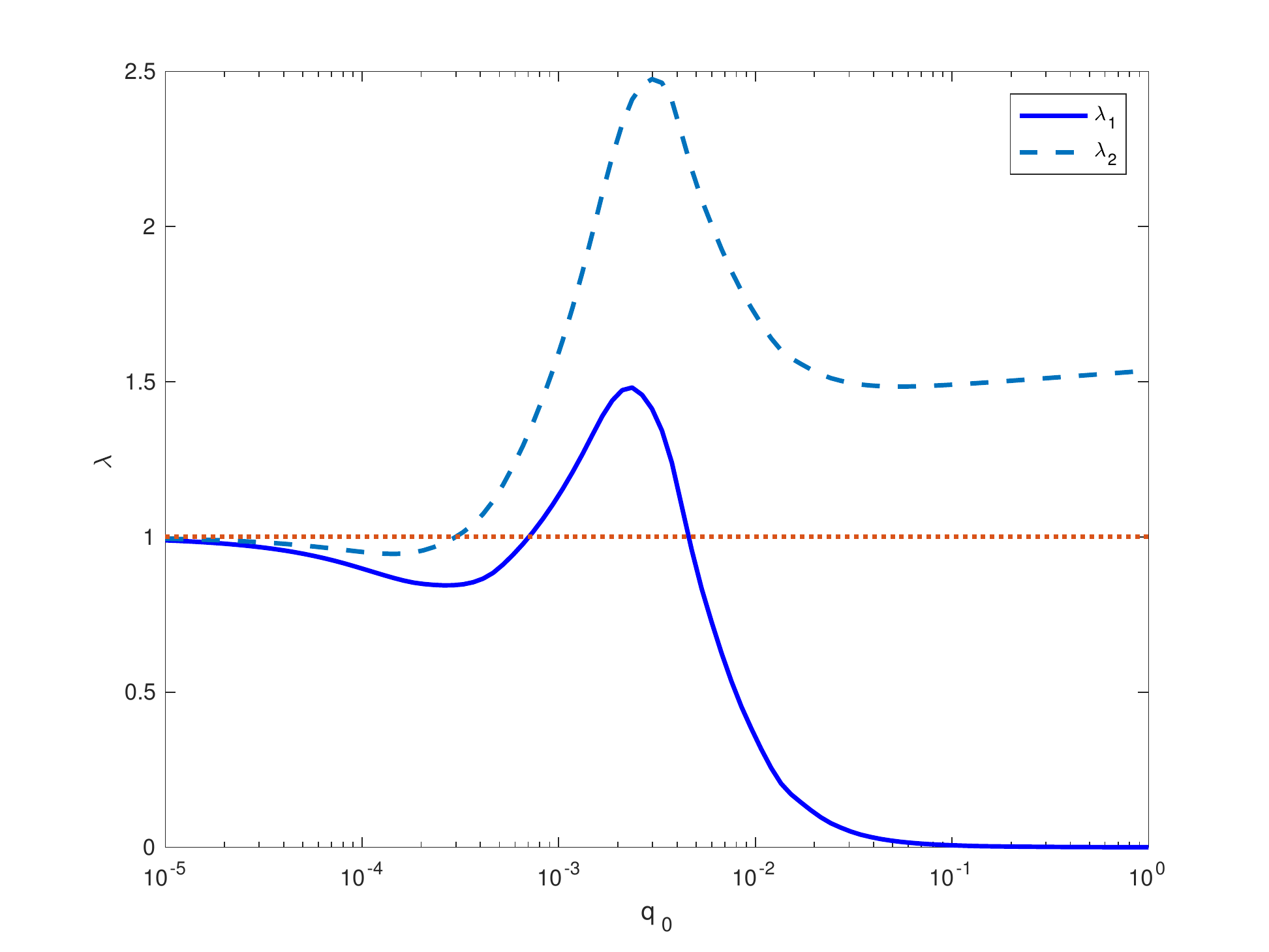} \label{fig:lambda-Q-V-60}}\qquad
\subfigure[$V = -110$]{\includegraphics[width = 0.47\textwidth]{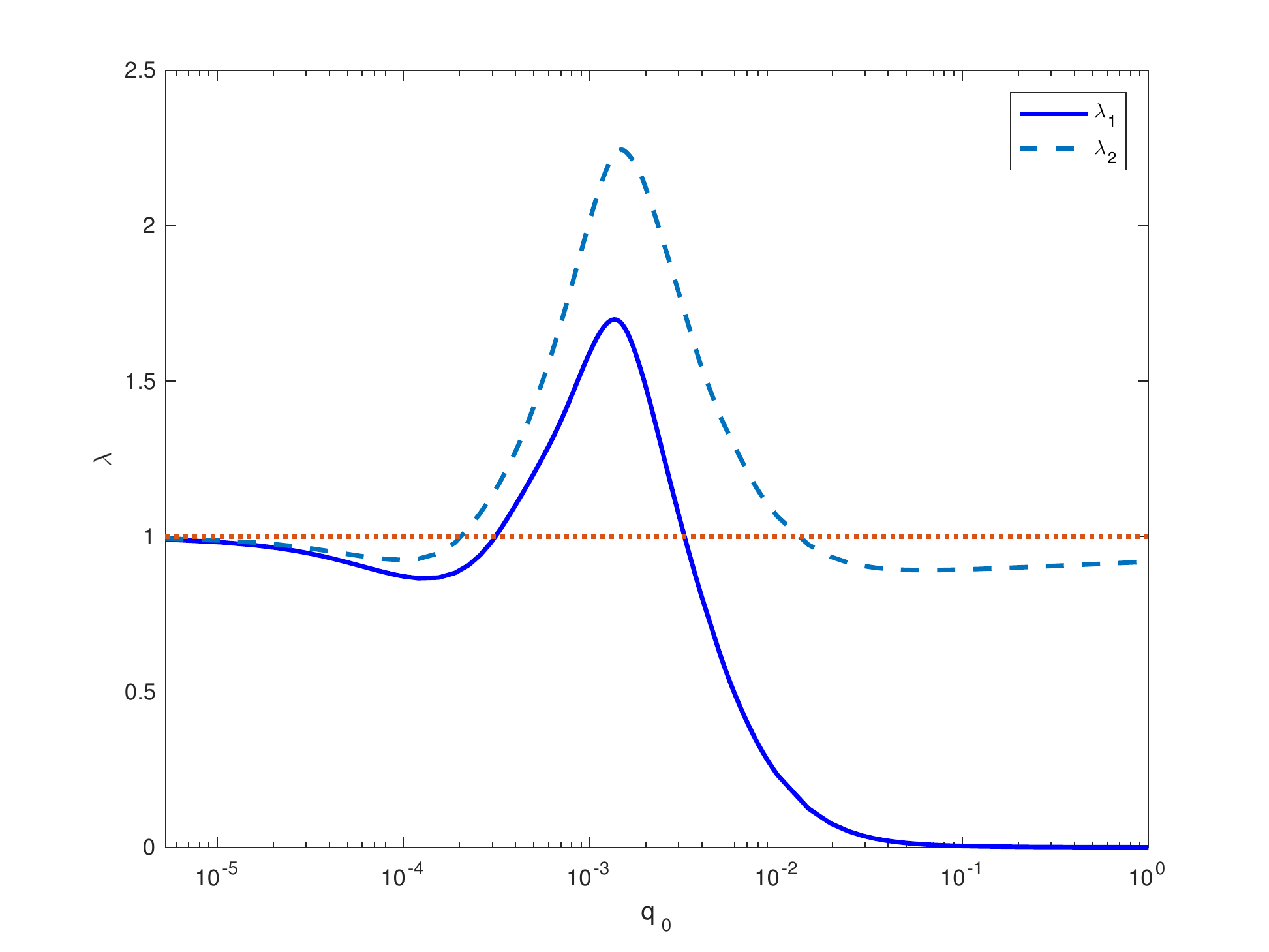} \label{fig:lambda-Q-V-110}}
\caption{$\lambda_1$ and $\lambda_2$ are plotted as functions of $q_0 = 2 Q_0$, with boundary values of the system (\ref{dimensionlessPNP}) chosen as $L = 0.008$, $R = 0.001$, and $V=50$, $10$, $-60$, or $-110$.}\label{fig:lambda-Q}
\end{figure}

In this subsection we fix $(L, R, h)$ and examine how permanent charge effects on fluxes are interacted with the role of the electric potential $V$.
Here we choose four values for $V$ and examine the dependence of $\lambda_k$ on $Q_0$.
Fig.~\ref{fig:lambda-Q-V50}, \ref{fig:lambda-Q-V10}, \ref{fig:lambda-Q-V-60} and \ref{fig:lambda-Q-V-110} show the trend of $\lambda_1$ and $\lambda_2$ for $V = 50$, $10$, $-60$, and $-110$, respectively.

In Fig.~\ref{fig:lambda-Q-V10}, starting from $\lambda_1 = \lambda_2 = 1$ at $Q_0=0$, one can observe that $\lambda_1 < 1$ and $\lambda_2 > 1$ for $q_0 \equiv 2 Q_0 \in (0,1]$. This is the case coinciding with one's intuition: the positive permanent charge helps promote the flux of the negative ion species, and inhibits that of the positive ion species. However, Fig.~\ref{fig:lambda-Q-V50}, \ref{fig:lambda-Q-V-60} and \ref{fig:lambda-Q-V-110} demonstrate that $\lambda_1 < 1 < \lambda_2$ is not the only case that we can obtain. For example, in Fig.~\ref{fig:lambda-Q-V50} (with $V = 50$), both $\lambda_1$ and $\lambda_2$ stay above 1 for small $Q_0$. As $Q_0$ increases, $\lambda_1$ crosses the value $1$ and becomes
less than $1$ and then $\lambda_2$ follows. In Fig. \ref{fig:lambda-Q-V-60} ($V = - 60$) and \ref{fig:lambda-Q-V-110} ($V = - 110$), the situation becomes more complicated. For small $Q_0$, we have $\lambda_1 < \lambda_2 < 1$. As $Q_0$ increases, the graph of $\lambda_1$ and $\lambda_2$ can cross the value $1$ multiple times. Indeed, $\lambda_1$ goes up above 1 and then becomes less than 1 for large $Q_0$. On the other hand, $\lambda_2$ goes up above 1 and
then stays there for $V = -60$ or becomes less than 1 for $V = -110$ for large $Q_0$.
These results show all the three possibilities as described in \S\ref{SEC:recentresults}: $\lambda_1 < \lambda_2 < 1$, $\lambda_1 < 1 < \lambda_2$, as well as $1 < \lambda_1 < \lambda_2$.

We remark that, for $Q_0>0$ small, Fig.~\ref{fig:lambda-Q} agrees well with the analytical prediction in \cite{JLZ2015} that has been described in \S\ref{SEC:recentresults}. With $L = 0.008$, $R = 0.001$, $a = 1/3$, and $b = 2/3$, we have $V_1^0 = -V_2^0 = 18.97$, where $V_1^0$ and $V_2^0$ are defined as in (\ref{Vq12}). One can calculate that $\alpha = 0.07$, $\beta = 0.93$, and $\beta_1 = 0.89$, where $\alpha$ and $\beta$ are as described in (\ref{smallQs}) and $\beta_1$ is the root of $g(\beta)$ defined in (\ref{g-fun}). Thus the set of boundary conditions corresponds to the case where $t = L/R > 1$, $\alpha < \gamma(L/R)$ (see (\ref{preresult})), and $\beta > \beta_1$. As described in \S\ref{SEC:recentresults}, for small $Q_0$,
\[
\begin{cases}
1 < \lambda_1 < \lambda_2, & \mbox{ for } V > 18.97 \\
\lambda_1 < 1 < \lambda_2, & \mbox{ for } -18.97 < V < 18.97 \\
\lambda_1 < \lambda_2 < 1, & \mbox{ for } V < - 18.97 .
\end{cases}
\]
Fig.~\ref{fig:lambda-Q} shows that, for small $Q_0$,
$1 < \lambda_1 < \lambda_2$ in Fig.~\ref{fig:lambda-Q-V50} ($V = 50$),
$ \lambda_1 < 1 < \lambda_2$ in Fig.~\ref{fig:lambda-Q-V10} ($V = 10$),
and $ \lambda_1 < \lambda_2 < 1$ in Fig.~\ref{fig:lambda-Q-V-60} ($V = -60$) and \ref{fig:lambda-Q-V-110} ($V = -110$).
Thus, numerical results for $Q_0$ small agree well with those analytical results in \cite{JLZ2015}.

For large $Q_0$ ($q_0\equiv 2 Q_0$ near 1 in Fig.~\ref{fig:lambda-Q}), the numerical result shows that $\lambda_1 \sim 0$ for all $V$ while $\lambda_2 < 1$ for $V = 50$ (Fig.~\ref{fig:lambda-Q-V50}) and $V = -110$ (Fig.~\ref{fig:lambda-Q-V-110})
and above 1 for $V = 10$ (Fig.~\ref{fig:lambda-Q-V10}) and $V = -60$ (Fig.~\ref{fig:lambda-Q-V-60}).
This is consistent with the analytical analysis in \cite{ZEL2019} (cf. \S\ref{SEC:recentresults})
which states that, for large $Q_0$, $\lambda_1 < 1$ independent of $V$ and $\lambda_2$ can becomes above 1 for a bounded interval of $V$.

The non-monotone behaviors of $\lambda_k$'s in $Q_0$ shown in Fig.~\ref{fig:lambda-Q} are generally
hard for analysis to predict. Particularly interesting is the crossing of $\lambda_k=1$ over the range $(10^{-4}, 10^{-1})$ in Fig.~\ref{fig:lambda-Q-V-60} and Fig.~\ref{fig:lambda-Q-V-110} and may deserve further investigation about its mechanism.

Fig.~\ref{fig:lambda-Q-V10} and Fig.~\ref{fig:lambda-Q-V-110} together present strong evidence of highly nonlinear interaction between $Q_0$ and the boundary condition $V$ (since $L$ and $R$ are fixed for these studies).  More precisely, $\lambda_k(Q_0, 10)$ and $\lambda_k(Q_0,-110)$ have significantly different properties as functions of $Q_0$. Similarly, in the next subsection, when $Q_0$ is fixed at different values, $\lambda_k$ as functions of $V$ have significantly different behaviors.

\begin{figure}[htbp]
\subfigure[$V = 50$]{\includegraphics[width = 0.47\textwidth]{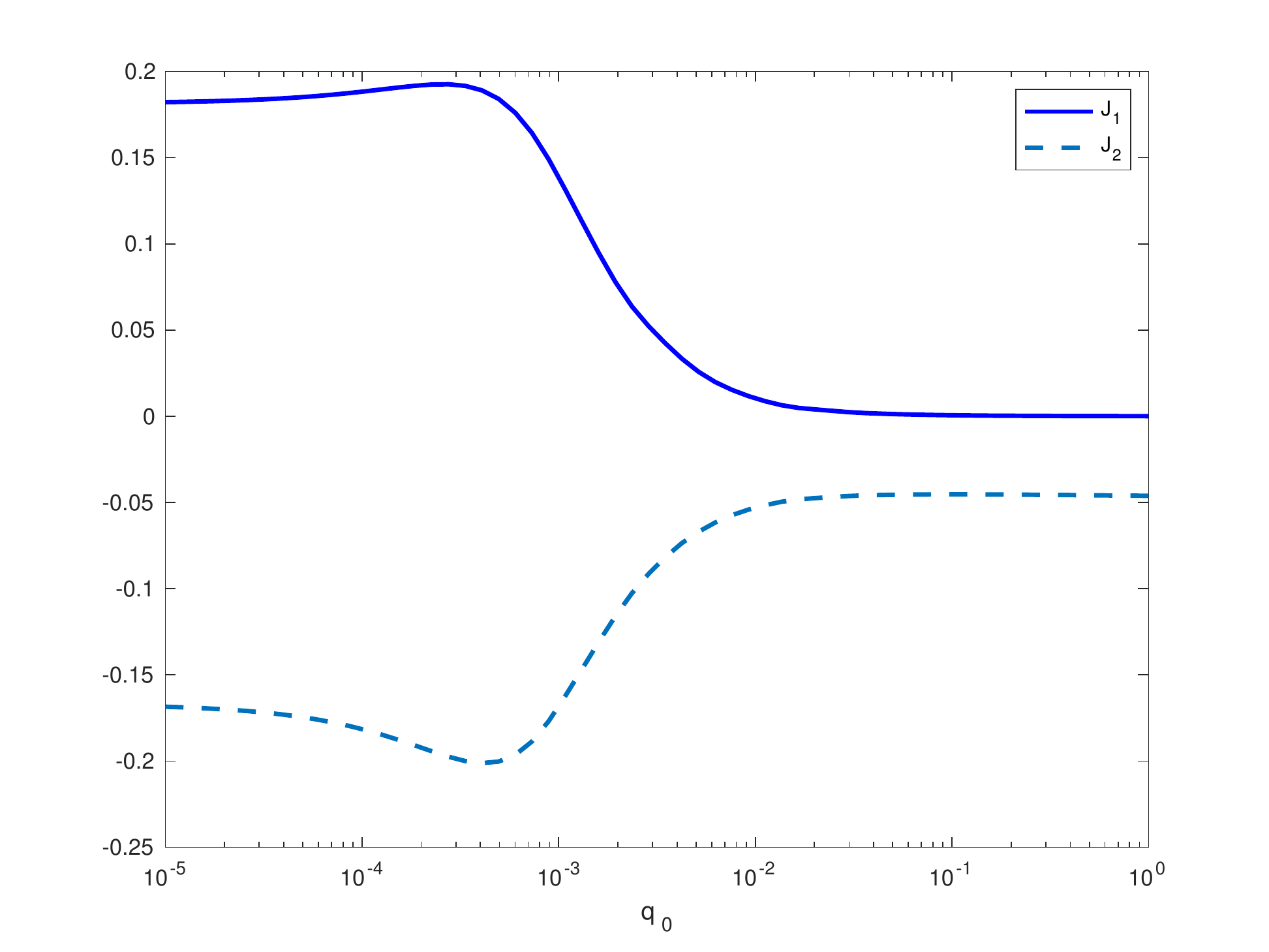} \label{fig:JV50}}\qquad
\subfigure[$V = 10$]{\includegraphics[width = 0.47\textwidth]{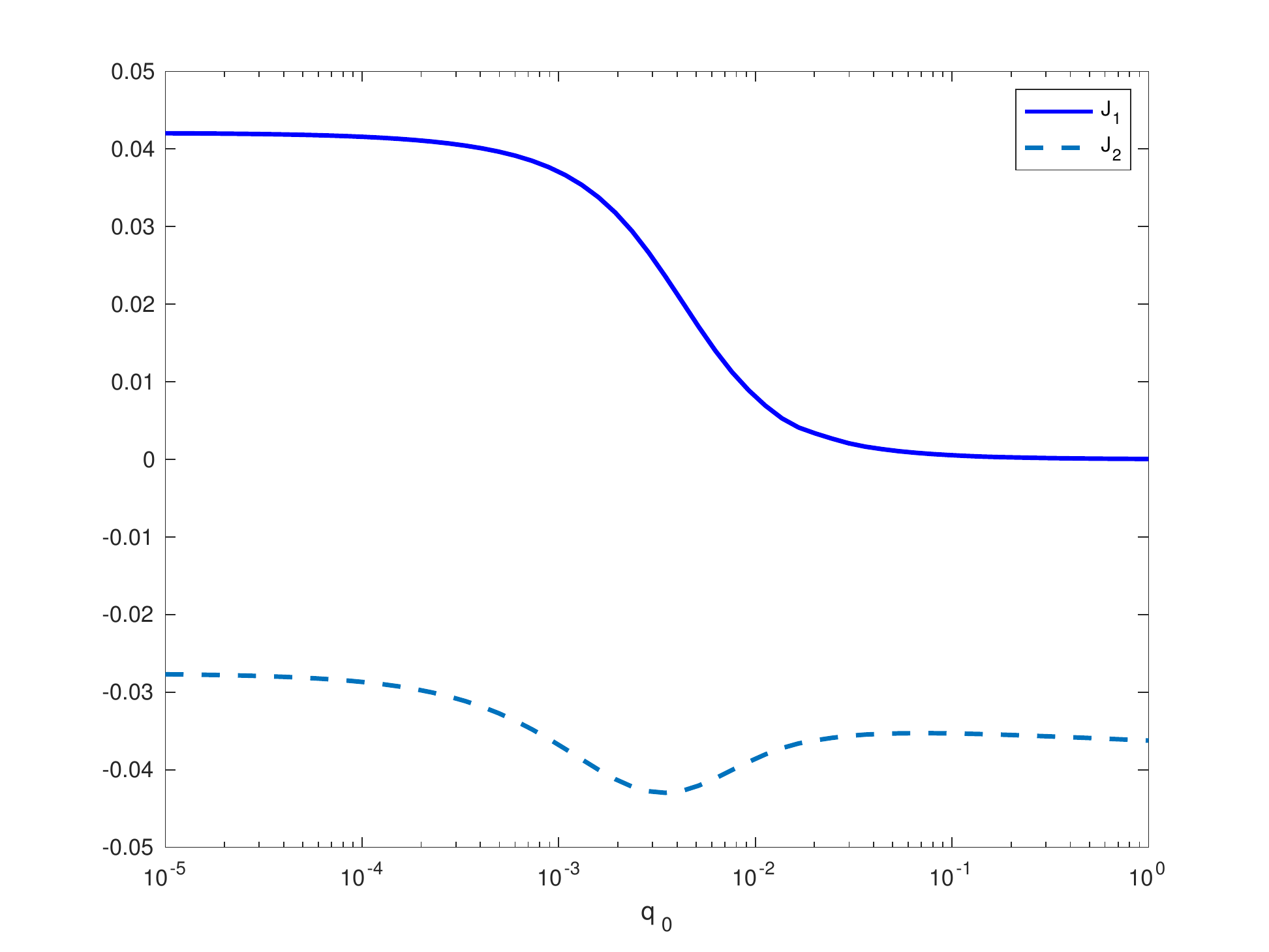} \label{fig:JV10}}
\subfigure[$V = -60$]{\includegraphics[width = 0.47\textwidth]{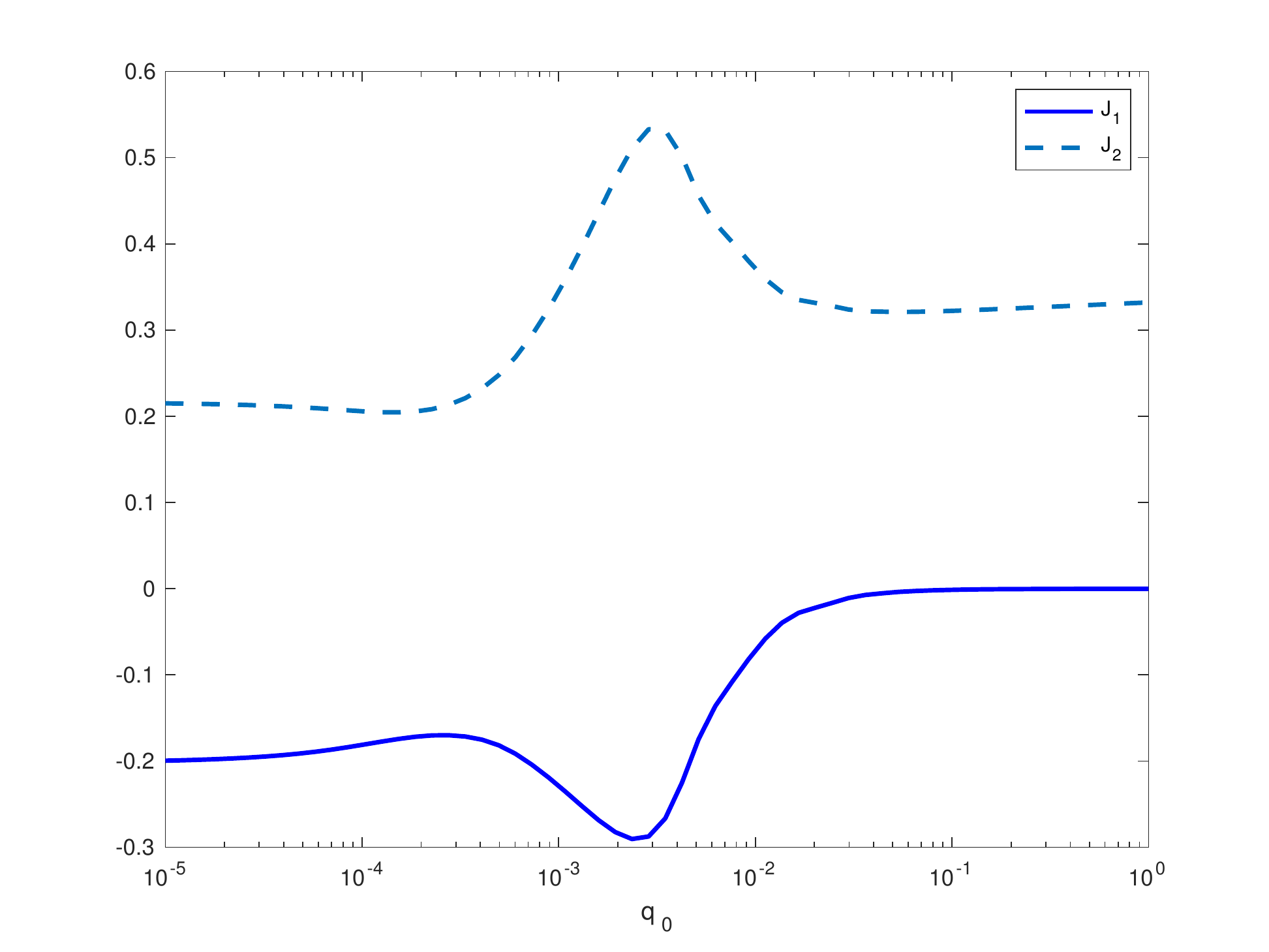} \label{fig:JV-60}}\qquad
\subfigure[$V = -110$]{\includegraphics[width = 0.47\textwidth]{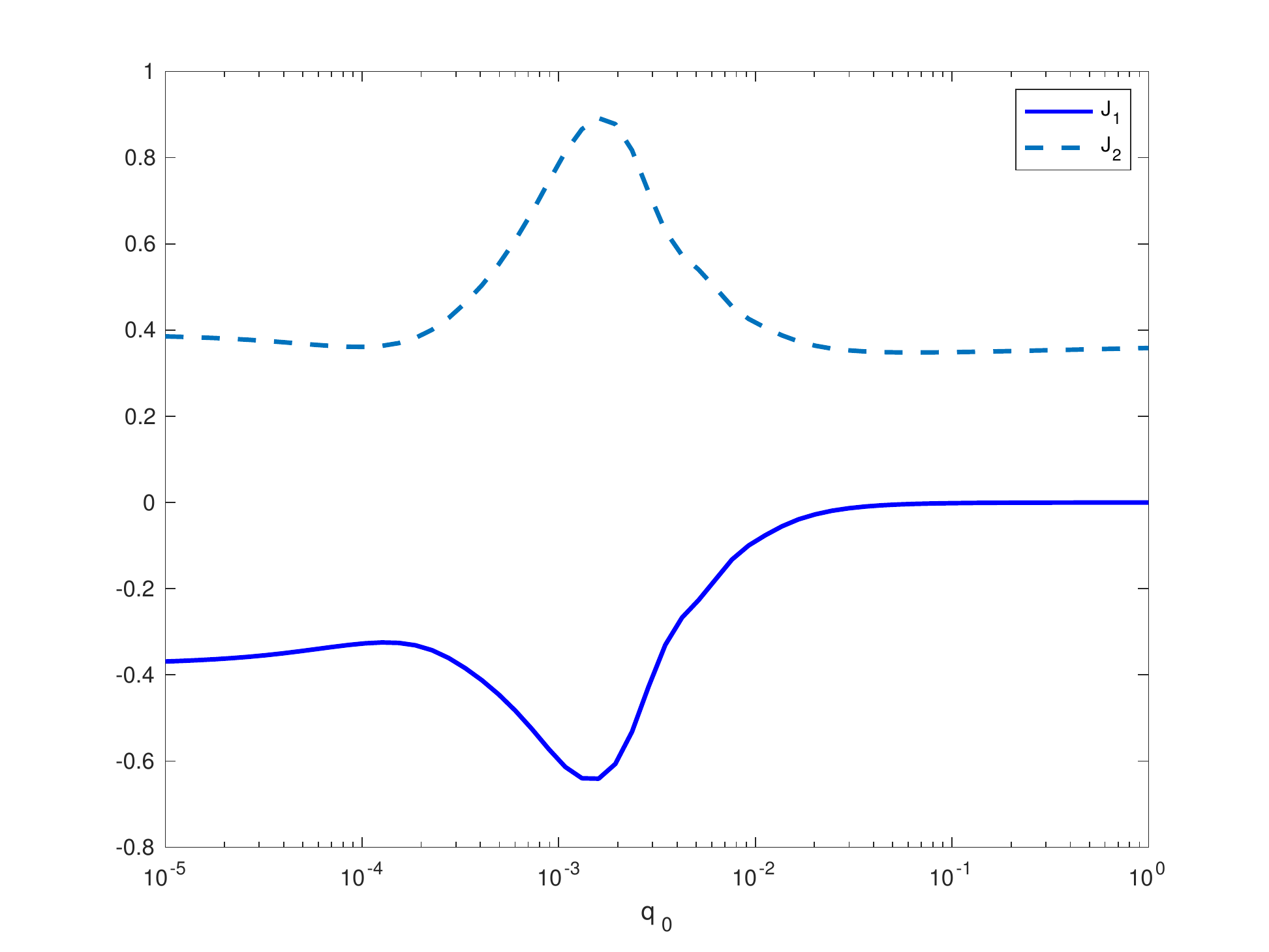} \label{fig:JV-110}}
\caption{$J_1$ and $J_2$ are plotted as functions of $q_0 = 2Q_0$, with boundary values of the system (\ref{dimensionlessPNP}) chosen as $L = 0.008$, $R = 0.001$, and $V=50$, $10$, $-60$, or $-110$.\label{fig:J-Q}}
\end{figure}

In Fig.~\ref{fig:J-Q}, $J_1$ and $J_2$ are plotted as functions of $Q_0$ associated with same $V$ values as in Fig.~\ref{fig:lambda-Q}.

\subsection{Dependence of $\lambda_1$ and $\lambda_2$ on $V$ for fixed $Q_0$}
\label{sec:Valternate}
\begin{figure}[htbp]
\subfigure[$2 Q_0 = 0.0001$]{\includegraphics[width = 0.47\textwidth]{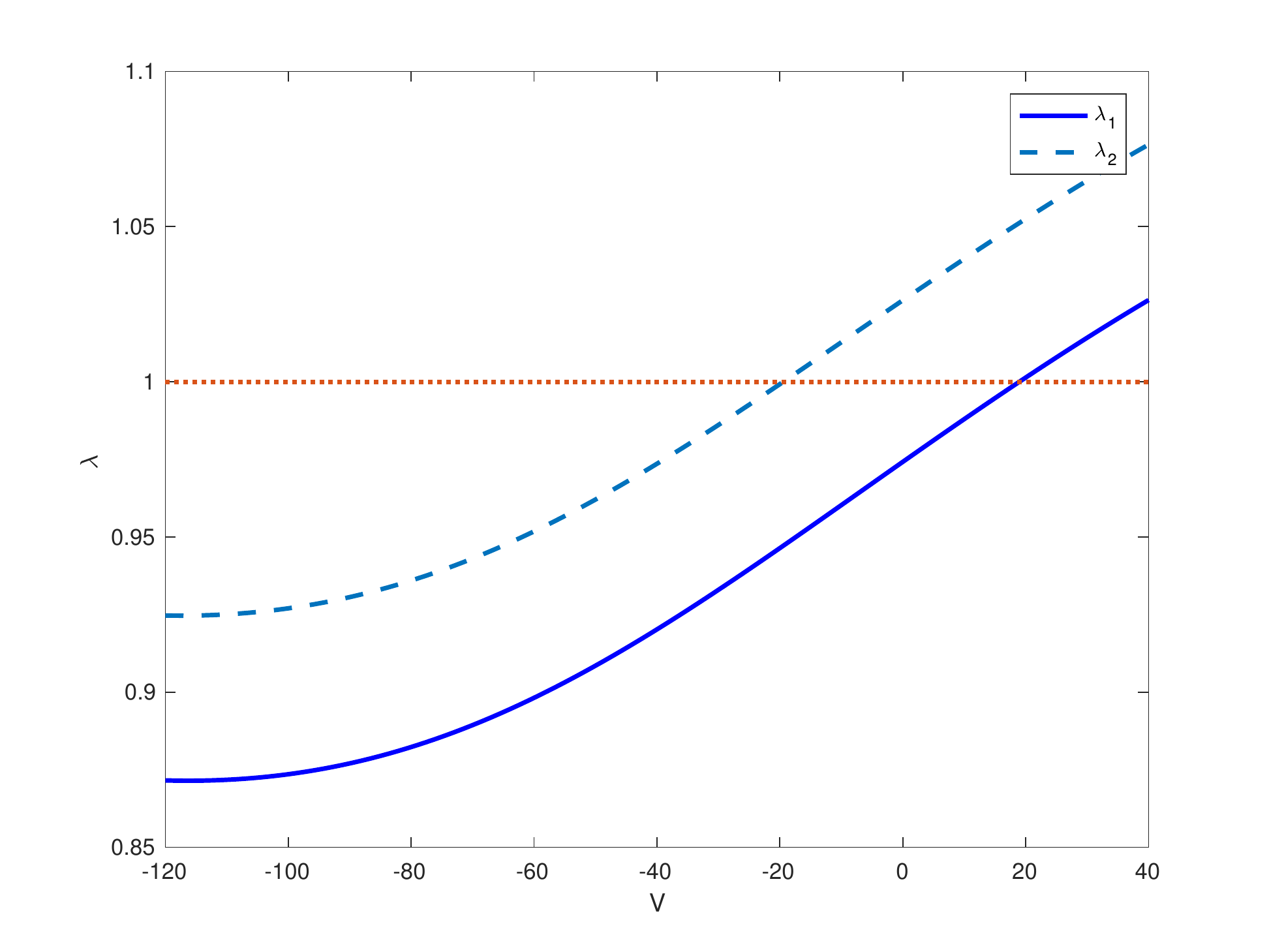} \label{fig:lambda-V-Q1e_4}}\qquad
\subfigure[$2 Q_0 = 0.00037$]{\includegraphics[width = 0.47\textwidth]{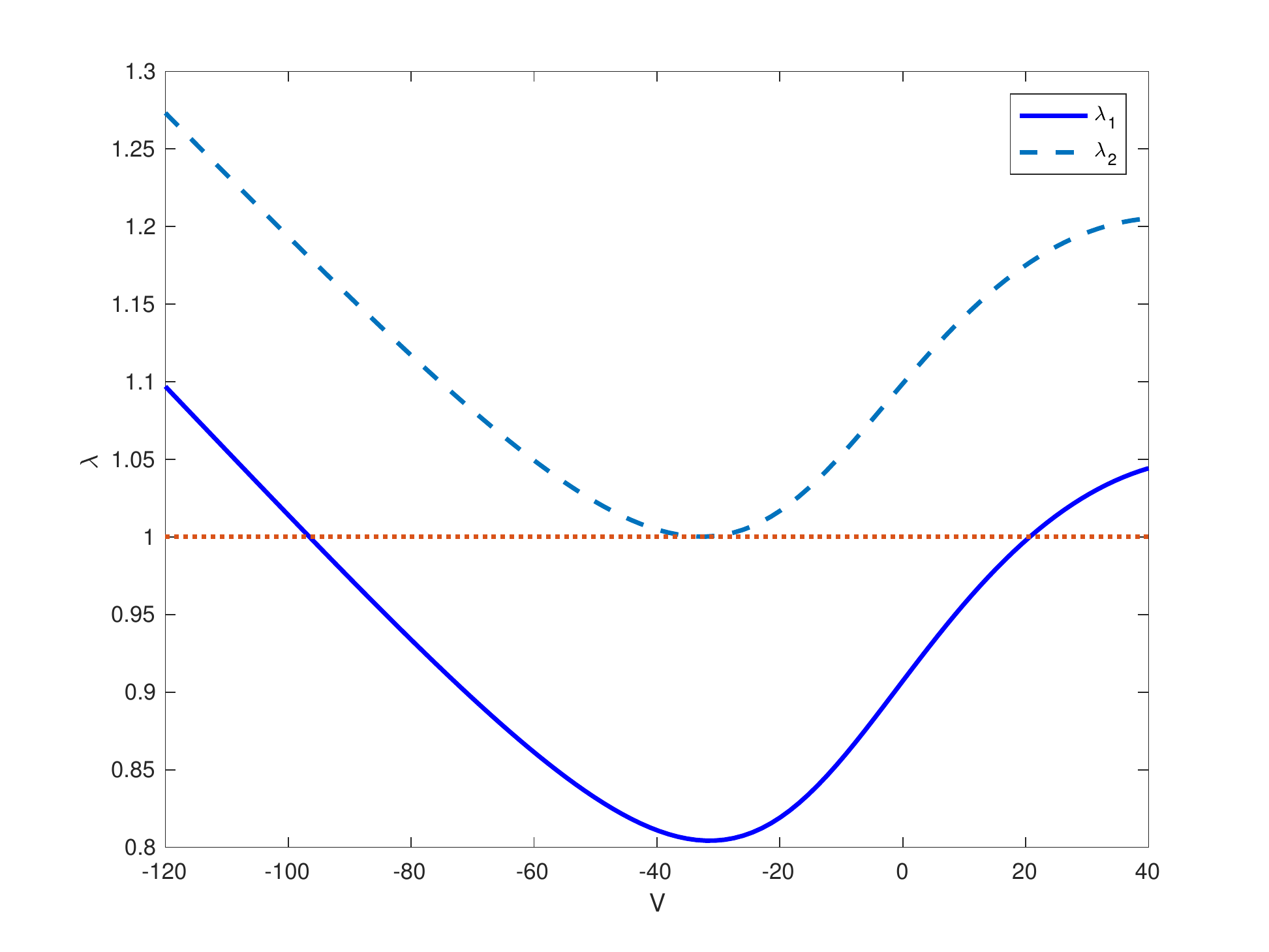} \label{fig:lambda-V-Q37e_5}}
\subfigure[$2 Q_0 = 0.00062$]{\includegraphics[width = 0.47\textwidth]{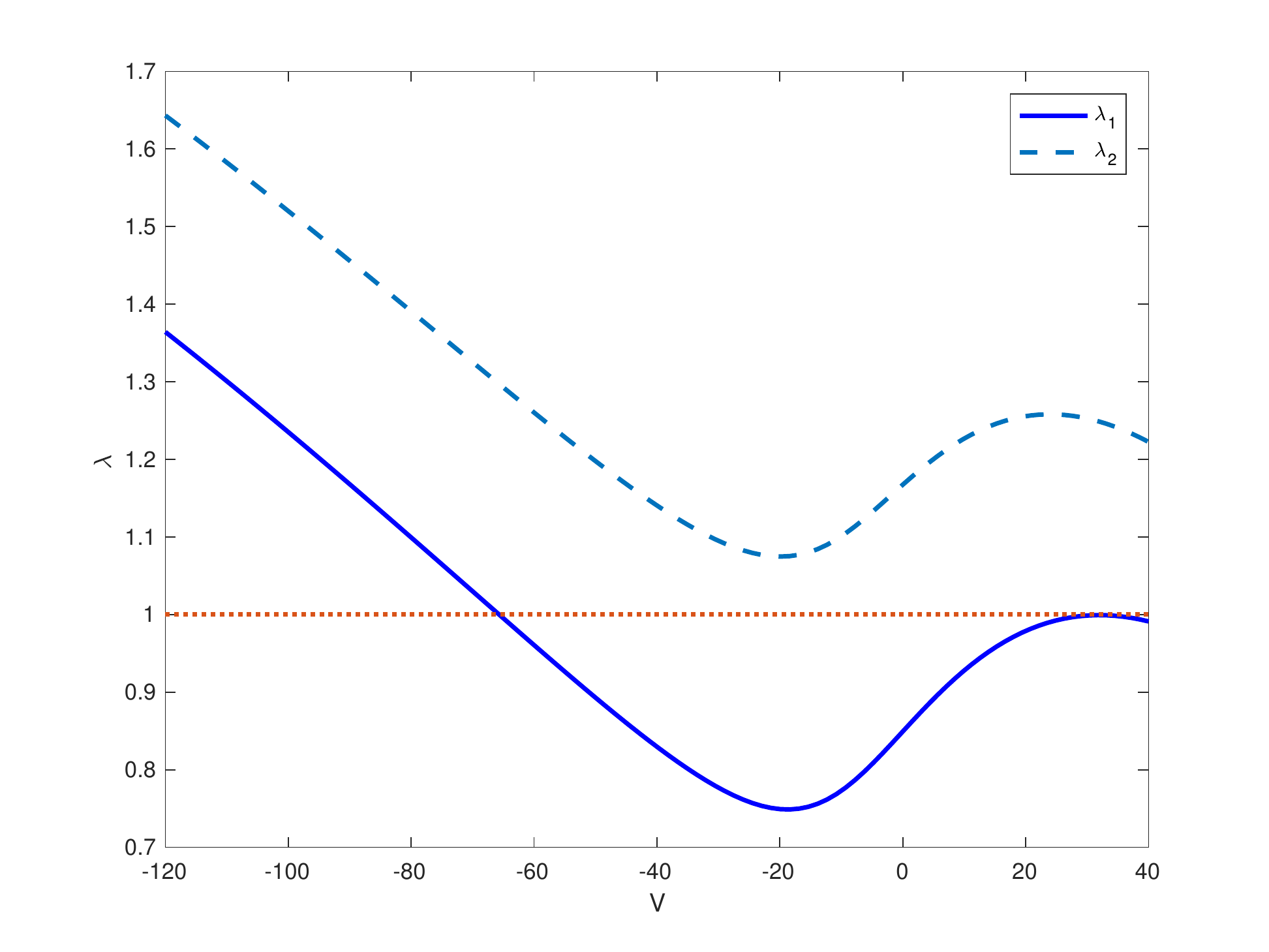} \label{fig:lambda-V-Q62e_5}}\qquad
\subfigure[$2 Q_0 = 0.04$]{\includegraphics[width = 0.47\textwidth]{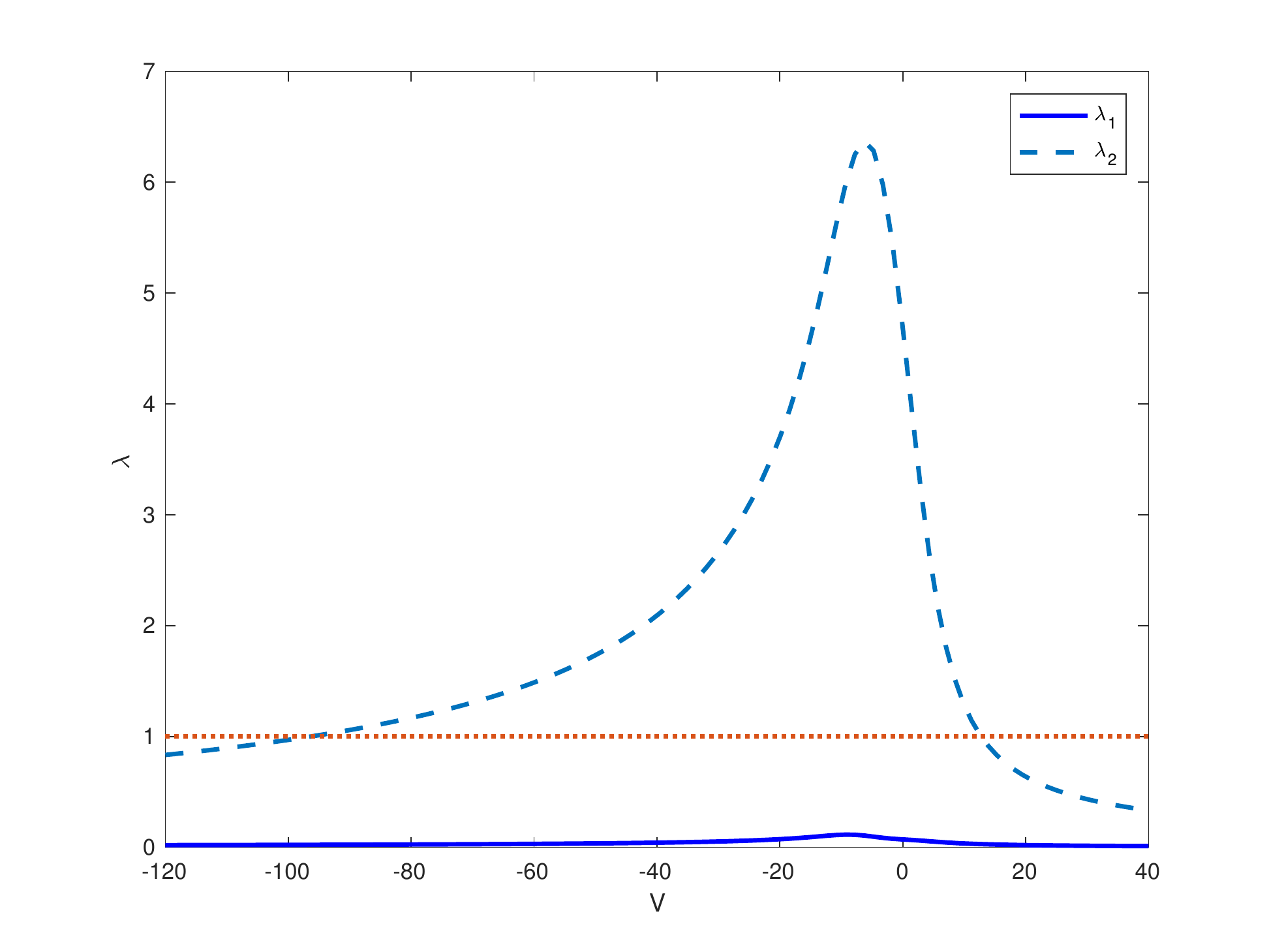} \label{fig:lambda-V-Q4e_2}}
\caption{$\lambda_1$ and $\lambda_2$ are plotted as functions of $V$, with boundary values of the system (\ref{dimensionlessPNP}) chosen as $L = 0.008$, $R = 0.001$, and different values of $Q_0$.\label{fig:lambda-V}}
\end{figure}

We now examine the dependence of $\lambda_1$ and $\lambda_2$ on $V$ for several fixed values of $Q_0$. In Fig.~\ref{fig:lambda-V}, they are plotted as functions of $V \in (-110, 40)$ for $q_0 \equiv 2 Q_0 = 0.0001$, $0.00037$, $0.00062$, and $0.04$. 


{\bf Monotonicity in $V$ for small $Q_0$.} From Fig.~\ref{fig:lambda-V-Q1e_4} one observes that for very small $2 Q_0=0.0001$, both $\lambda_1$ and $\lambda_2$ are monotone in $V$. This is consistent with the theoretical prediction made in \cite{JLZ2015} and the intuition that the flux ratios are dominated by the effects of $V$ when $Q_0$ is small. The $V$ values at the intersections of the curves with $\lambda = 1$ (which will be referred to as $V$-intercepts) are close to the theoretical values obtained in \cite{JLZ2015}, i.e.,
$V_1^0 = 18.97$, and $V_2^0 = -18.97$ (cf. \S\ref{sec:Qalternate}).

{\bf Saddle-node bifurcations of $\lambda_k=1$.} As $Q_0$ increases to $2 Q_0 = 0.00037$ and $0.00062$,  as shown in Fig.~\ref{fig:lambda-V-Q37e_5} and Fig.~\ref{fig:lambda-V-Q62e_5}, respectively, $\lambda_1$ and $\lambda_2$ become non-monotone in $V$ and can pass the value $1$ multiple times. This behavior is not predicted by the analysis in \cite{JLZ2015}. As a matter of fact, the analysis of \cite{JLZ2015} is valid only for very small $Q_0$ and thus cannot be used to predict the behaviors of $\lambda_1$ and $\lambda_2$ for $Q_0$ not very small.
As will be shown in \S\ref{sec:CompleteQV}, what value of $Q_0$ is called small is relative, depending on $L$ and $R$.
 Moreover, $2 Q_0=0.00037$ for Fig.~\ref{fig:lambda-V-Q37e_5} is near a saddle-node bifurcation value for the relation $\lambda_2(Q_0, V)=1$ in the sense that it has two roots in $V$ for $2 Q_0$ less than and close to $0.00037$ and no roots for $2 Q_0$ greater than and close to $0.00037$.  Similarly, $2 Q_0=0.00062$ for Fig.~\ref{fig:lambda-V-Q62e_5} is near a saddle-node bifurcation value for the relation $\lambda_1(Q_0,V)=1$.   

{\bf A selectivity regime.} It is also interesting to see that in Fig.~\ref{fig:lambda-V-Q4e_2} ($2 Q_0=0.04$), $\lambda_1$ is close to zero, indicating that large positive permanent charges inhibit the flow of cation.  On the other hand,  for $V \in (V_1^\infty, V_2^\infty) \approx (-90,15)$, $\lambda_2>1$ so the anion flux is enhanced, and for $V\in [-120,-90) \cup (15,40]$, $\lambda_2<1$ so the anion flux is reduced (even if $Q_0$ is positive and relatively large). This observation verifies the analysis described in \S\ref{SEC:recentresults} for relatively large $Q_0$. It is also interesting to compare it with the behavior in  Fig.~\ref{fig:lambda-V-Q62e_5} for a smaller $2 Q_0=0.00062$ where $\lambda_2>1$ for $V\in [-120,40]$.

Moreover,   around $V=-5$, say $V\in [-15,0]$ in Fig.~\ref{fig:lambda-V-Q4e_2}, $\lambda_2$ is much greater than $1$, which shows  strong selectivity for anion. 

\begin{figure}[htbp]
\subfigure[$2 Q_0 = 0.0001$]{\includegraphics[width = 0.47\textwidth]{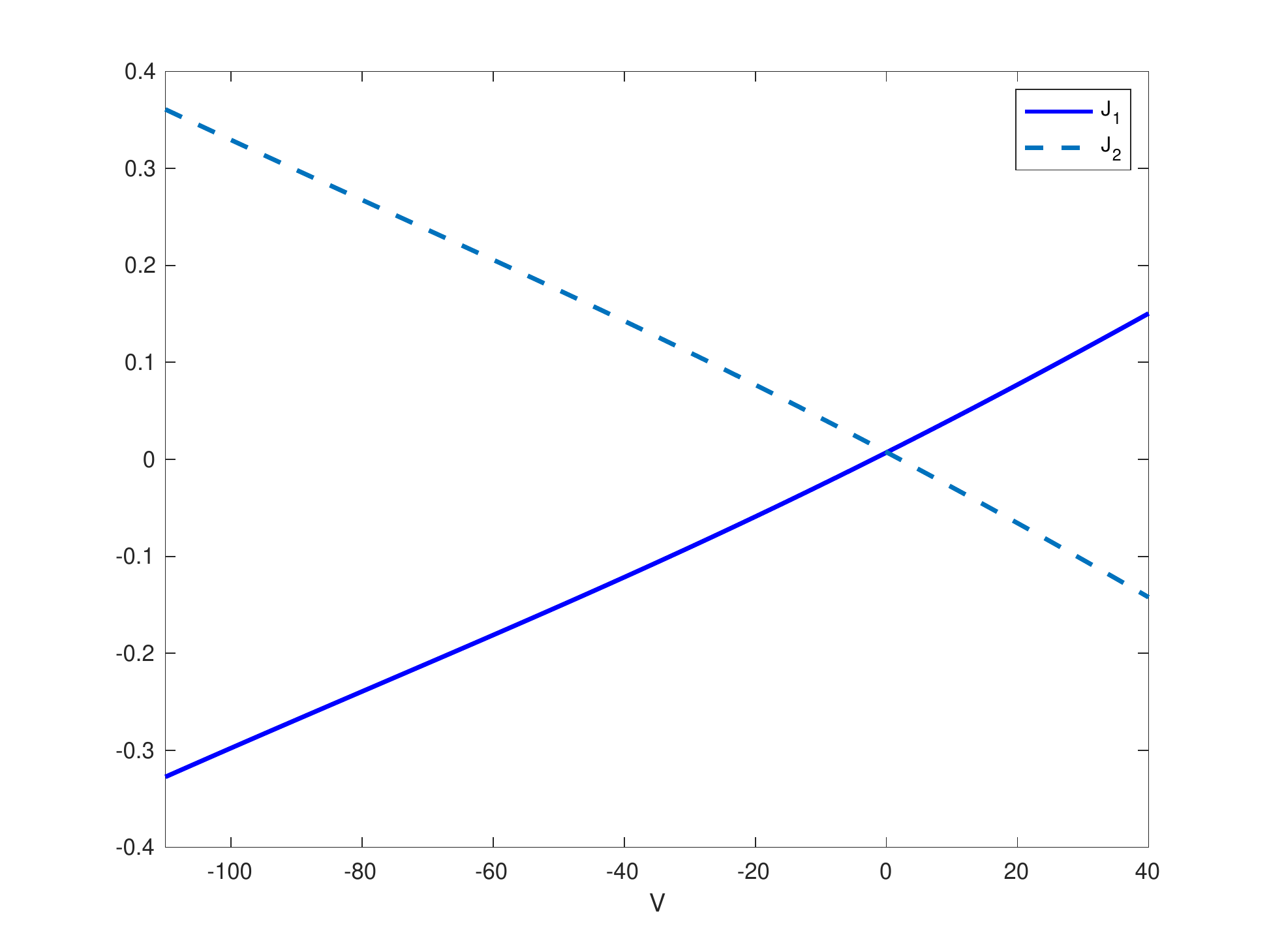} \label{fig:JQ1e_4}}\qquad
\subfigure[$2 Q_0 = 0.00037$]{\includegraphics[width = 0.47\textwidth]{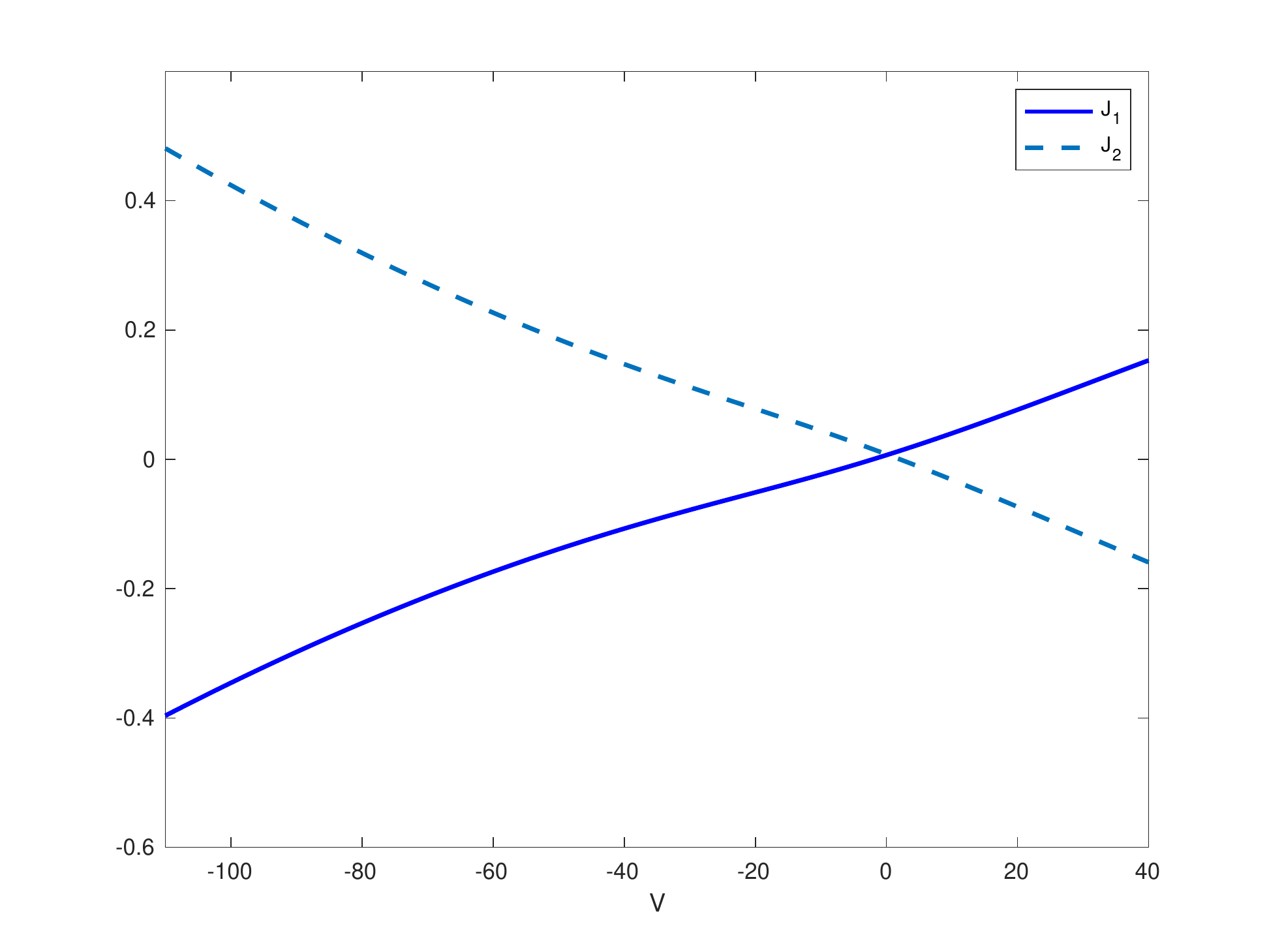} \label{fig:JQ37e_5}}
\subfigure[$2 Q_0 = 0.00062$]{\includegraphics[width = 0.47\textwidth]{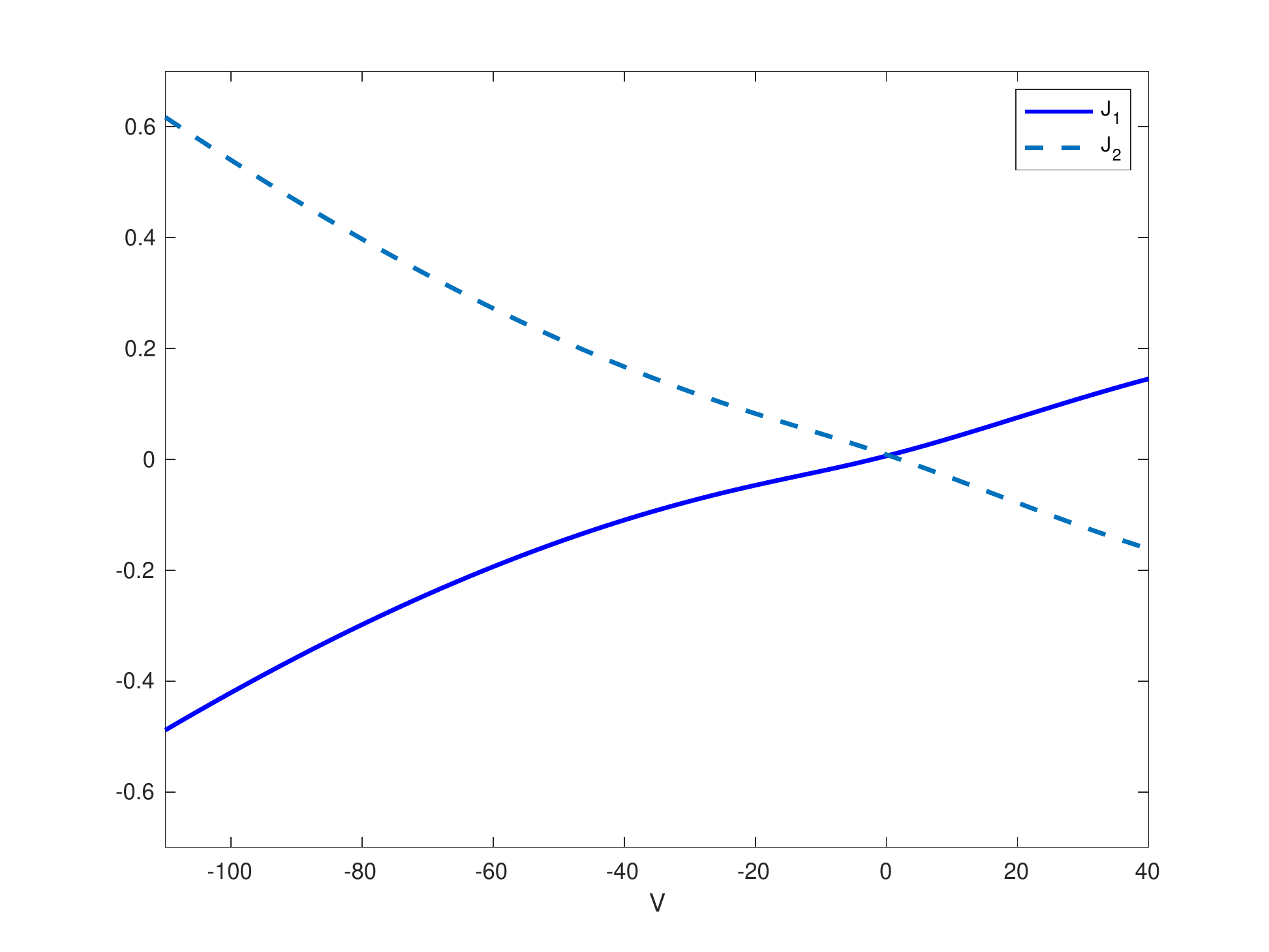} \label{fig:JQ62e_5}}\qquad
\subfigure[$2 Q_0 = 0.04$]{\includegraphics[width = 0.47\textwidth]{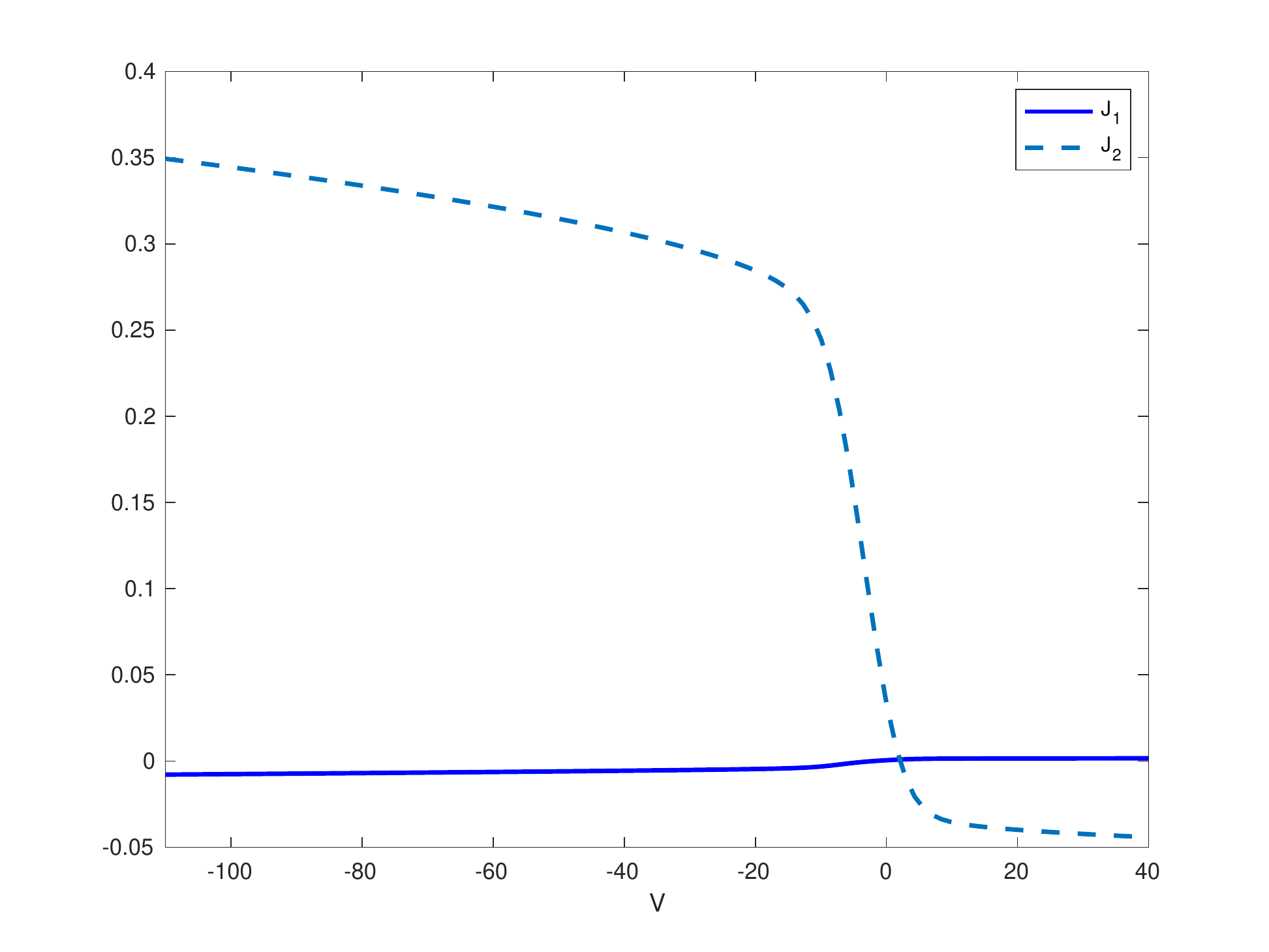} \label{fig:JQ4e_2}}
\caption{$J_1$ and $J_2$ are plotted as functions of $V$, with boundary values of the system (\ref{dimensionlessPNP}) chosen as $L = 0.008$, $R = 0.001$, and different values of $Q_0$.\label{fig:J-V}}
\end{figure}

 In order to investigate the $J$-$V$ relation, in Fig.~\ref{fig:J-V} we present $J_1$ and $J_2$ as functions of $V$, with the same values of $Q_0$ as in Fig.~\ref{fig:lambda-V}. One can see that for small $Q_0$ values, such as $2 Q_0 = 0.0001$ as shown in Fig.~\ref{fig:JQ1e_4}, $J_k$ is very close to a linear function of $V$. As we increase the value of $Q_0$, the $J$-$V$ maintains less linearity. This observation is consistent with our knowledge of the $I$-$V$ relation associated with different permanent charge densities. Moreover, one can observe that to some extent, for large negative voltage values, such as $V = -110$, increasing the positive permanent charge density can even enhance the positive flux, when $2 Q_0$ increases from $ 0.0001$ to $0.00062$ as shown in Fig.~\ref{fig:JV-60} and \ref{fig:JV-110} or Fig.~\ref{fig:JQ1e_4},~\ref{fig:JQ37e_5}, and ~\ref{fig:JQ62e_5}. However, if we increase this value too much, for example, to $2 Q_0 = 0.04$ as in Fig.~\ref{fig:JQ4e_2}, not only the positive flux, but also the negative flux, are prohibited for $V = -110$. The abrupt change of $J_2$ from $V = 0$ to $V = -20$ is consistent to the selectivity regime we have observed in Fig.~\ref{fig:lambda-V-Q4e_2}.

Further studies on these rich phenomena and how the interaction of $V$ with $Q_0$ influences $\lambda_1$ and $\lambda_2$ would be interesting.

\subsection{A complete bifurcation diagram of $\lambda_k(Q_0,V)$}
\label{sec:CompleteQV}

\begin{figure}
\centering
\includegraphics[scale = 0.6]{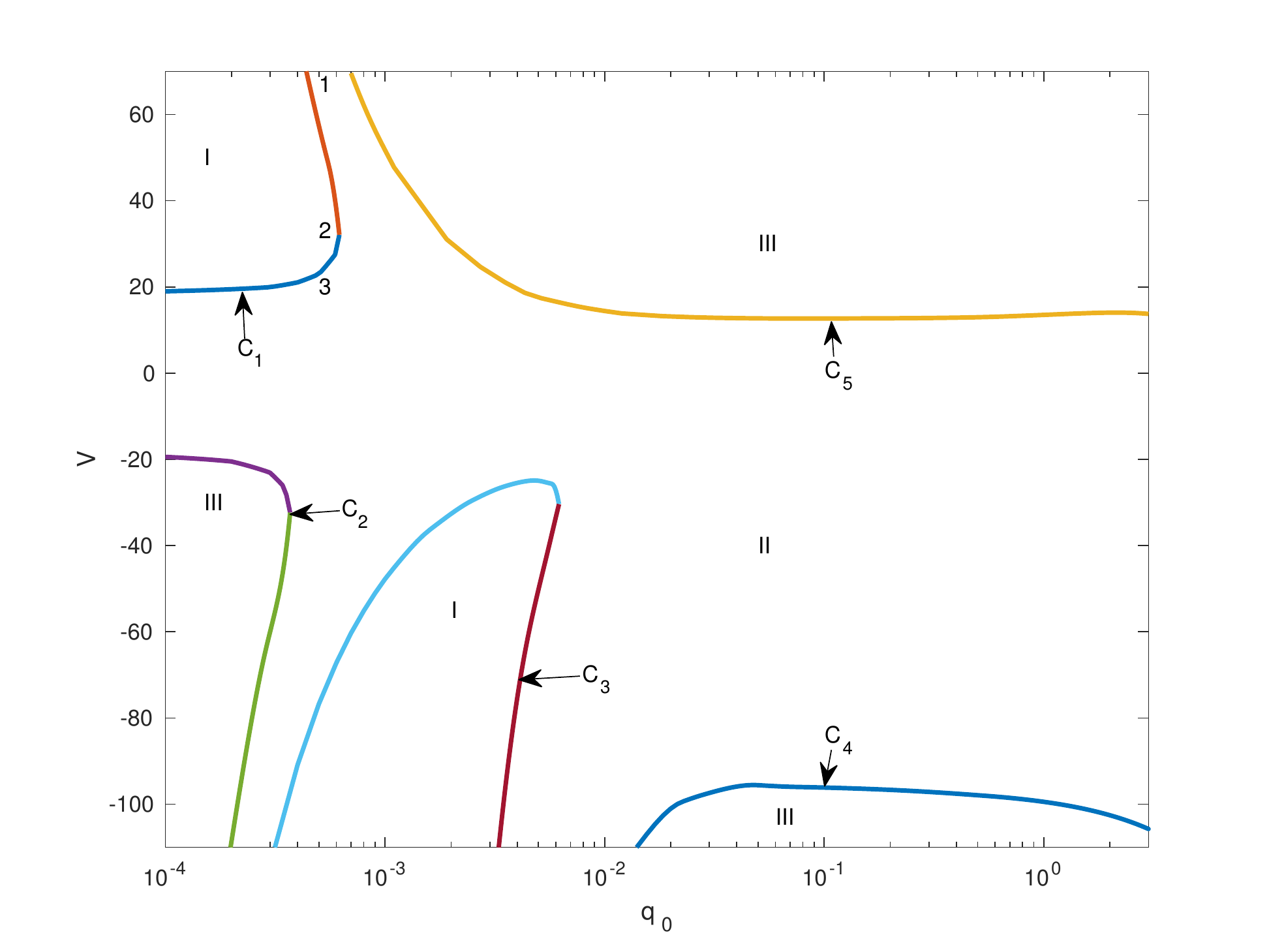}
\caption{A complete bifurcation diagram for the case with the boundary conditions $L = 0.008$ and $R = 0.001$. Region I: $1 < \lambda_1 < \lambda_2$; Region II: $\lambda_1 < 1 < \lambda_2$; Region III: $ \lambda_1 < \lambda_2 < 1$. The curves $C_1$ and $C_3$ are for $\lambda_1=1$ while $C_2$, $C_4$, and $C_5$ are for $\lambda_2=1$. Internal dynamics at Point 1, 2, and 3 are shown in Fig.~\ref{fig:Dynamics}. Here, $q_0 = 2 Q_0$.}
\label{L8R1complete}
\end{figure}

We now are ready to study a complete bifurcation diagram of $\lambda_1(Q_0, V)$ and $\lambda_2(Q_0, V)$, with fixed boundary concentrations. For a couple of fixed $(L, R)$, we work on the domain $(0,3) \times (-110,70)$ for $(q_0, V)$, where $q_0\equiv 2 Q_0 = 3$ is large enough. Our goal is to determine the partition of the domain into regions with $ 1 < \lambda_1 < \lambda_2$ (Region I), $\lambda_1 < 1 < \lambda_2$ (Region II), and $\lambda_1 < \lambda_2 < 1$ (Region III), respectively. Fig.~\ref{L8R1complete} contains the numerical experiment for $(L, R) = (0.008, 0.001)$, where the boundary curves $C_1$, ..., $C_5$ for $\lambda_1 = 1$ or $\lambda_2 = 1$ are plotted.

We first remark that Fig.~\ref{L8R1complete} is consistent with the analytical analysis described in \S\ref{SEC:recentresults}. In particular, for small $Q_0$, the situation switches from $1 < \lambda_1 < \lambda_2$ to $ \lambda_1 < 1 < \lambda_2$ as $V$ decreases and passes $V_1^0 = 18.97$ and from $\lambda_1 < 1 < \lambda_2$ to $ \lambda_1 < \lambda_2 < 1$ as $V$ further decreases and passes $V_2^0 = -18.97$. Moreover, the numerical results show that $V \to V_1^0$ along the curve $C_1$ ($\lambda_1 = 1$) as $Q_0 \to 0$. Similarly, $V \to V_2^0$ along $C_2$ ($\lambda_2 = 1$) as $Q_0 \to 0$. Furthermore, we recall that, by definition, $\lambda_1 = \lambda_2 = 1$ for $Q_0 = 0$ and all $V$. $C_1$ along with the $V$-axis makes part of the $\lambda_1 = 1$ contour while $C_2$ along with the $V$-axis makes part of the $\lambda_2 = 1$ contour. From Fig.~\ref{L8R1complete} we can also see that, for relatively large $Q_0$, the domain contains basically two regions, Region II ($\lambda_1 < 1 < \lambda_2$) and Region III ($\lambda_1 < \lambda_2 < 1$), which is consistent with the analytical results described in \S\ref{SEC:recentresults}.
It is clear from the figure how different the behaviors of $\lambda_1$ and $\lambda_2$ are for small and large
$Q_0$ and how these are connected in the $Q_0$-$V$ plane. It is unlikely that this type of result can be obtained using the existing analytical techniques.

As mentioned in \S\ref{SEC:recentresults}, there exist critical potential values that $V$ converges to along the parametric boundaries $\lambda_1 = 1$ and $\lambda_2 = 1$ for small and large $Q_0$. These values provide predictable boundaries of the range of the permanent charge that enhances $J_1$ or $J_2$. This property has been demonstrated in Fig.~\ref{L8R1complete}. We remark that the critical potential values that we have obtained in the computation
are slightly different from the analytical values due to the fact that the regularization has been used in (\ref{Qregular}) and (\ref{Aregular}) for numerical computation whereas the analytical work in \cite{JLZ2015} and \cite{ZEL2019} has been based on non-smooth permanent charges and channel shape as in (\ref{shape}) and (\ref{permanentcharge}). Nevertheless, Fig.~\ref{L8R1complete} still demonstrates the limiting properties, and is a reliable reference for ion species selectivities with large permanent charge densities. 

It is impressive how complicated  the domain is partitioned into the regions $1 < \lambda_1 < \lambda_2$, $\lambda_1 < 1 < \lambda_2$, and $\lambda_1 < \lambda_2 < 1$. Multiple bifurcations on the curves $\lambda_1 = 1$ and $\lambda_2 = 1$ can be observed, which are consistent with the results in \S\ref{sec:Valternate}. In fact, if we treat $\lambda_k$ as functions of the variables $Q_0$ and $V$, and the $Q_0$-$V$ space in Fig.~\ref{L8R1complete} as a domain of the functions $\lambda_k(Q_0,V)$, Figs.~\ref{fig:lambda-Q} and \ref{fig:lambda-V} can be treated as curves obtained by cutting through the surfaces $\lambda_k(Q_0,V)$ along $Q_0$ axis direction or $V$ axis direction. For example, Fig.~\ref{fig:lambda-V-Q37e_5} essentially shows the intersecting surfaces of $\lambda_k(Q_0,V)$ and the plane $2 Q_0 = 0.00037$. Notice that the vertical line $2 Q_0 = 0.00037$ is approximately tangent to the curve $C_2$ (on which $\lambda_2 = 1$) in Fig.~\ref{L8R1complete} at the saddle-node bifurcation point $(2 Q_0, V) \approx (0.00037,-30)$.
The figure shows that there are other saddle-node bifurcation points.
Further studies, both biological and mathematical, on the cause of the bifurcation are highly demanded.

\begin{figure}
\centering
\includegraphics[scale = 0.6]{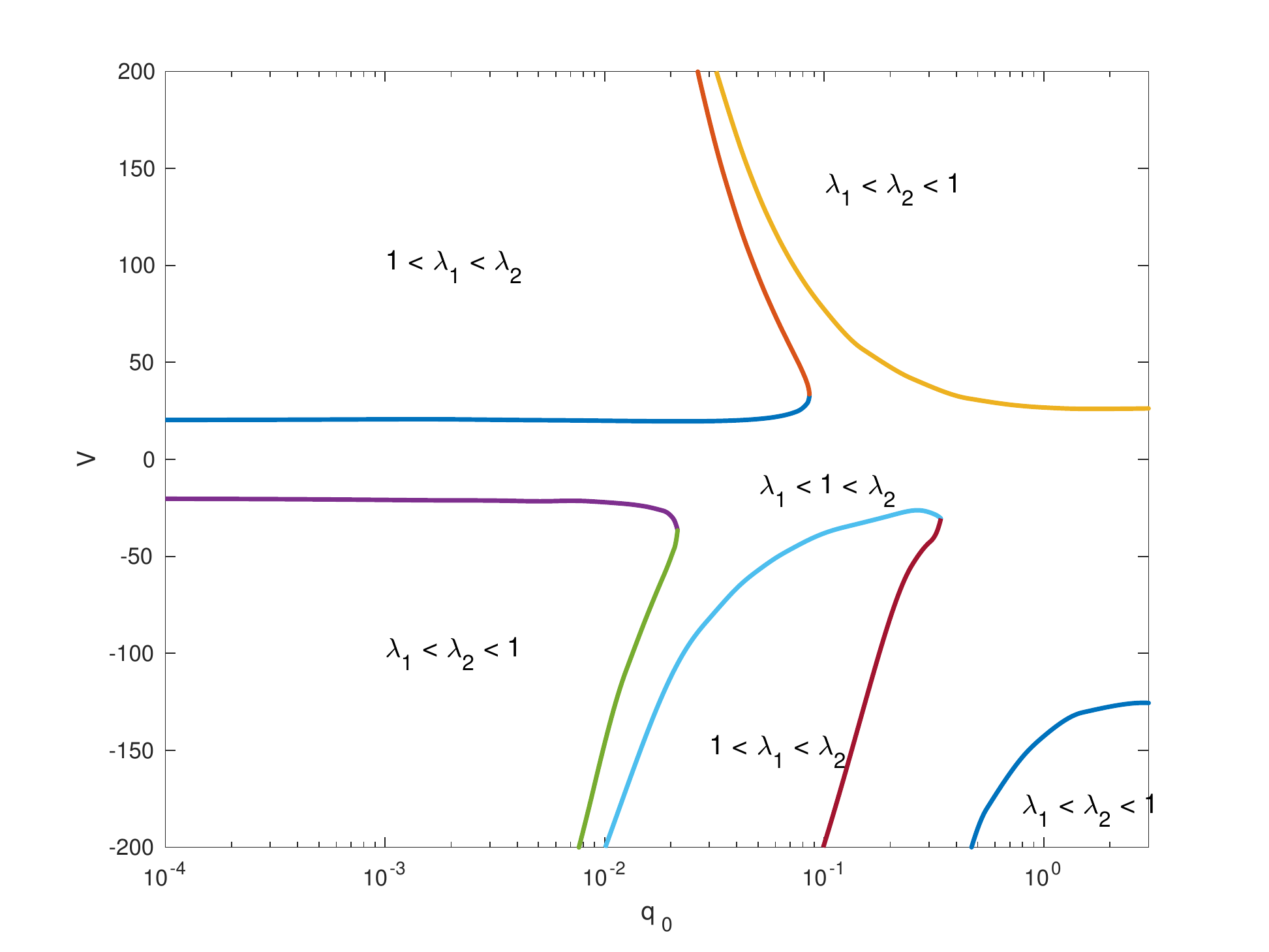}
\caption{A complete bifurcation diagram for the case with the boundary conditions $L = 0.5$ and $R = 0.1$.
Here, $q_0 = 2 Q_0$.}
\label{L5R1complete}
\end{figure}

We emphasize that we have not particularly chosen boundary values to obtain such bifurcations. Similar properties of Fig.~\ref{L8R1complete} have been observed with other boundary conditions. For example, Fig.~\ref{L5R1complete} shows the complete $Q_0$-$V$ diagram obtained with $L = 0.5$ and $R = 0.1$. Despite the difference in the quantities, one can observe similar bifurcations and limiting properties as in Fig.~\ref{L8R1complete}.

\subsection{Internal dynamics and $J$-$V$ relation}
\label{internaldynamics}

A reasonable approach to understand the cause of the bifurcations in Fig.~\ref{L8R1complete} is to investigate the channel dynamics when the pair $(Q_0,V)$ passes the bifurcation. For this purpose, we present the numerical results of the channel dynamics at the Points 1, 2, and 3 in Fig.~\ref{L8R1complete}. It should be pointed out that it is unclear if we can make an observation from the flux dynamics that could be related to the cause of the bifurcation or the non-monotone nature of $\lambda_k$. We look forward to further studies on this topic.

\begin{figure}[htbp]
\subfigure[$J_1$, $J_2$, and $I$ at Point $1$]{\includegraphics[width = 5cm]{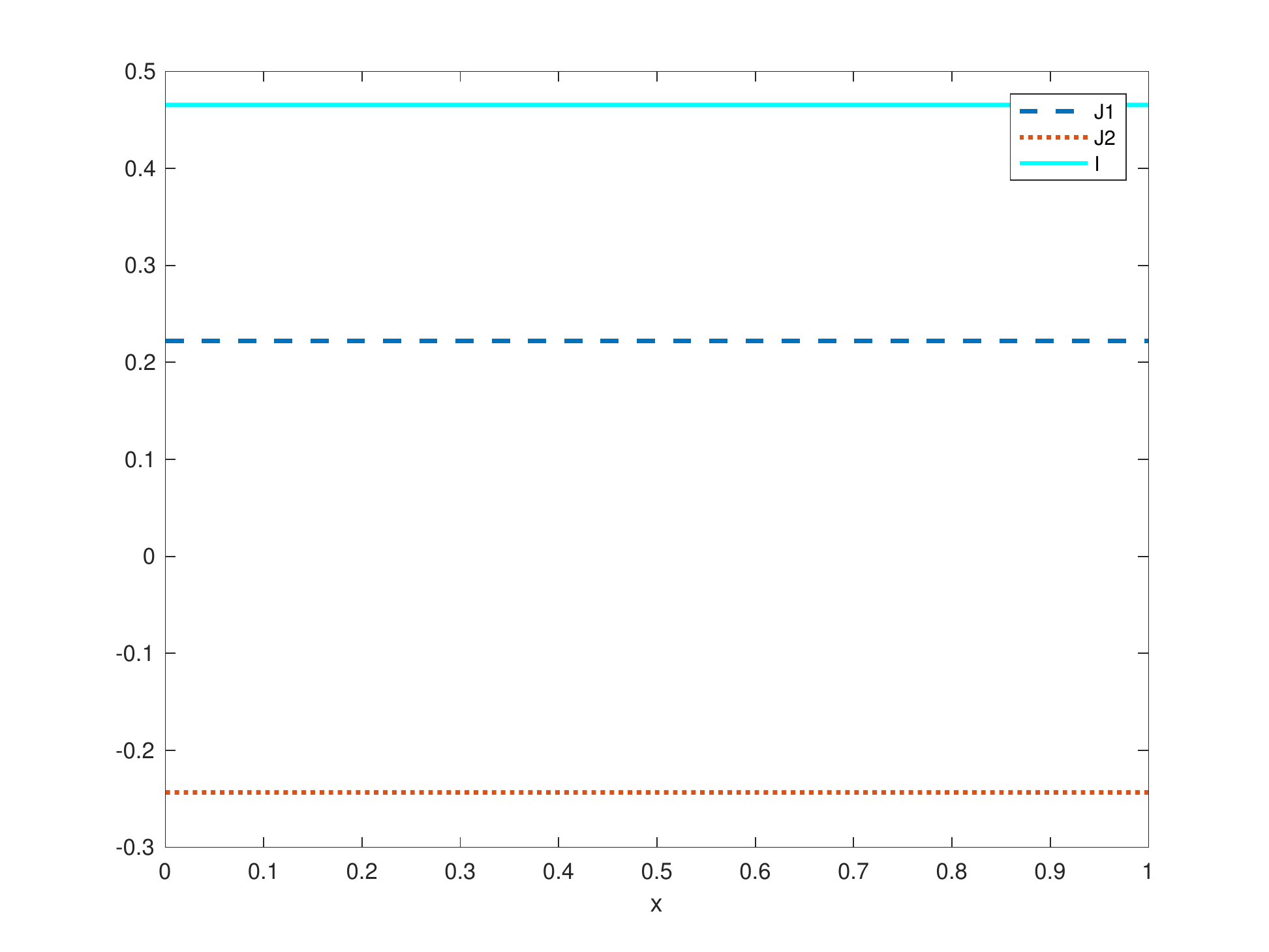} \label{fig:dynamics1}}
\subfigure[$J_1$, $J_2$, and $I$ at Point $2$]{\includegraphics[width = 5cm]{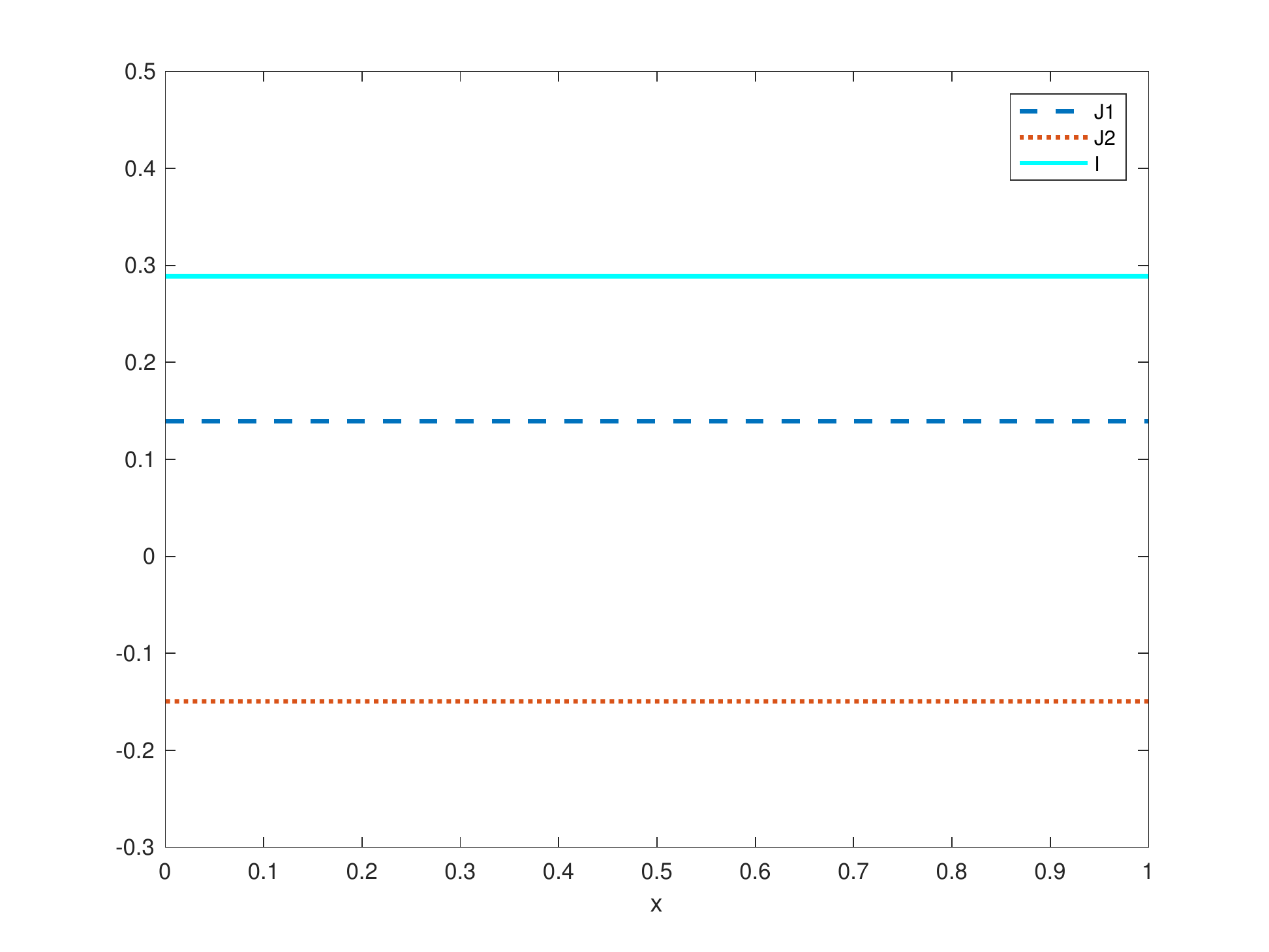} \label{fig:dynamics2}}
\subfigure[$J_1$, $J_2$, and $I$ at Point $3$]{\includegraphics[width = 5cm]{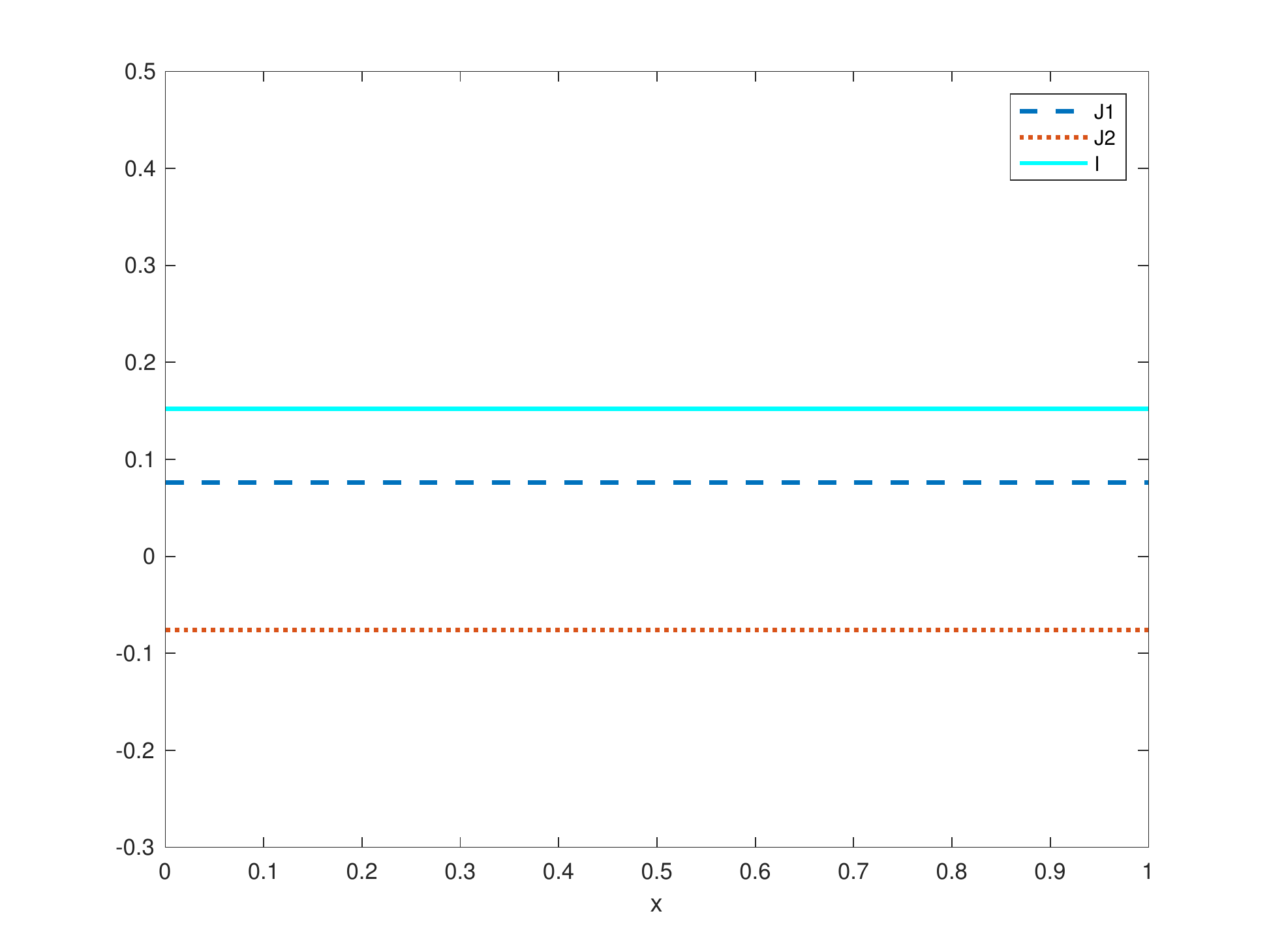} \label{fig:dynamics3}}

\subfigure[$\mu_1$ and $\mu_2$ at Point $1$]{\includegraphics[width = 5cm]{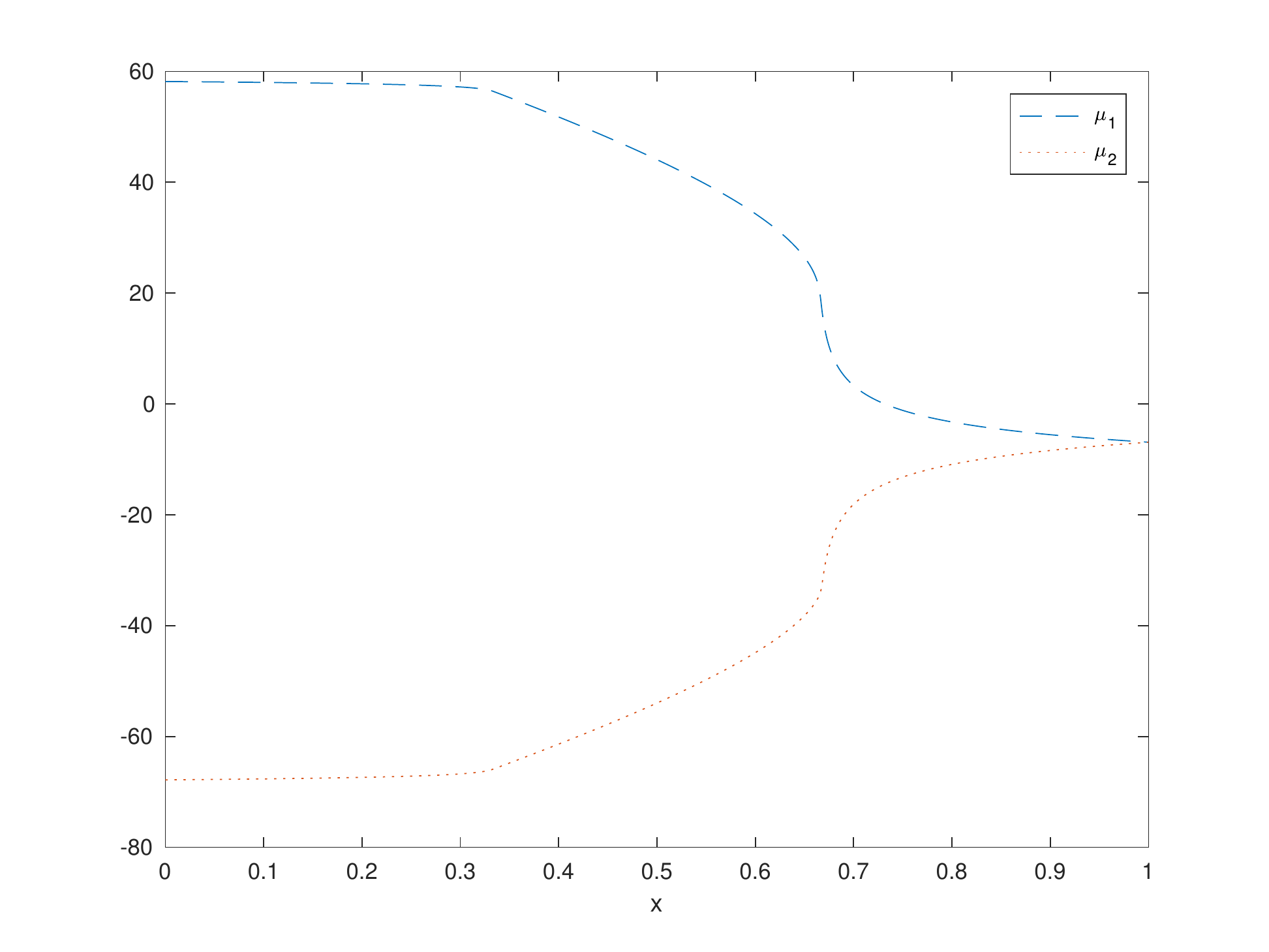} \label{fig:dynamics4}}
\subfigure[$\mu_1$ and $\mu_2$ at Point $2$]{\includegraphics[width = 5cm]{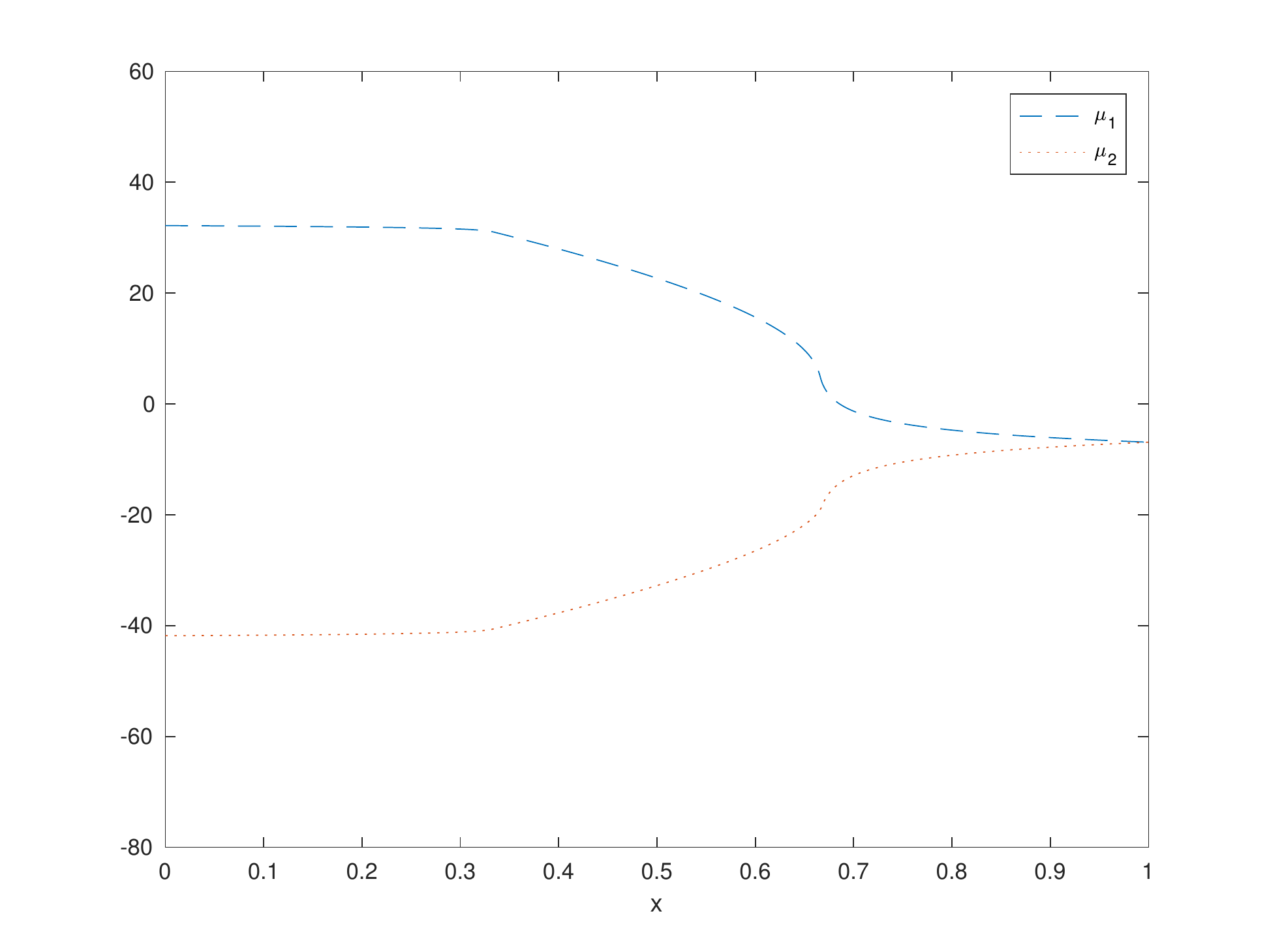} \label{fig:dynamics5}}
\subfigure[$\mu_1$ and $\mu_2$ at Point $3$]{\includegraphics[width = 5cm]{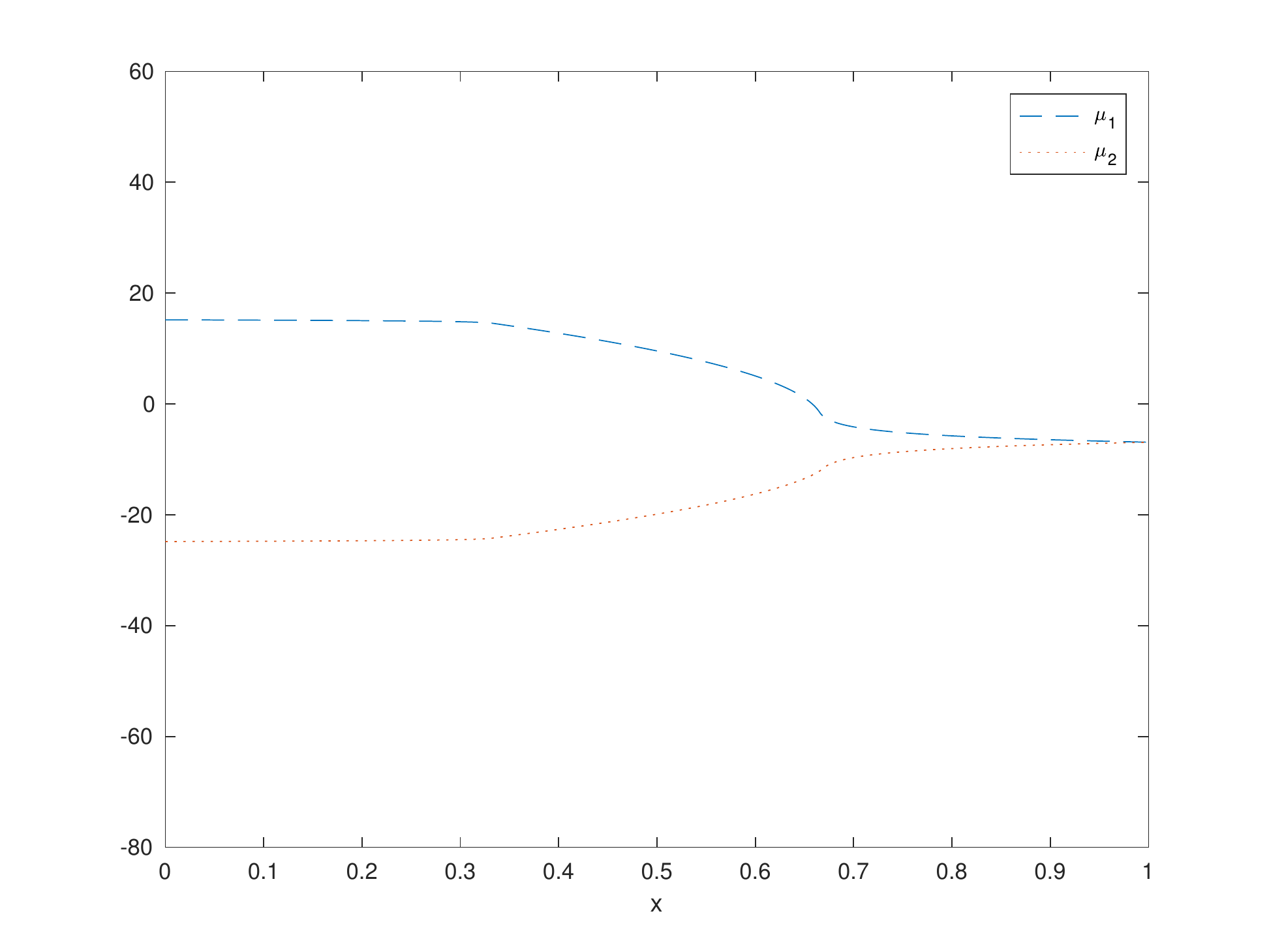} \label{fig:dynamics6}}

\subfigure[$c_1$ and $c_2$ at Point $1$]{\includegraphics[width = 5cm]{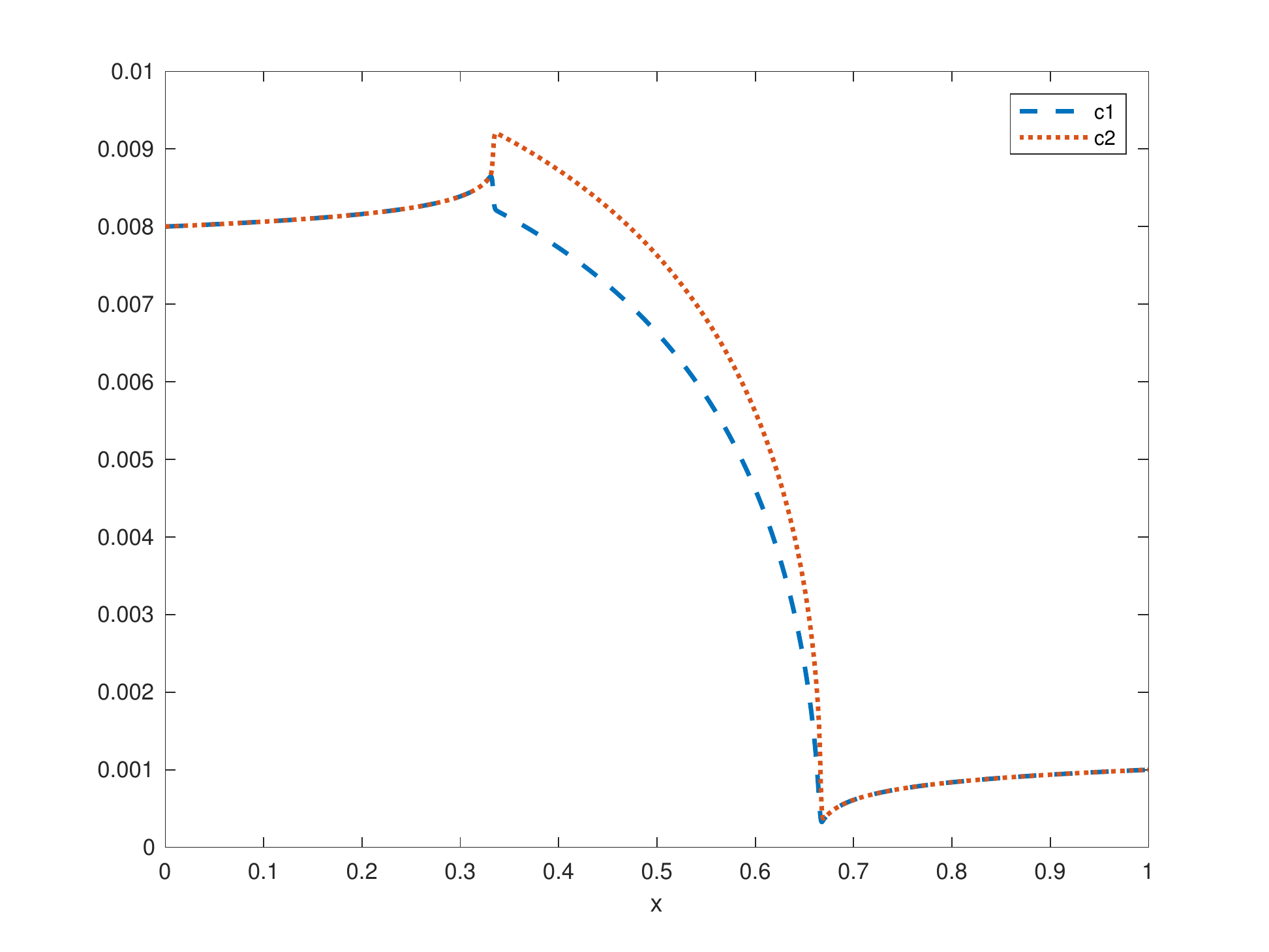} \label{fig:dynamics7}}
\subfigure[$c_1$ and $c_2$ at Point $2$]{\includegraphics[width = 5cm]{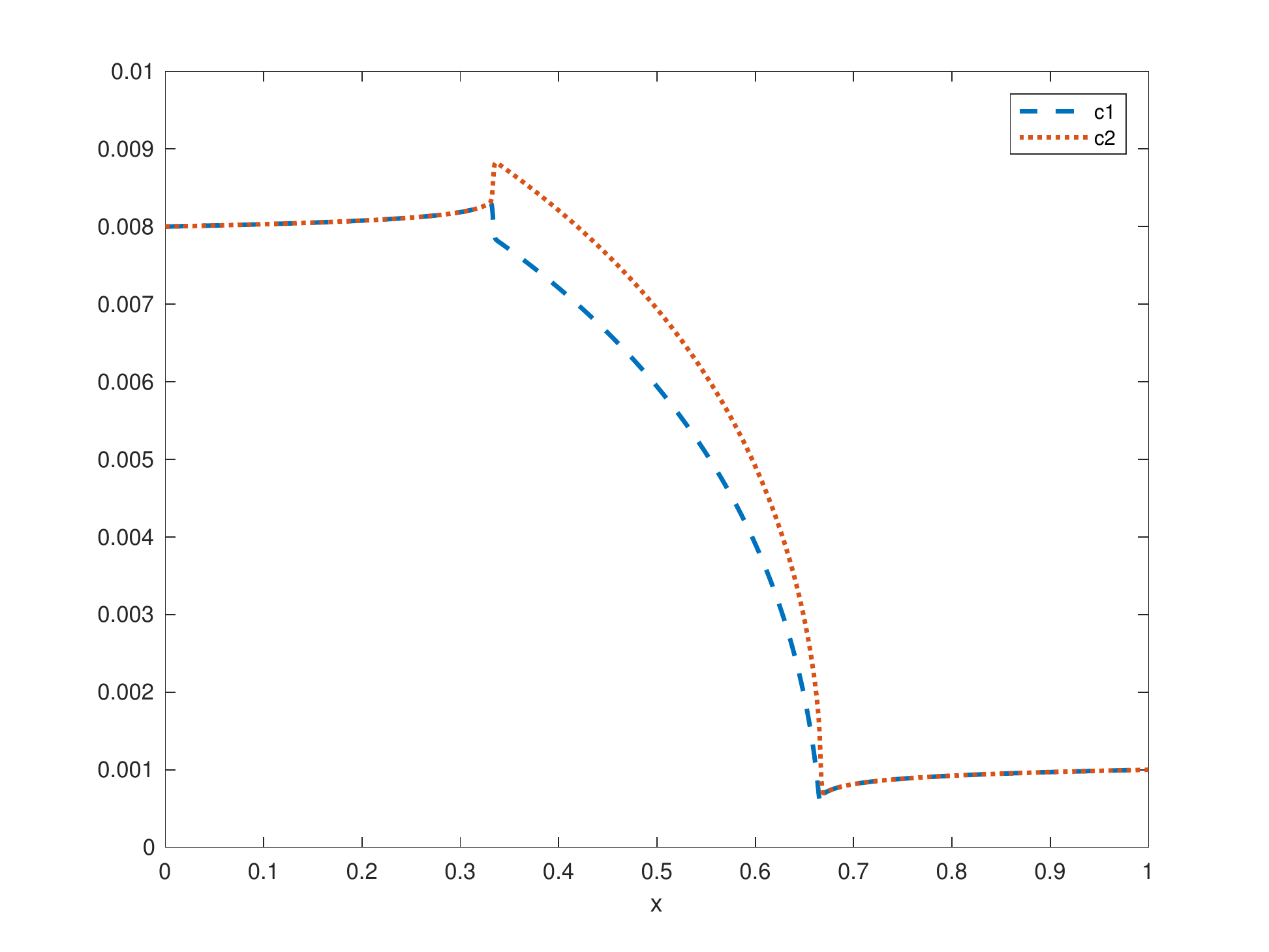} \label{fig:dynamics8}}
\subfigure[$c_1$ and $c_2$ at Point $3$]{\includegraphics[width = 5cm]{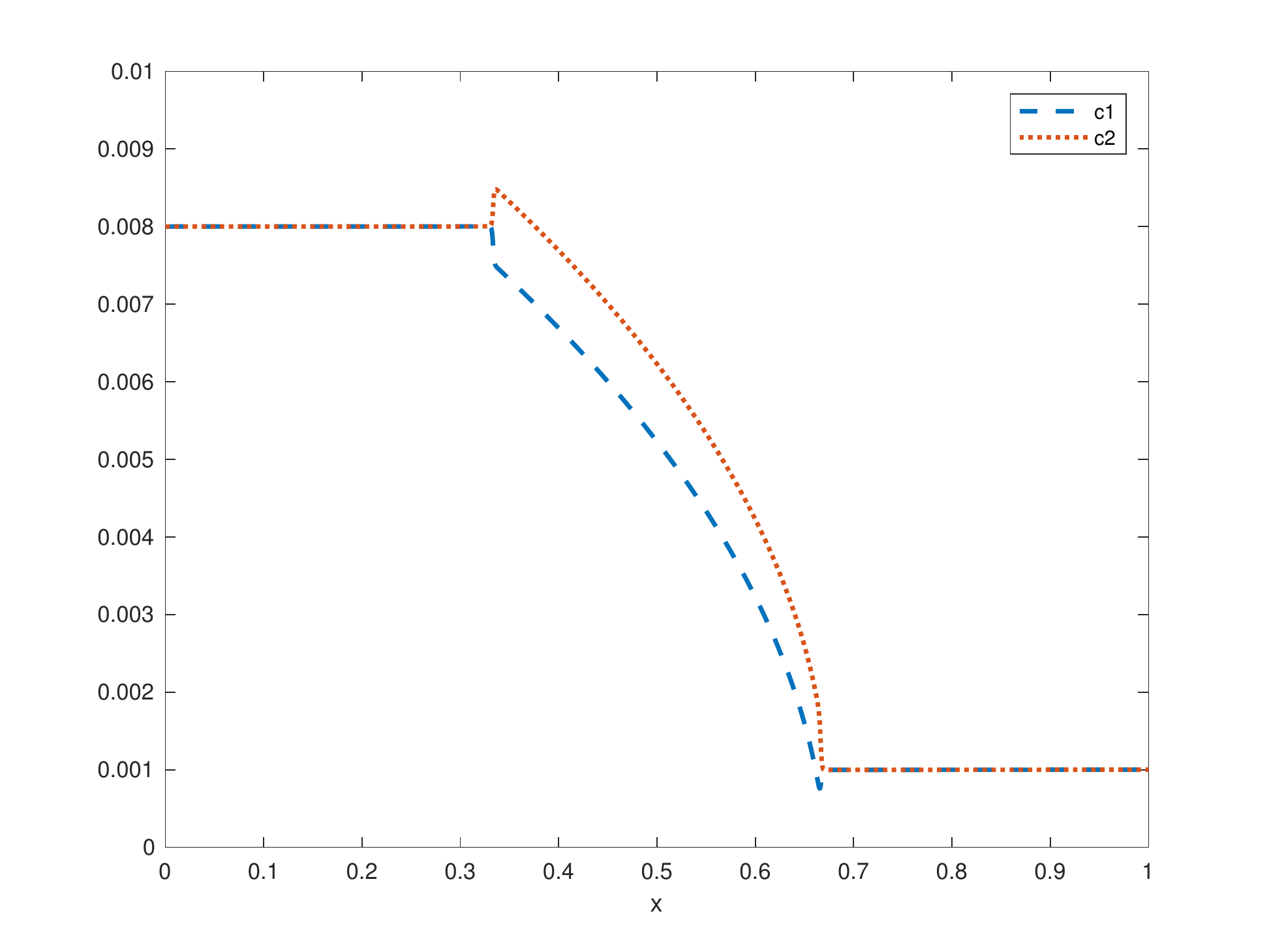} \label{fig:dynamics9}}

\subfigure[$\phi$ at Point $1$]{\includegraphics[width = 5cm]{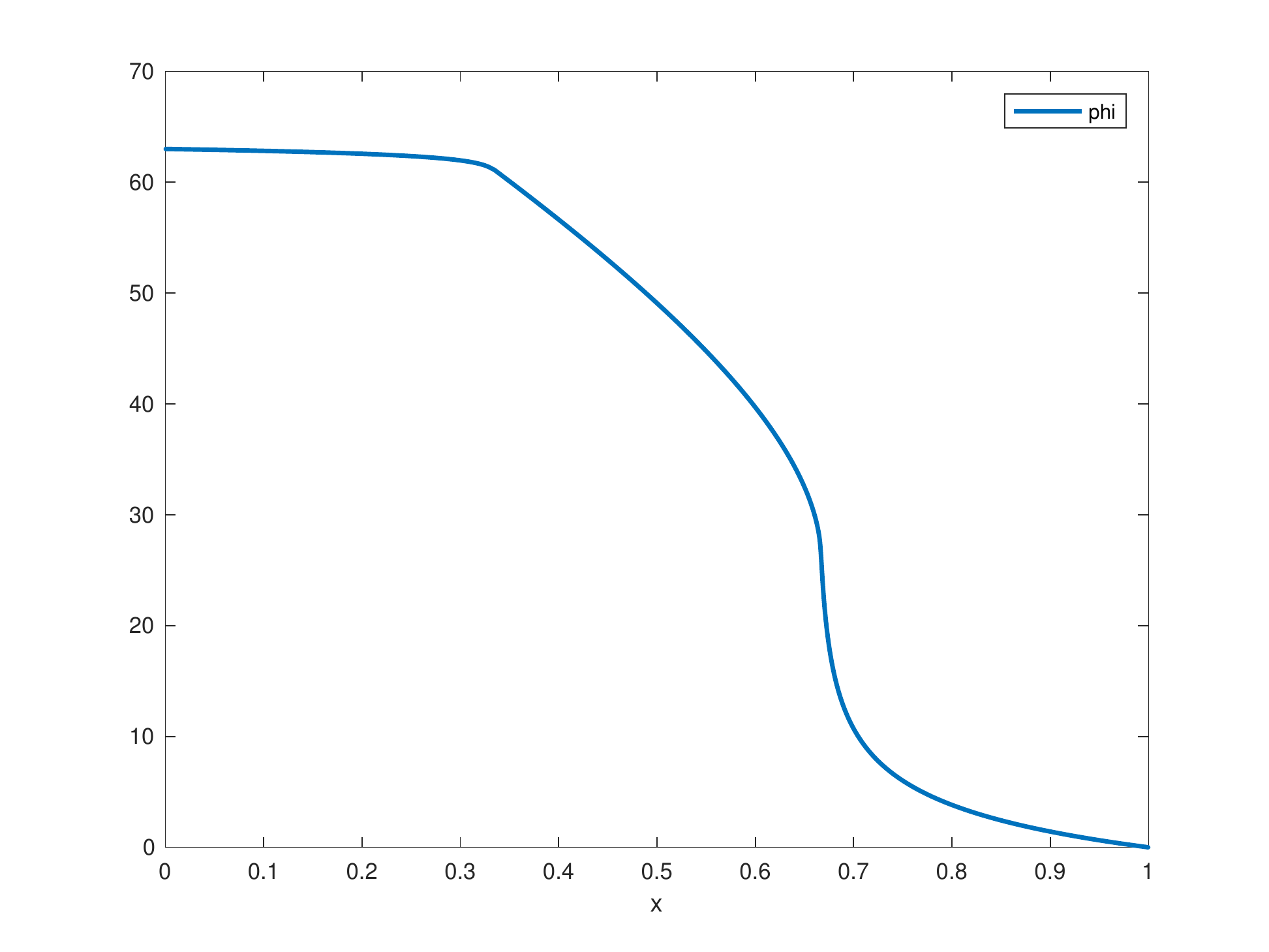} \label{fig:dynamics10}}
\subfigure[$\phi$ at Point $2$]{\includegraphics[width = 5cm]{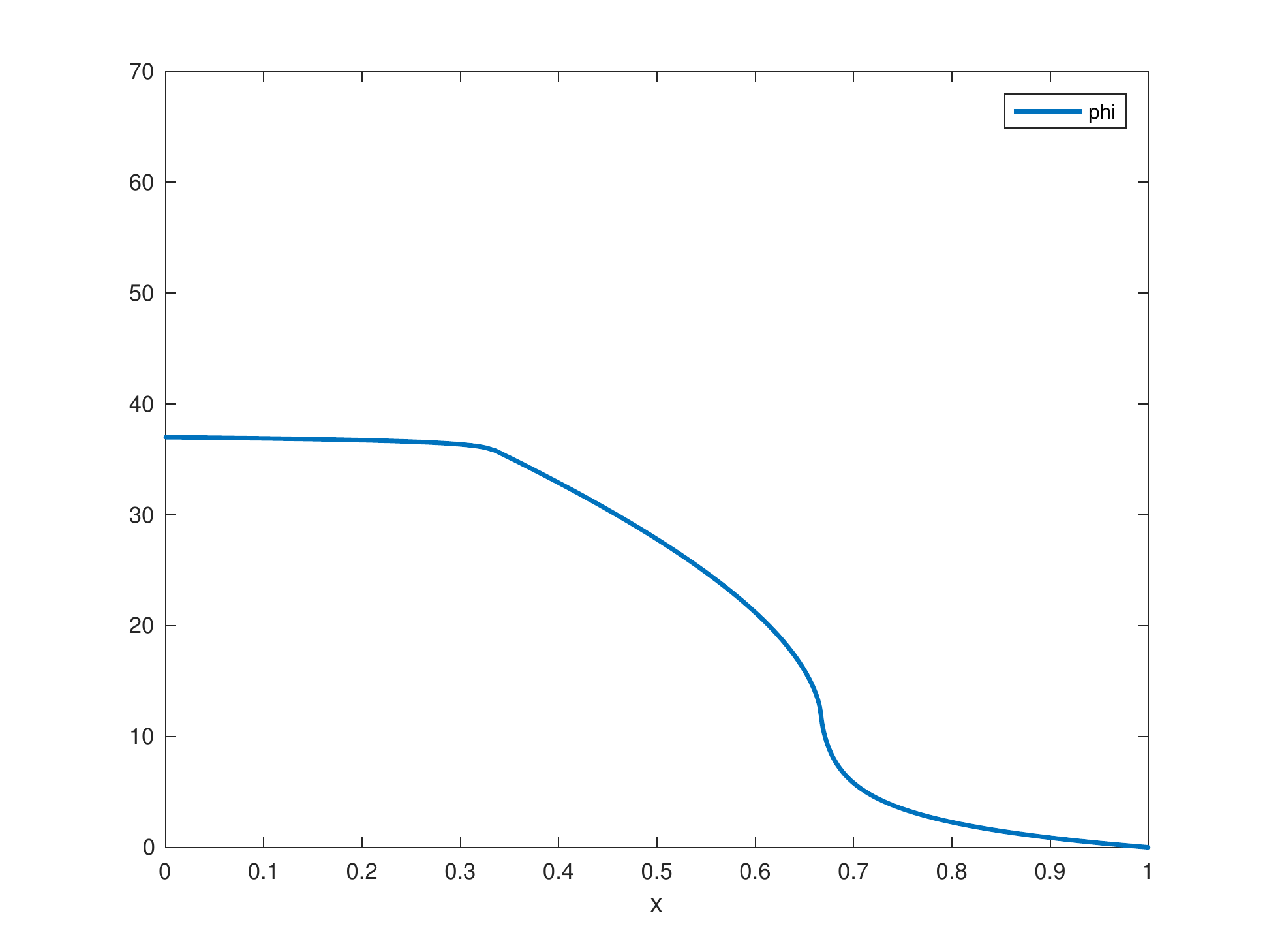} \label{fig:dynamics11}}
\subfigure[$\phi$ at Point $3$]{\includegraphics[width = 5cm]{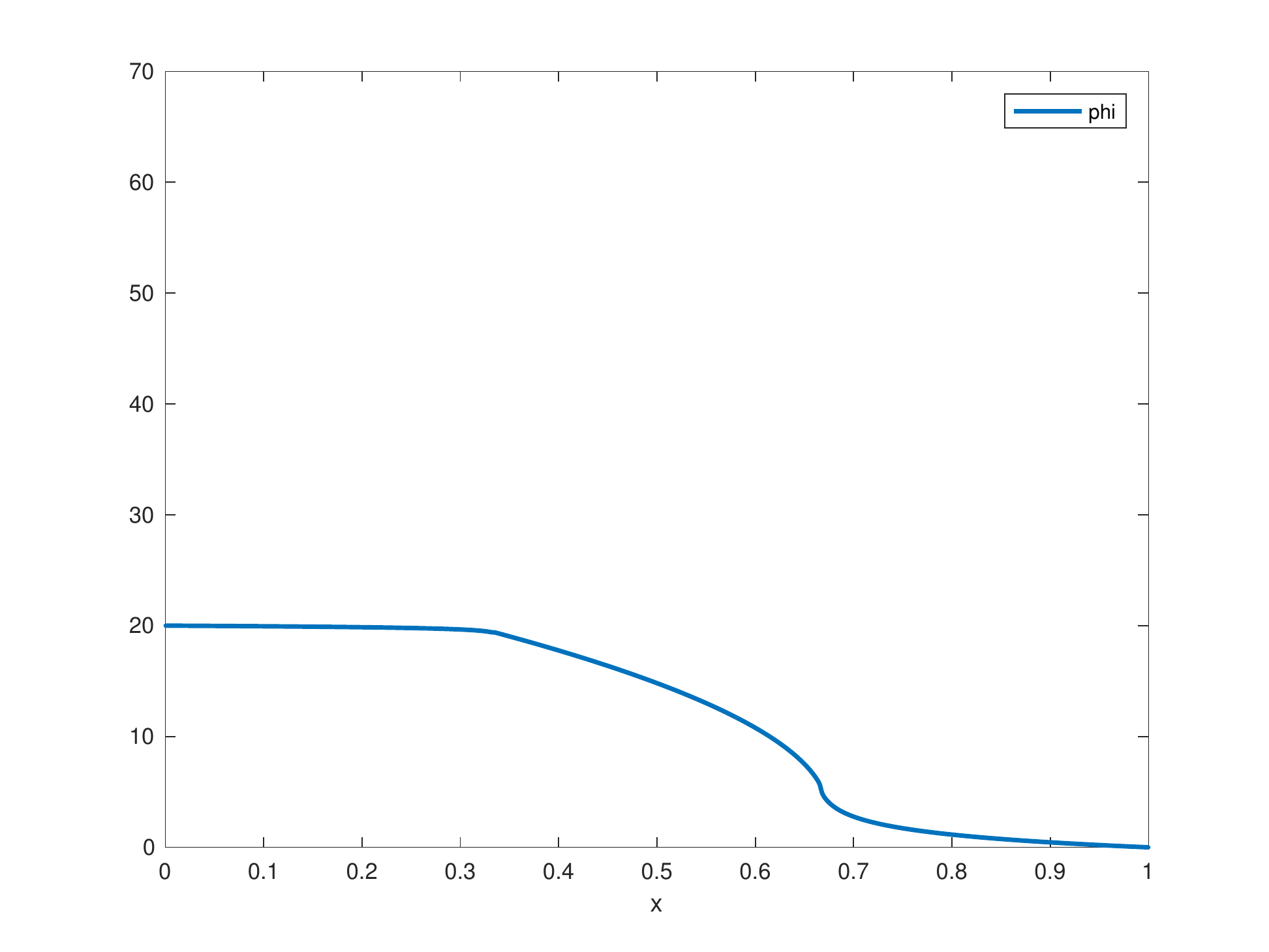} \label{fig:dynamics12}}
\caption{Internal dynamics at Point $1$ $(0.0005, 67)$, Point $2$ $(0.0005, 33)$, and Point $3$ $(0.0005, 20)$ in Fig.~\ref{L8R1complete}. \label{fig:Dynamics}}
\end{figure}

In Fig.~\ref{fig:Dynamics}, we present the flux and current, electrochemical potential, concentration of each species, and electric potential at Point $1$ $(0.0005, 67)$, Point $2$ $(0.0005, 33)$, and Point $3$ $(0.0005, 20)$ in Fig.~\ref{L8R1complete}. One can observe from Fig.~\ref{L8R1complete} that these three points are on a path where $\lambda_1$ is first less than $1$, then greater than $1$, then less than $1$ again. 

\begin{figure}[htbp]
\subfigure[$2 Q_0 = 0.0005$.]{\includegraphics[width = 0.45\textwidth]{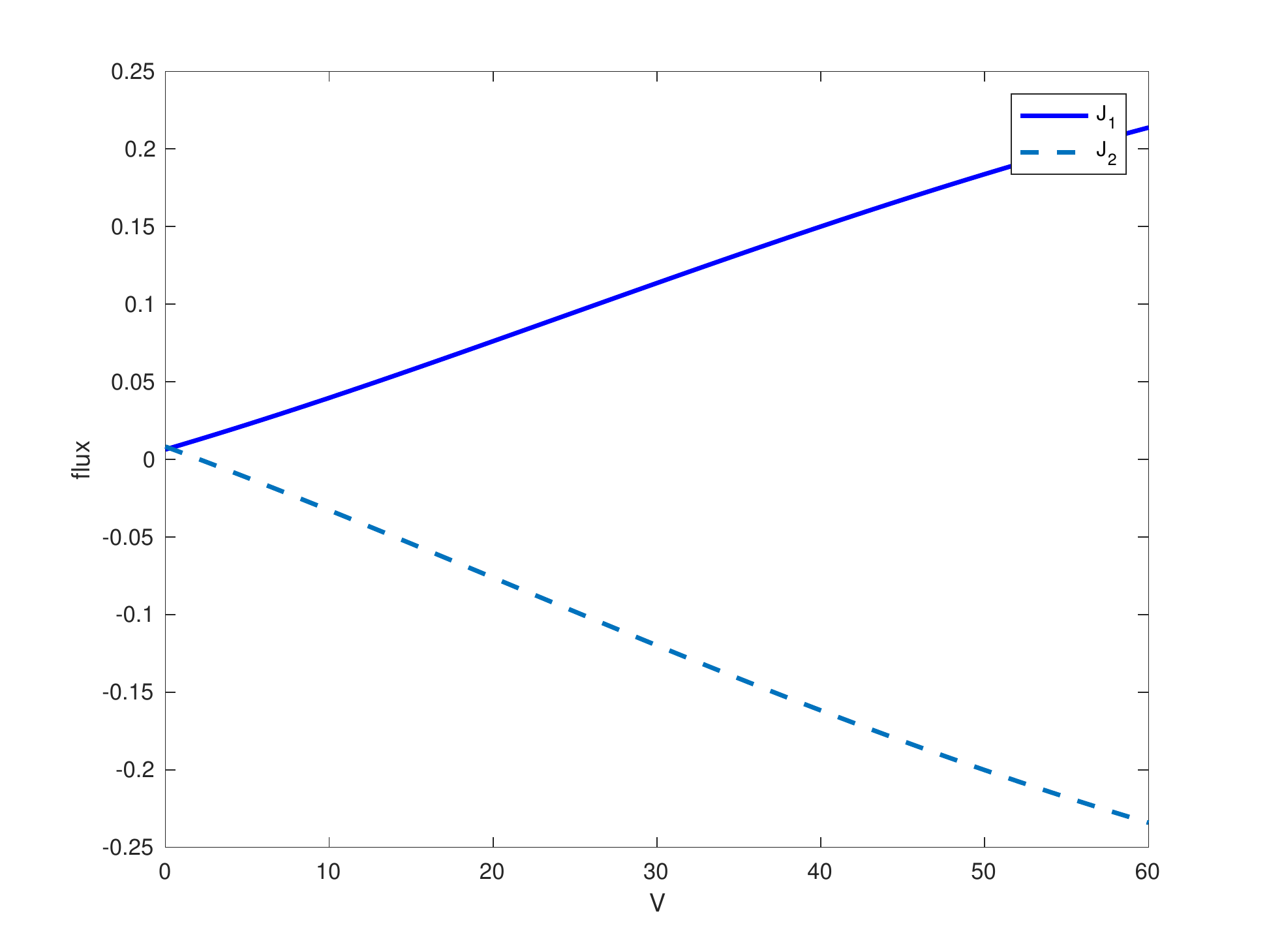} \label{fig:flux5e_4}}\quad
\subfigure[$2 Q_0 = 0.003$.]{\includegraphics[width = 0.45\textwidth]{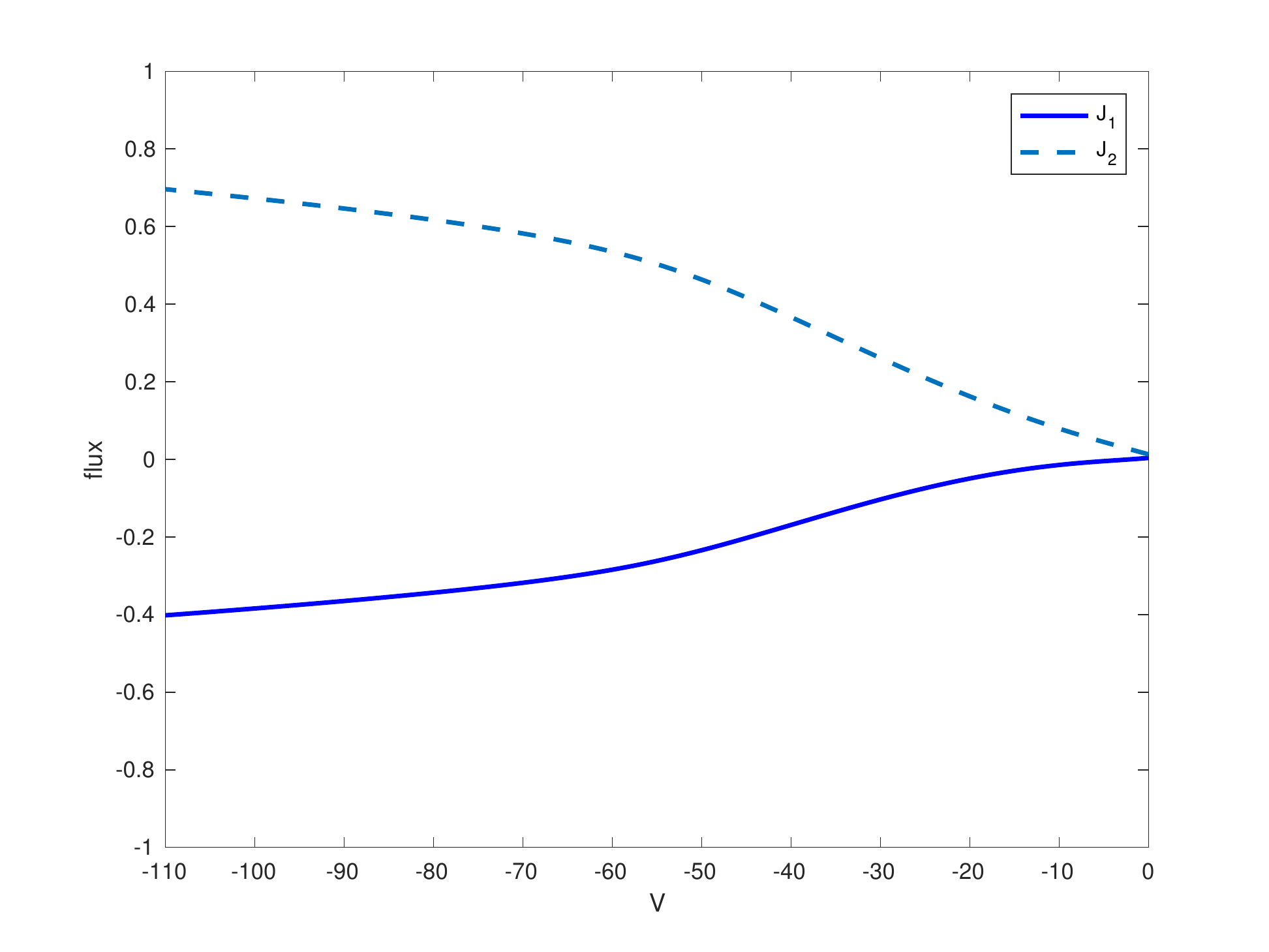} \label{fig:flux3e_3}}\quad
\caption{$J_1$ and $J_2$ as functions of $V$ with $2 Q_0 = 0.0005$ and $0.003$. \label{fluxplot}}
\end{figure}

In Fig.~\ref{fig:dynamics1}-\ref{fig:dynamics3}, we observe that as $V$ increases, both the positive and negative species are enhanced. While we know that $\lambda_1 > 1$ at Point $2$ while $\lambda_1 < 1$ at Points $1$ and $3$, this trend is clearly displayed in flux. To further study the trend of the $J$-$V$ relation, we present the plots of $J_1$ and $J_2$ as functions of $V$ in Fig.~\ref{fig:flux5e_4} and Fig.~\ref{fig:flux3e_3}, where $2 Q_0 = 0.0005$ and $0.003$, respectively.
One observes that for small $Q_0$ such as $2 Q_0 = 0.0005$, the fluxes are almost linear in $V$. However, as we increase $Q_0$ to $2 Q_0 = 0.003$, the fluxes are not linear any more in $V$, as shown in Fig.~\ref{fig:flux3e_3}. This observation is consistent with the linear $I$-$V$ relation for small permanent charge densities and the saturation effect for large permanent charge densities. Moreover, because $J_1$ and $J_2$ are nearly linear in $V$ when $Q_0 = 0$ (cf. (\ref{JKsmallQ})), the nonlinearity as shown in Fig.~\ref{fig:flux3e_3} is consistent with the fact that $\lambda_1$ and $\lambda_2$ are not monotone with respect to $V$.

In Fig.~~\ref{fig:dynamics4}-\ref{fig:dynamics6}, we observe that the voltage enhances the electrochemical potentials of both species. In Fig.~\ref{fig:dynamics7}-\ref{fig:dynamics9} and Fig.~\ref{fig:dynamics10}-\ref{fig:dynamics12}, sharp layers of $\phi$ and $c_k$ are demonstrated at the edge of neck of the channel.

\subsection{The hard-sphere case}
\label{sec:HS}

So far we have treated ions as point charges and assumed the excess component $\mu_k^{ex} = 0$. However, ions have positive volumes, and can have large influence on channel selectivity especially when they have large sizes and densities. In this section, we apply a hard sphere (HS) PNP model \cite{QLCL2016} (also see \cite{JL2012,LTZ2012, SL2018}), derived by taking Taylor expansions of the excess chemical potential of the Hard-Sphere model, which has been obtained based on Rosenfeld's fundamental measure theory (FMT) \cite{FMT}.
We now have
\[ \bar{\mu}_k = \bar{\mu}_k^{id} + \bar{\mu}_k^{HS}, \]
where the dimensionless $\bar{\mu}_k^{HS}$ is given (cf. (\ref{HS})) by
\begin{equation}
\begin{split}
\bar{\mu}_k^{HS} (x)  & =  - \ln \Big( 1 - \Sigma_j  \frac{4}{3} \pi\bar{r}_j^3c_j(x)\Big) + \frac{\bar{r}_k \Sigma_j 4 \pi \bar{r}_j^2 c_j(x)}{1 - \Sigma_j \frac{4}{3} \pi \bar{r}_j^3 c_j(x)} \\ & \quad + \frac{4 \pi \bar{r}_k^2 \Sigma_j \bar{r}_j c_j(x)}{1 - \Sigma_j \frac{4}{3} \pi \bar{r}_j^3 c_j(x)}
+ \frac{4}{3} \pi \frac{\bar{r}_k^3 \Sigma_j c_j(x)}{1 - \Sigma_j\frac{4}{3} \pi \bar{r}_j^3 c_j(x)},
\label{muLHS}
\end{split}
\end{equation}
where $\bar{r}_j=r_jC_0^{1/3}$ is the dimensionless radius of the $j$th hard sphere ion species.

\begin{figure}[htbp]
\subfigure[$\mu^{ex} = 0$]{\includegraphics[width = 0.47\textwidth]{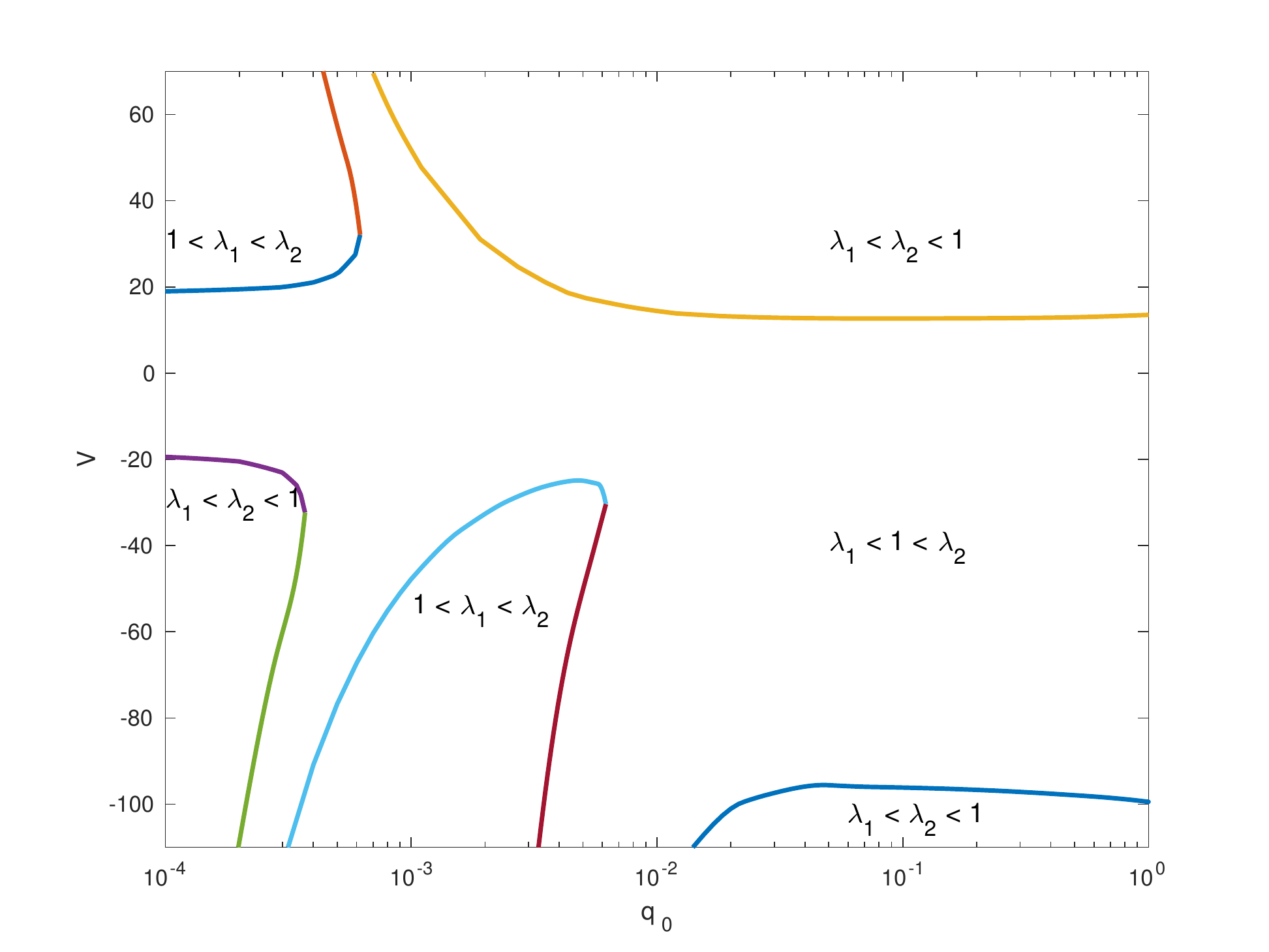} \label{fig:noLHSplot}}\qquad
\subfigure[$\mu^{ex}$ chosen as in (\ref{muLHS})]{\includegraphics[width = 0.47\textwidth]{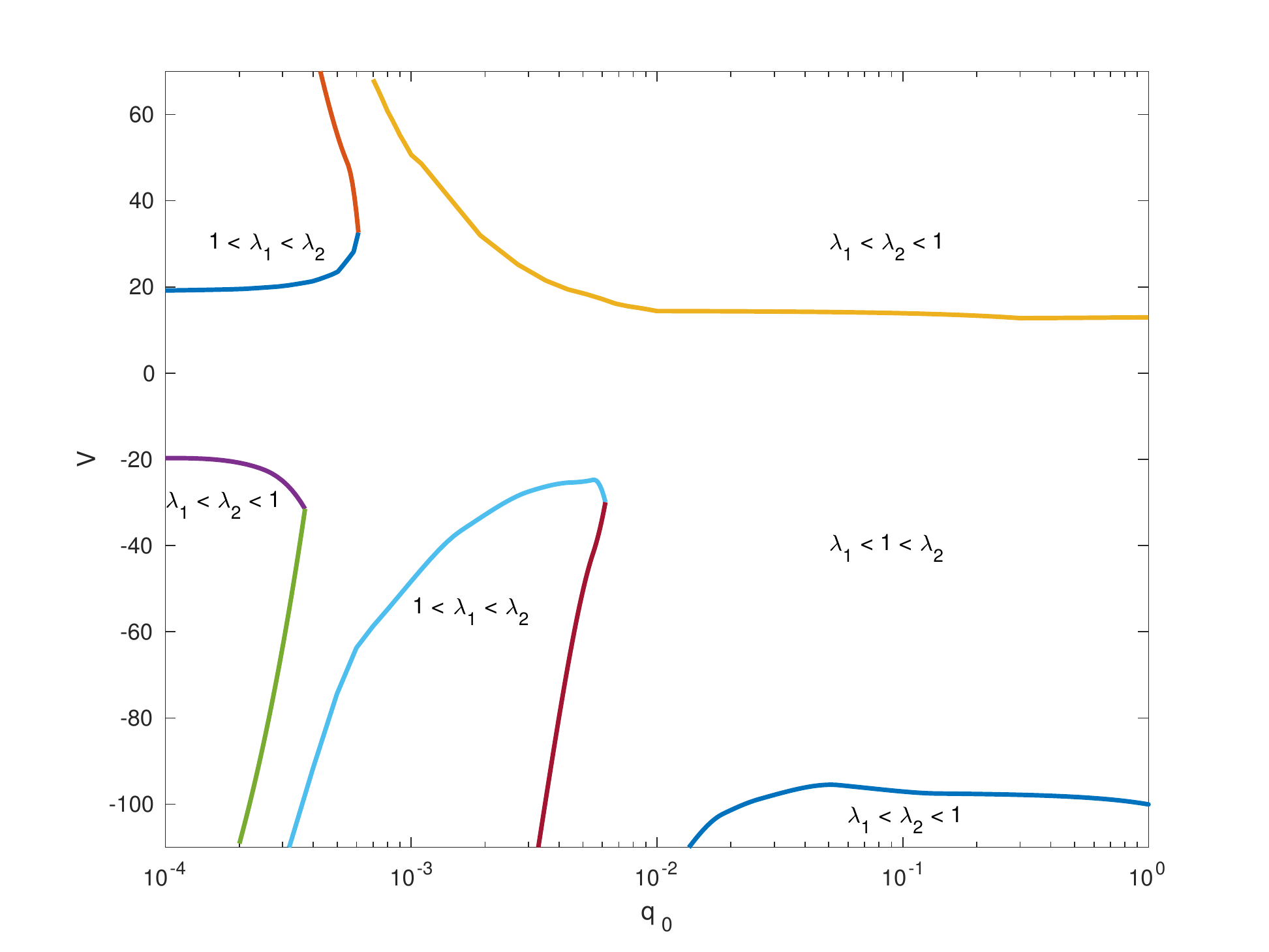} \label{fig:LHSplot}}
\caption{Comparison of bifurcation diagrams with $\mu^{ex} = 0$ and $\mu^{ex}$ chosen as in (\ref{muLHS}), with $L = 0.008$ and $R = 0.001$. Here, $q_0 = 2 Q_0$.
\label{fig:LHScomplete}}
\end{figure}
\begin{figure}[htbp]
\subfigure[ $\mu^{ex} = 0$]{\includegraphics[width = 0.47\textwidth]{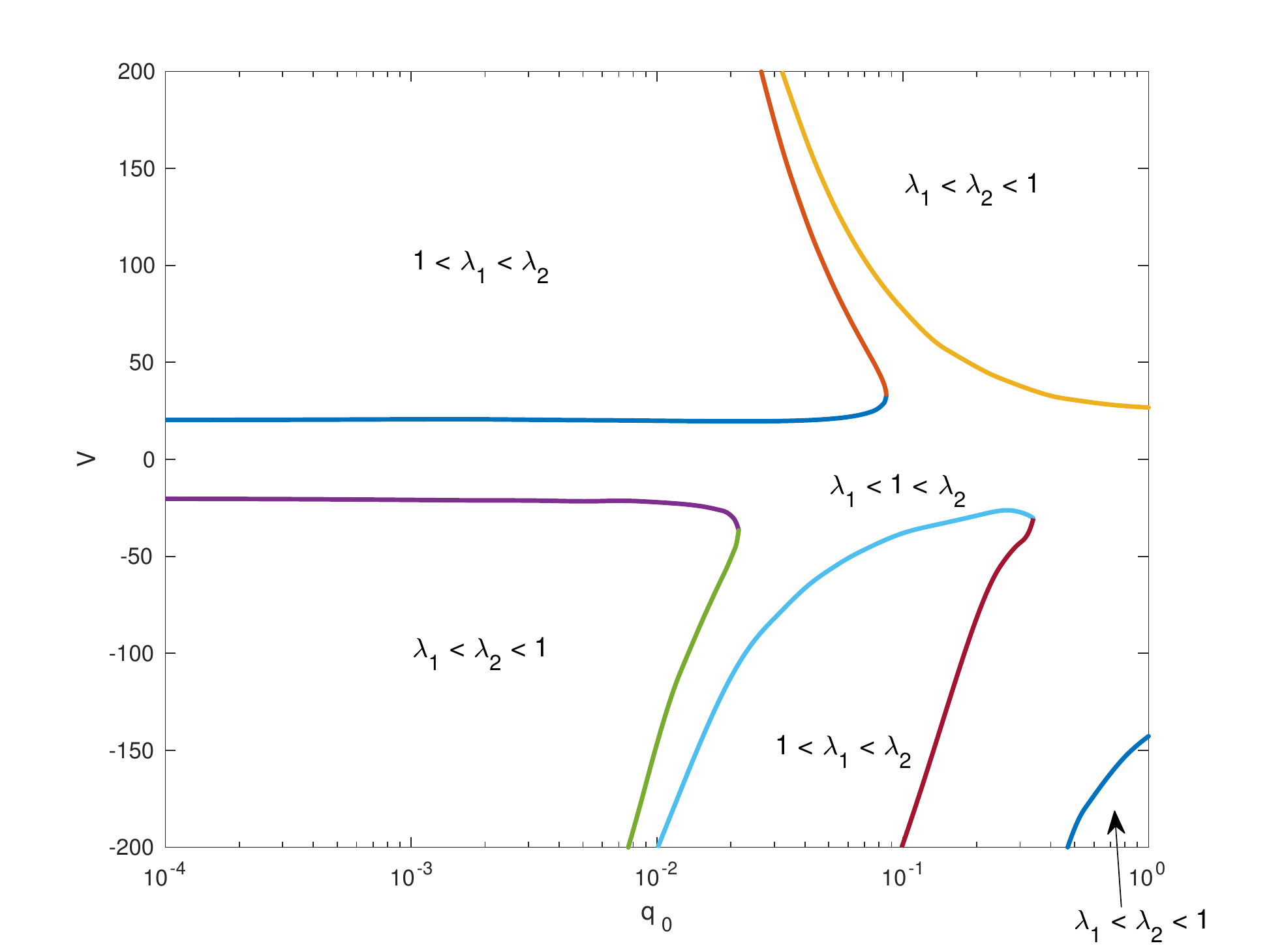} \label{fig:noLHSplotL5R1}}\qquad
\subfigure[$\mu^{ex}$ chosen as in (\ref{muLHS})]{\includegraphics[width = 0.47\textwidth]{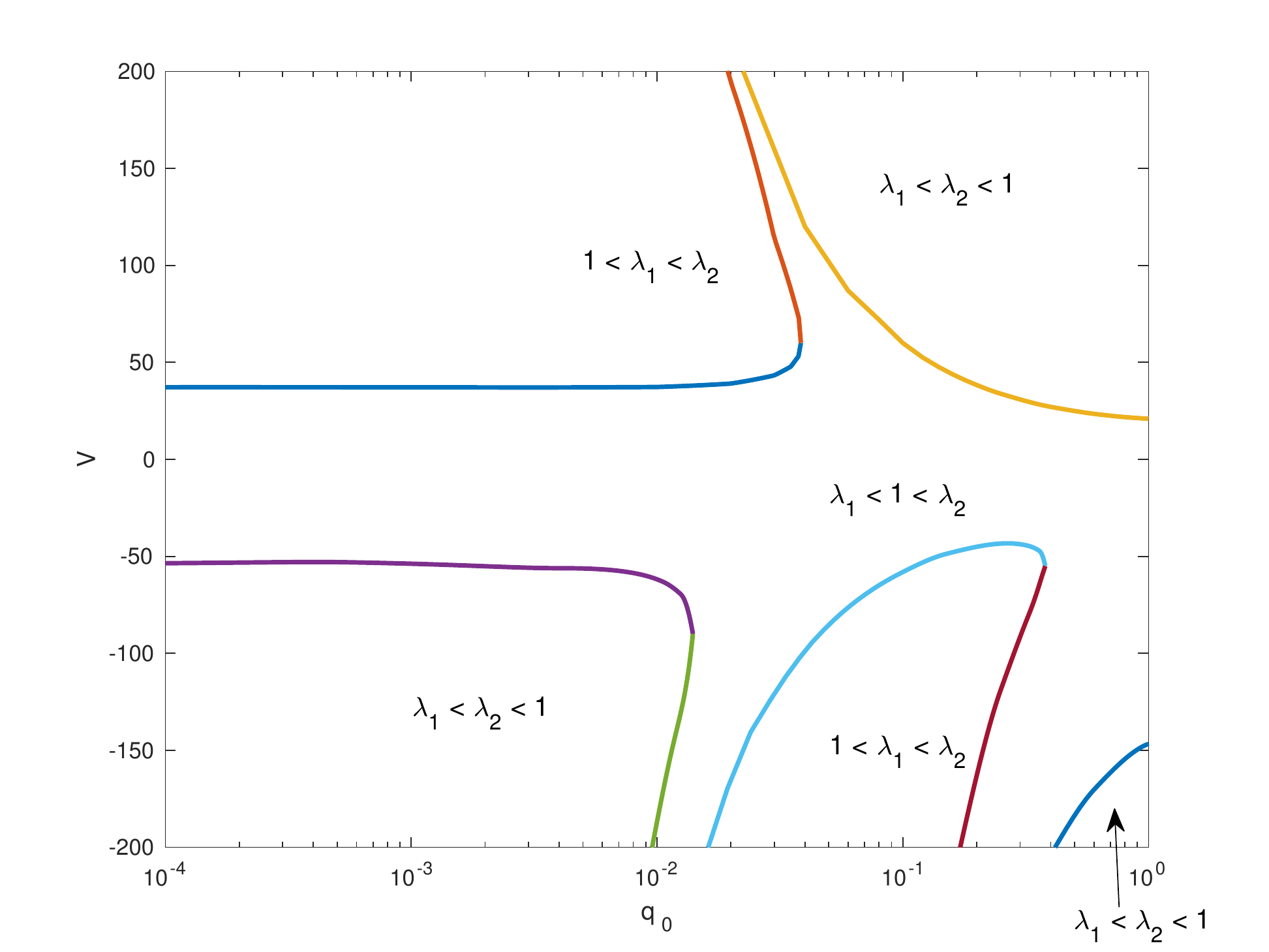} \label{fig:LHSplotL5R1}}
\caption{Comparison of bifurcation diagrams with $\mu^{ex} = 0$ and $\mu^{ex}$ chosen as in (\ref{muLHS}), with $L = 0.5$ and $R = 0.1$. Here, $q_0 = 2 Q_0$.
\label{fig:LHSL5R1}}
\end{figure}

In Fig.~\ref{fig:LHSplot}, we present the $Q_0$-$V$ bifurcation diagram with boundary conditions $L = 0.008$ and $R = 0.001$ and the radii of the species chosen as $\bar{r}_1 = 0.2$ and $\bar{r}_2 = 0.4$. Fig.~\ref{fig:LHSplot} is qualitatively and quantitively similar to Fig.~\ref{fig:noLHSplot}, which is associated with $\bar{\mu}^{ex} = 0$. If we change the boundary conditions to be $L = 0.5$ and $R = 0.1$, and compare the $Q_0$-$V$ bifurcation diagrams obtained with and without excess component, we can easily observe difference in quantities, but similar bifurcation properties, as shown in Fig.~\ref{fig:LHSL5R1}. This difference in quantity for large concentrations is expectable because, when the concentration is larger, more collision happens between ions, in which case the effects of the ion sizes play a more important role.

\begin{figure}[htbp]
\subfigure[ $\Sigma_j \frac{4}{3} \pi \bar{r}_j^3 c_j(x)$]{\includegraphics[width = 5cm,height = 4.2cm]{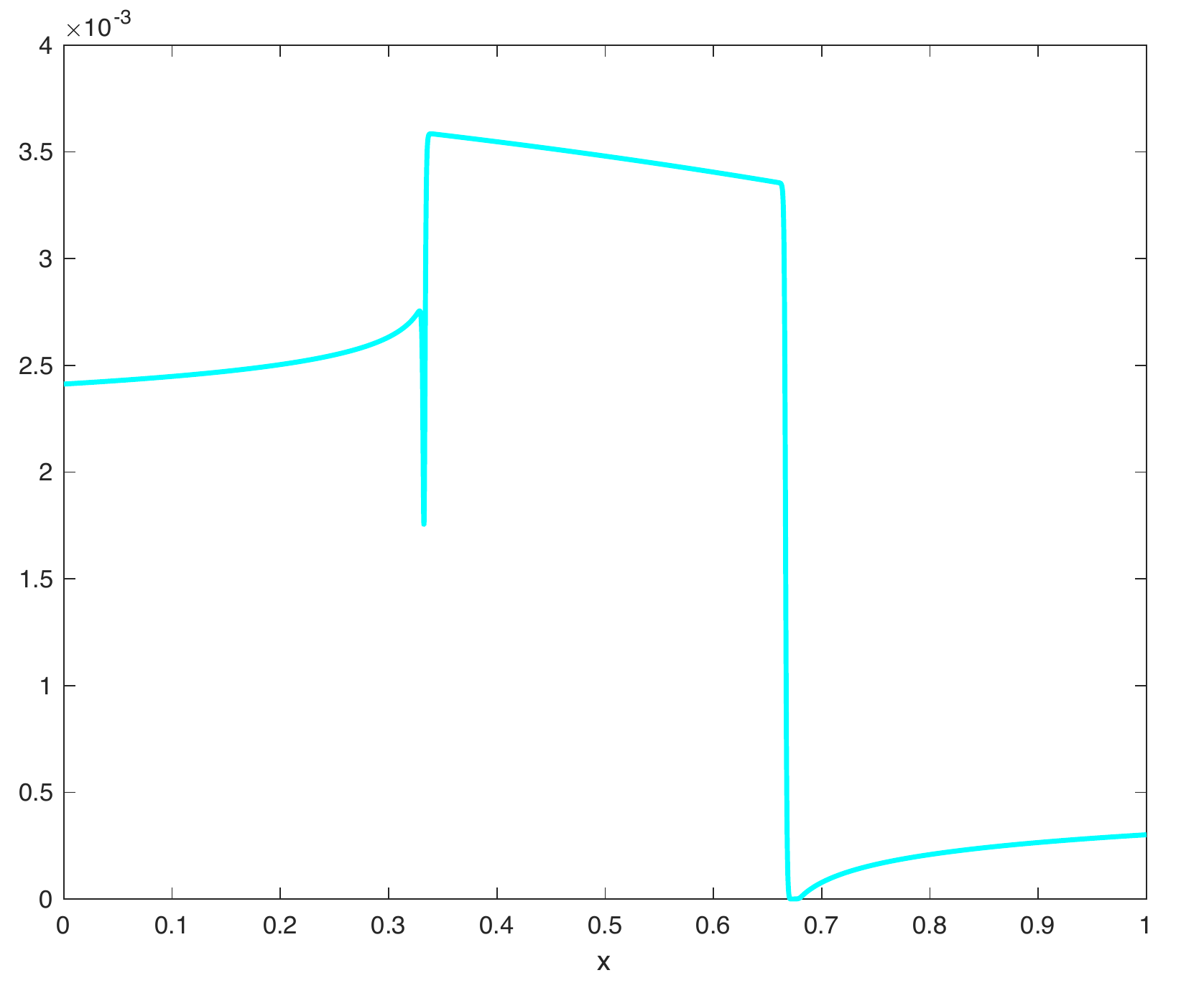} \label{fig:Q5e_2V30portion}}
\subfigure[$\bar{\mu}_k^{id}$]{\includegraphics[width = 5cm,height = 4.2cm]{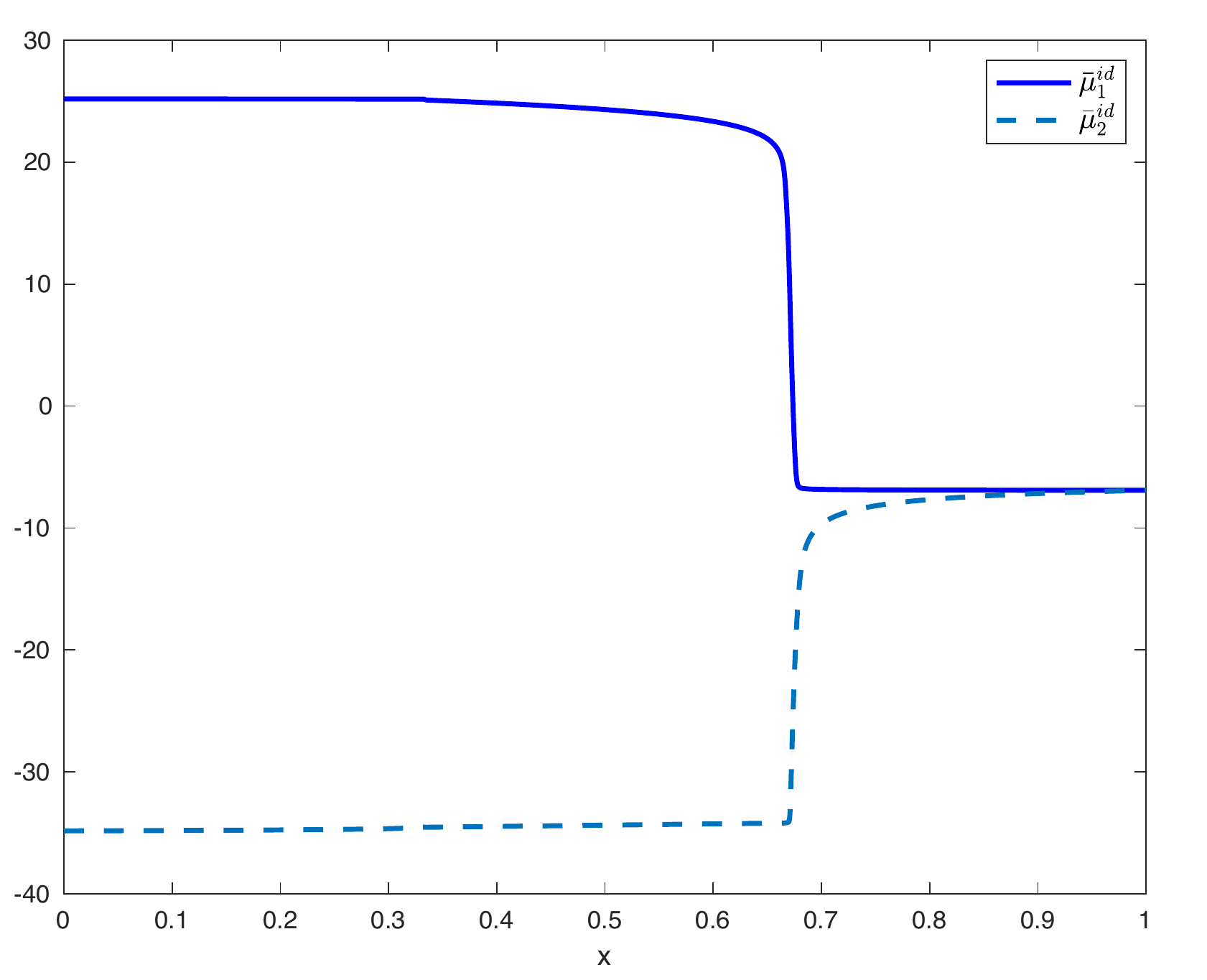} \label{fig:Q5e_2V30muid}}
\subfigure[$\bar{\mu}_k^{ex}$]{\includegraphics[width = 5cm,height = 4.2cm]{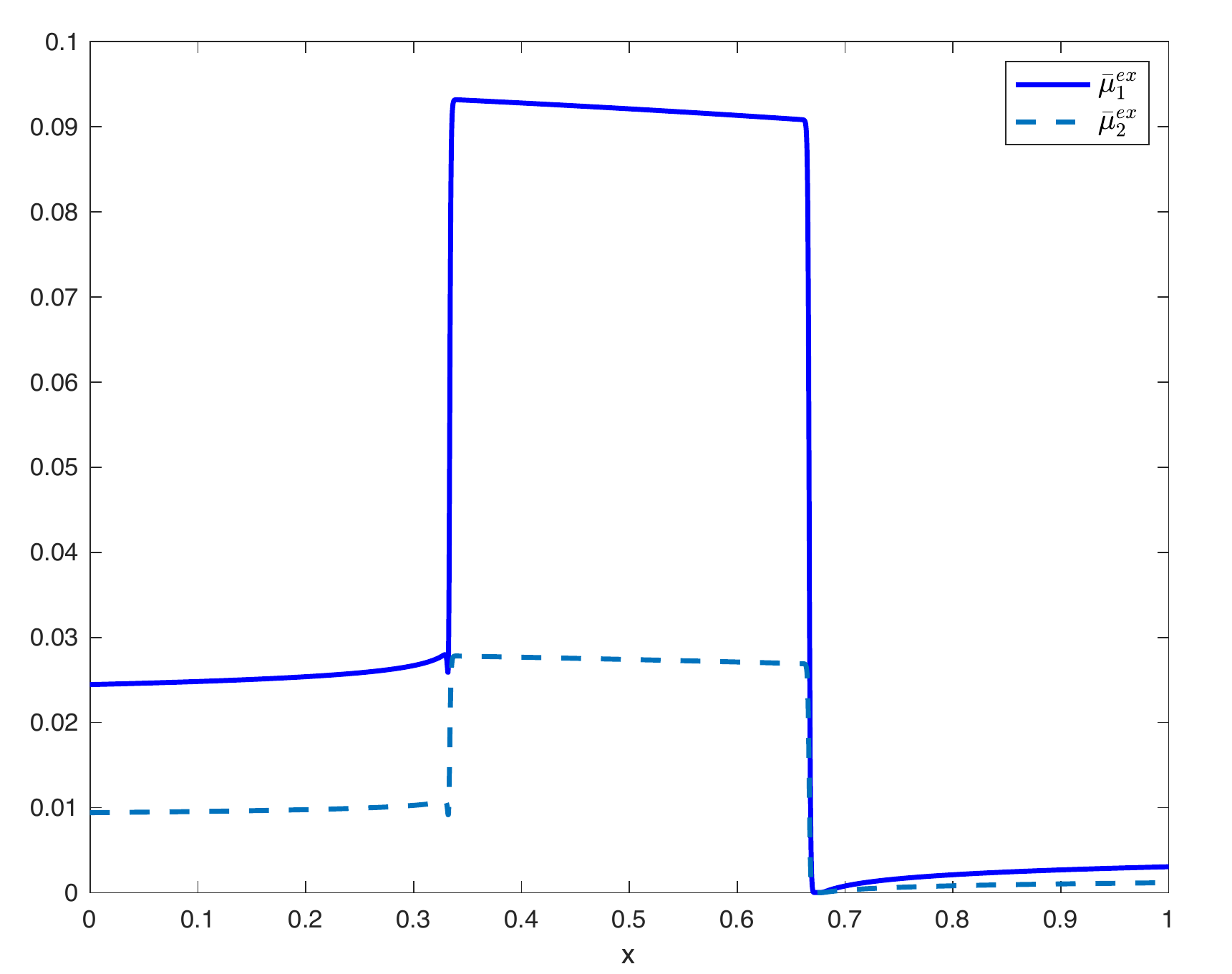} \label{fig:Q5e_2V30muex}}
\caption{$\Sigma_j \frac{4}{3} \pi \bar{r}_j^3 c_j(x)$, $\bar{\mu}_k^{id}$ and $\bar{\mu}_k^{ex}$, with $L = 0.008$, $R = 0.001$, $V = 30$, and $2 Q_0 = 0.05$.
\label{fig:LHSL8R1detail}}
\end{figure}

\begin{figure}[htbp]
\subfigure[ $\Sigma_j \frac{4}{3} \pi \bar{r}_j^3 c_j(x)$]{\includegraphics[width = 5cm,height = 4.2cm]{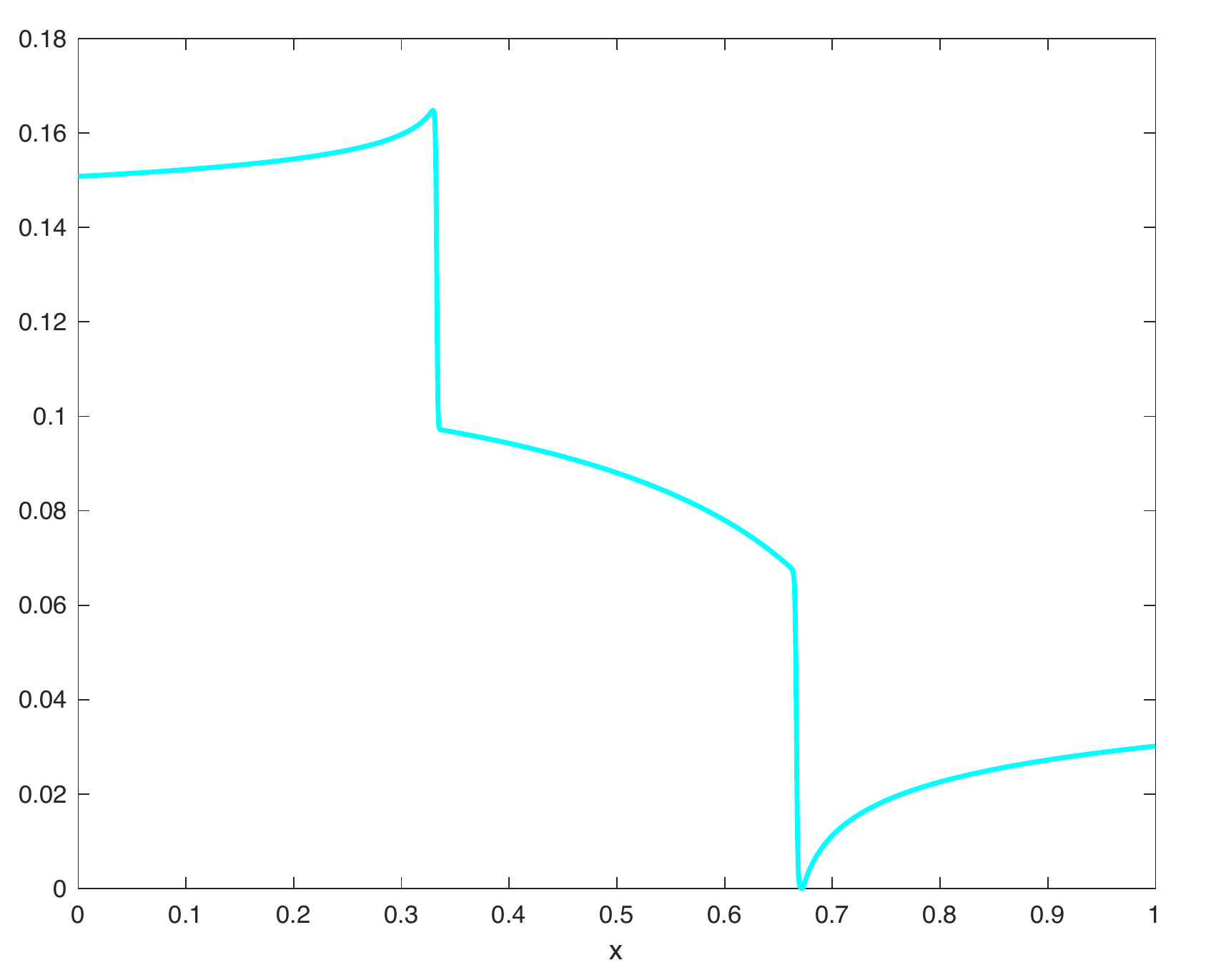} \label{fig:Q1V30portion}}
\subfigure[$\bar{\mu}_k^{id}$]{\includegraphics[width = 5cm,height = 4.2cm]{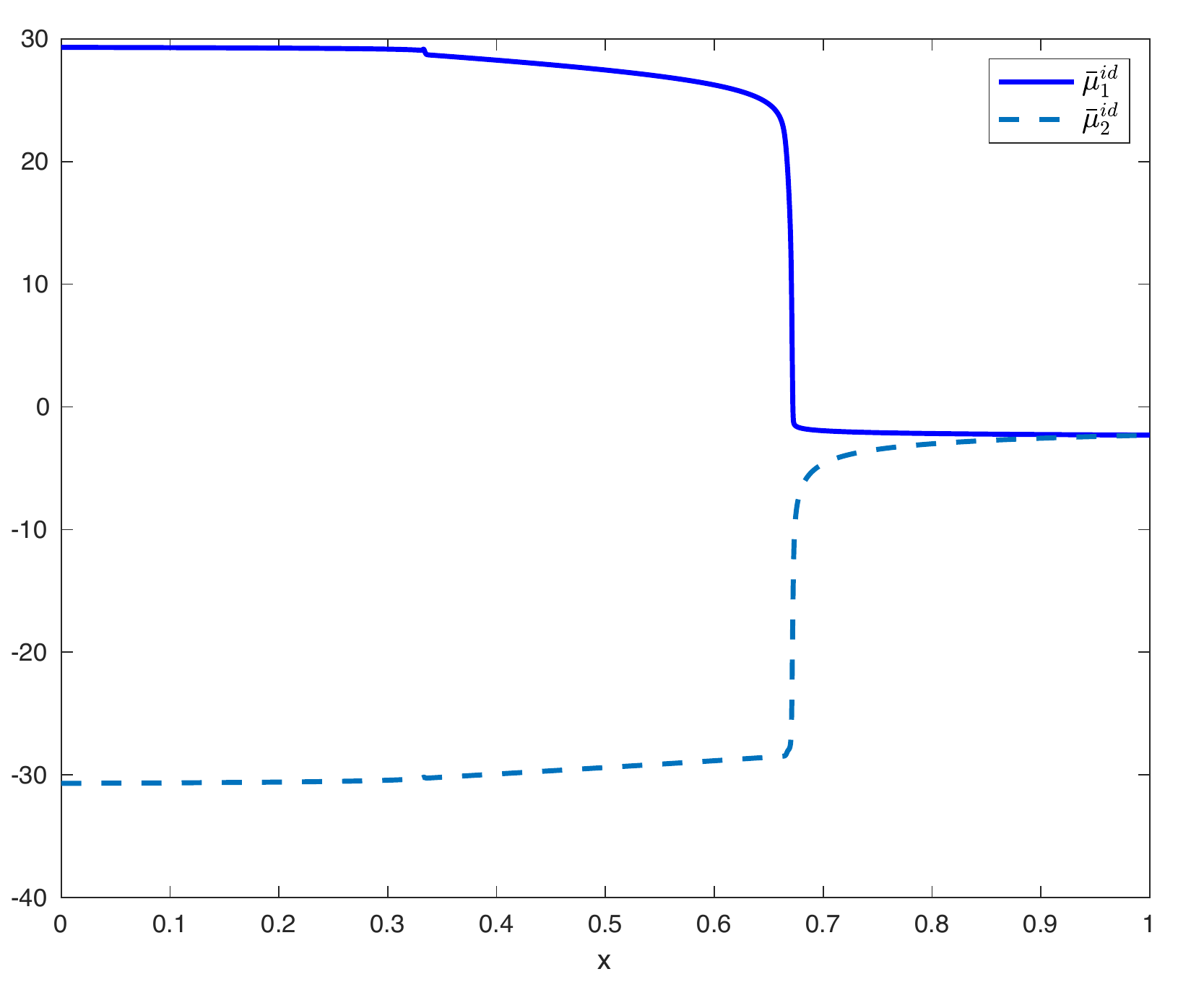} \label{fig:Q1V30muid}}
\subfigure[$\bar{\mu}_k^{ex}$]{\includegraphics[width = 5cm,height = 4.2cm]{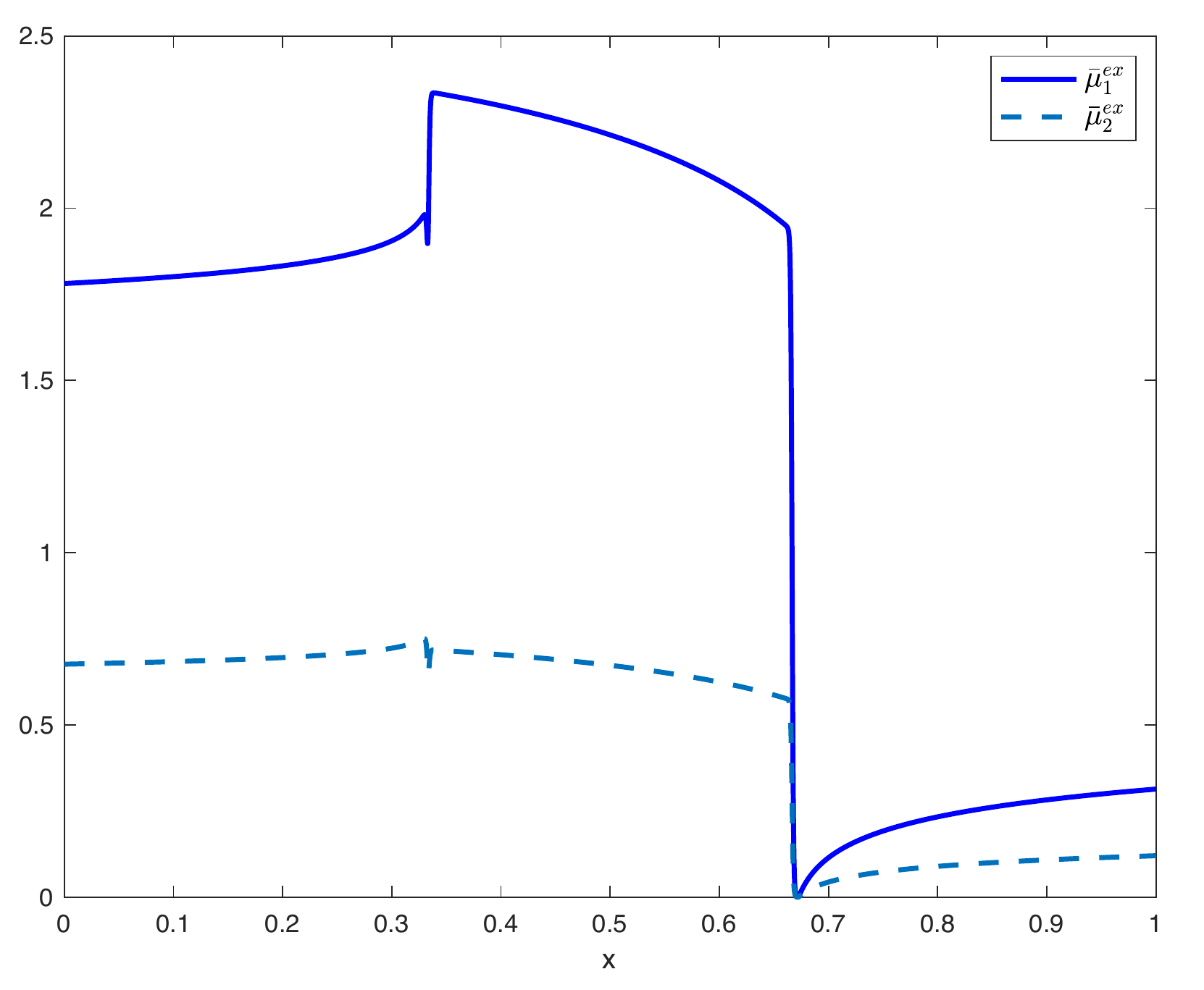} \label{fig:Q1V30muex}}
\caption{$\Sigma_j \frac{4}{3} \pi \bar{r}_j^3 c_j(x)$, $\bar{\mu}_k^{id}$ and $\bar{\mu}_k^{ex}$, with $L = 0.5$, $R = 0.1$, $V = 30$, and $2 Q_0 = 1$. \label{fig:LHSL5R1detail}}
\end{figure}

\begin{figure}[htbp]
\centering
\subfigure[ $\bar{\mu}_k$ when $\bar{\mu}^{ex}_k = 0$]{\includegraphics[width = 5cm,height = 4.2cm]{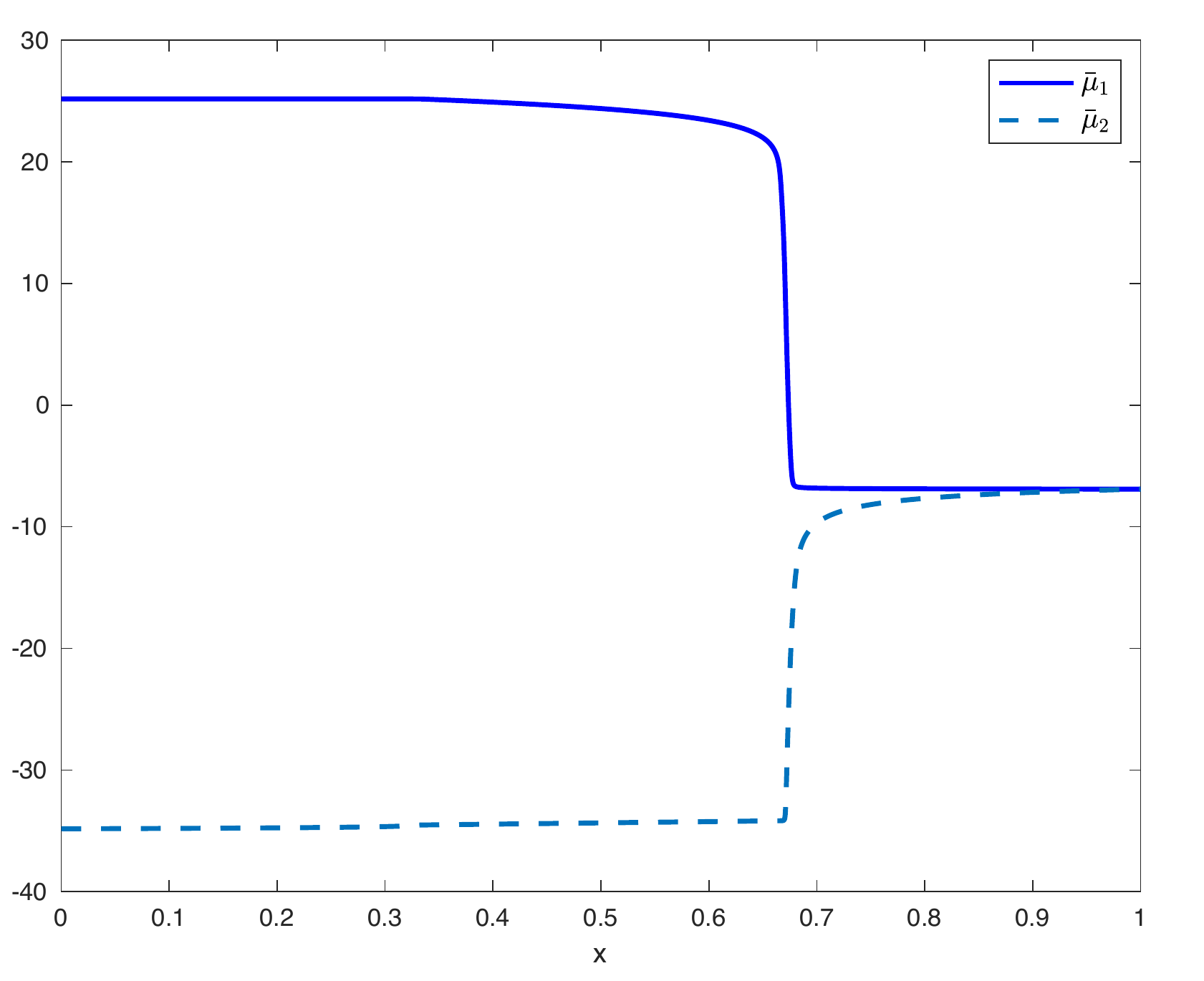} \label{fig:smallLmu}}
\subfigure[$\bar{\mu}_k$ when $\bar{\mu}^{ex}_k = \bar{\mu}_k^{HS}$]{\includegraphics[width = 5cm, height = 4.2cm]{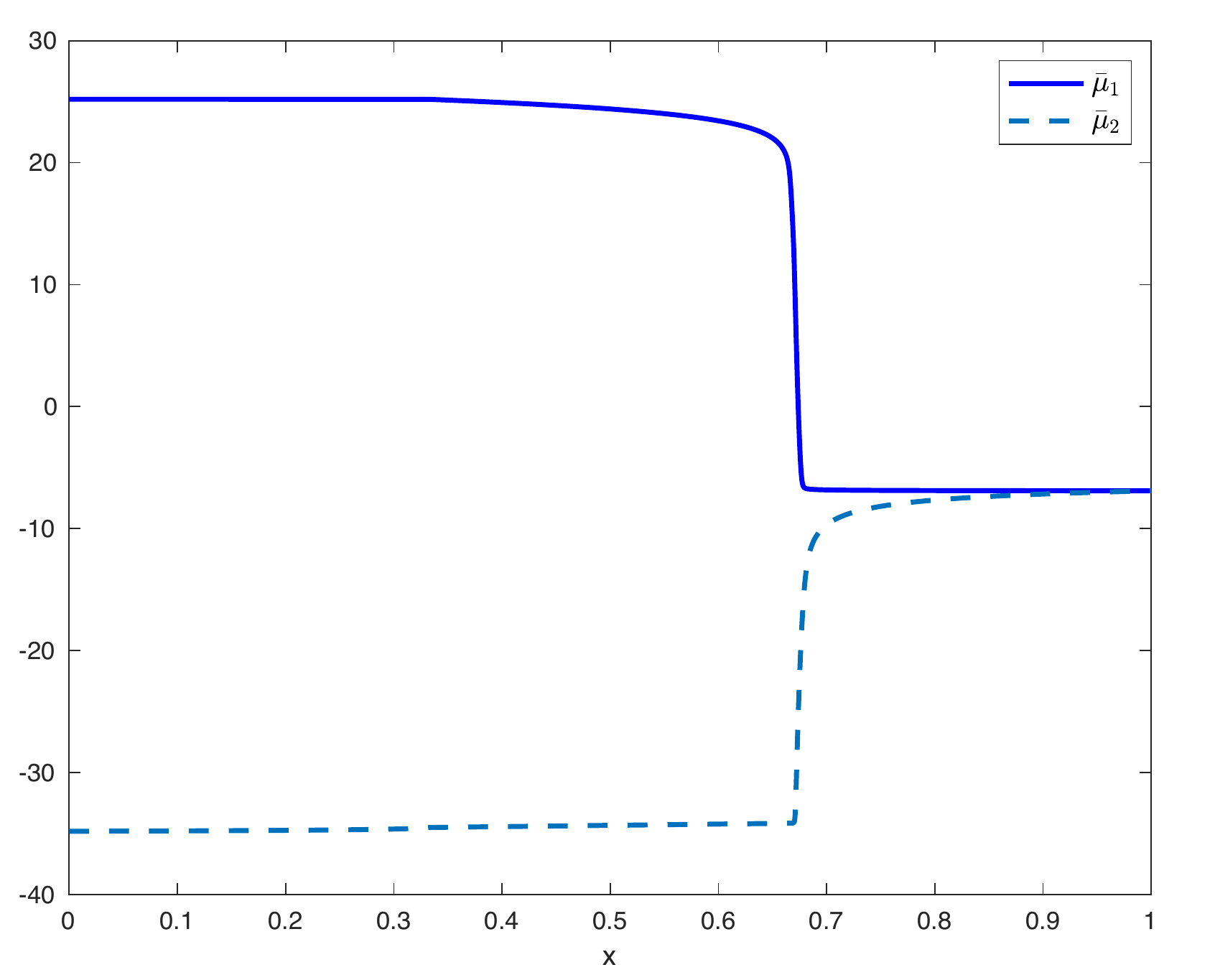} \label{fig:smallLmuLHS}}
\caption{The quantities of $\bar{\mu}_k$ for $L = 0.008$, $R = 0.001$, $V = 30$, and $2 Q_0 = 0.05$ with and without the HS term. \label{fig:smallL}}
\end{figure}

\begin{figure}[htbp]
\centering
\subfigure[ $\bar{\mu}_k$ when $\bar{\mu}^{ex}_k = 0$]{\includegraphics[width = 5cm,height = 4.2cm]{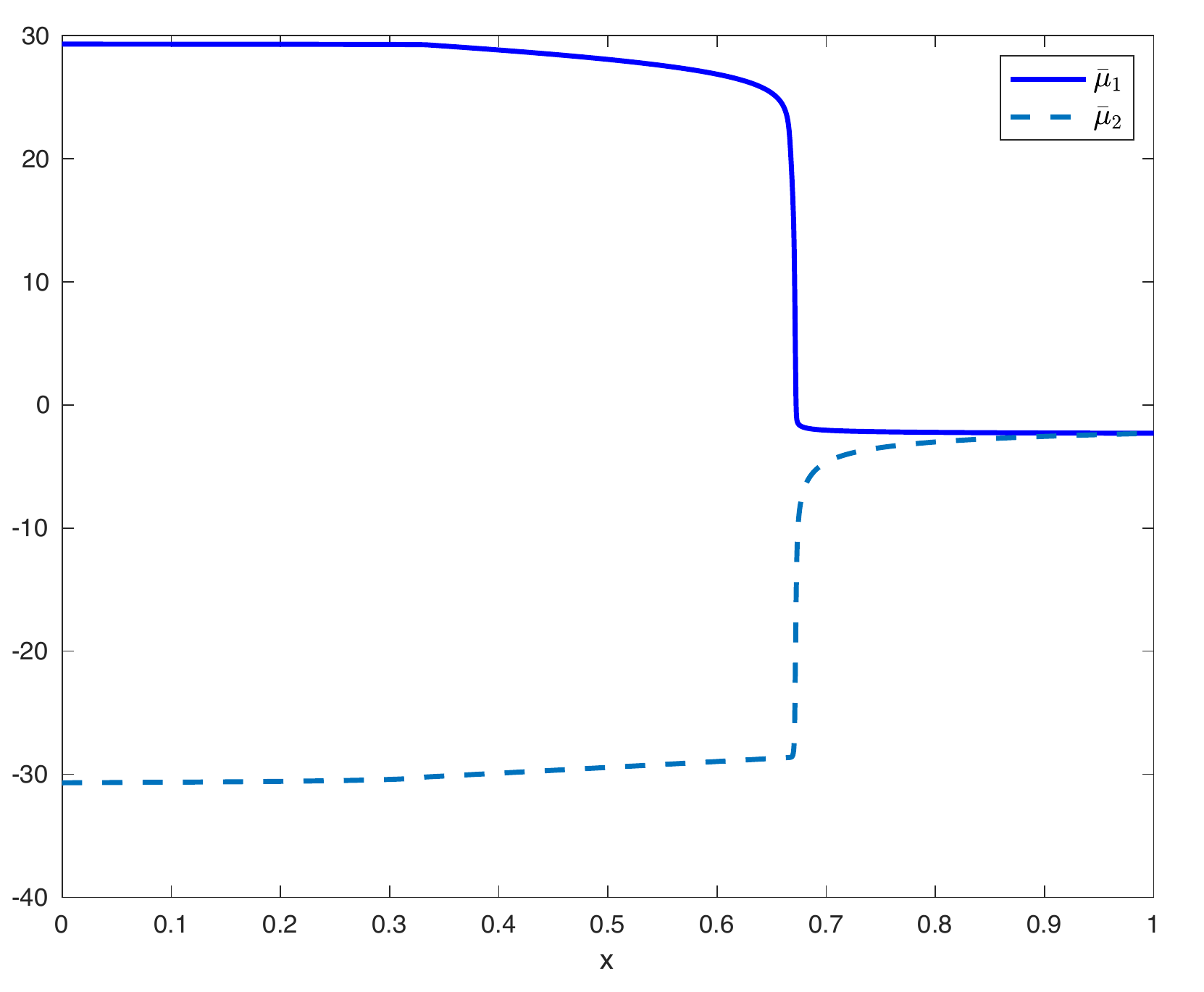} \label{fig:largeLmu}}
\subfigure[$\bar{\mu}_k$ when $\bar{\mu}^{ex}_k = \bar{\mu}_k^{HS}$]{\includegraphics[width = 5cm,height = 4.2cm]{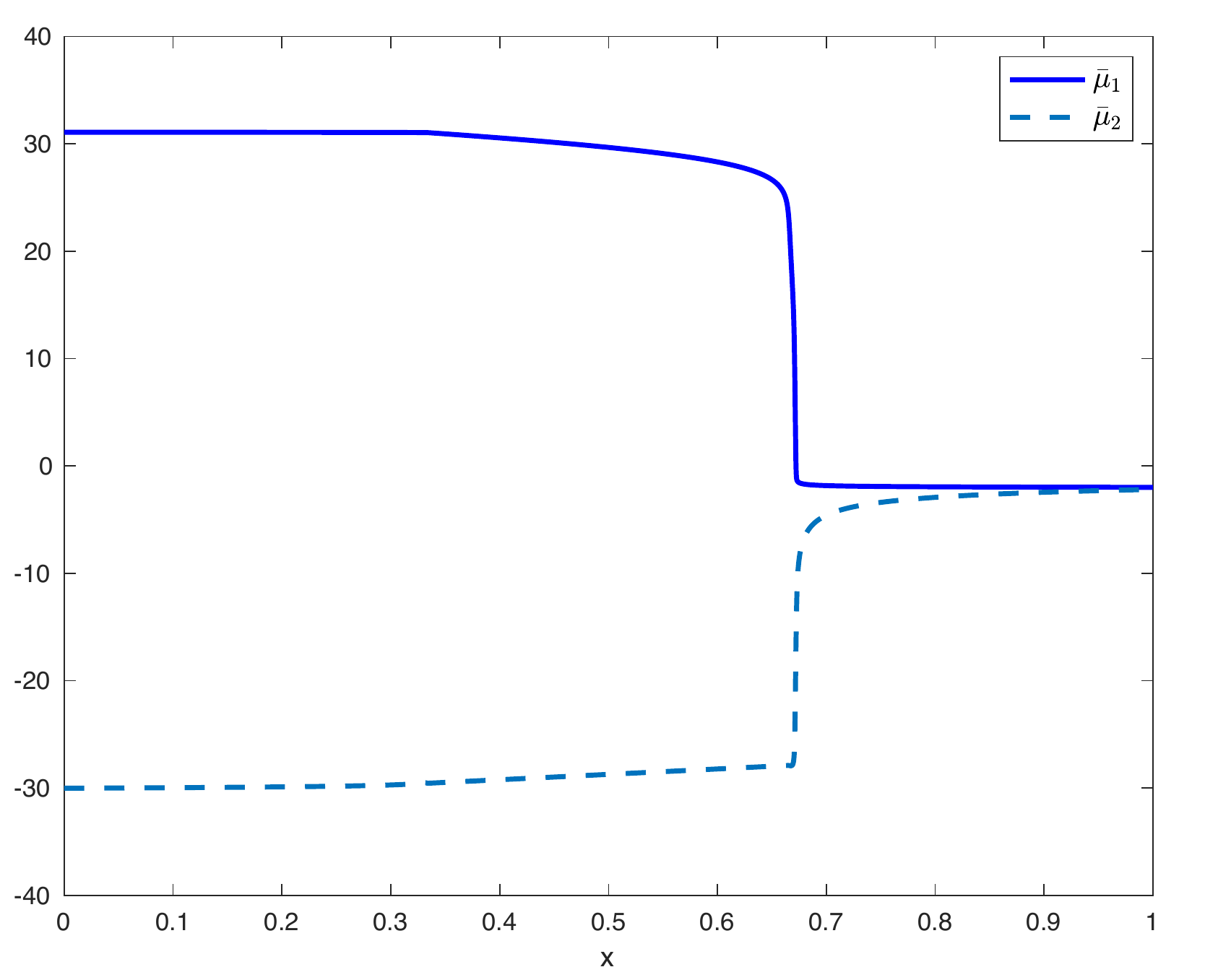} \label{fig:largeLmuLHS}}
\caption{The quantities of $\bar{\mu}_k$ for $L = 0.5$, $R = 0.1$, $V = 30$, and $2 Q_0 = 1$ with and without the HS term.\label{fig:largeL}}
\end{figure}

To quantify the effects of the ion sizes, we present the plots of the term $\Sigma_j \frac{4}{3} \pi \bar{r}_j^3 c_j(x)$, the ideal component $\bar{\mu}^{id}_k$, and the excess component $\bar{\mu}^{ex}_k$, in Fig.~\ref{fig:LHSL8R1detail} and Fig.~\ref{fig:LHSL5R1detail}, with different boundary conditions. As one can observe, $\Sigma_j \frac{4}{3} \pi \bar{r}_j^3 c_j(x)$ is almost negligible (compared to $1$) in Fig.~\ref{fig:Q5e_2V30portion}, with $L = 0.008$, $R = 0.001$, $V = 30$, and $2 Q_0 = 0.05$, and is relatively large in Fig.~\ref{fig:Q1V30portion}, where $L = 0.5$, $R = 0.1$, $V = 30$, and $2 Q_0 = 1$. As a result, $\bar{\mu}_k^{ex}$ in Fig.~\ref{fig:Q5e_2V30muex} is small compared to the ideal component $\bar{\mu}_k^{id}$ as shown in Fig.~\ref{fig:Q5e_2V30muid}, and does not make a significant difference in the flux ratio. By contrast, $\bar{\mu}_k^{ex}$ in Fig.~\ref{fig:Q1V30muex} is relatively large, correspondingly, with the same boundary conditions $L = 0.5$ and $R = 1$, and one can observe that the HS term has made a difference in the flux ratio diagram, as shown in Fig.~\ref{fig:LHSL5R1}. In order to demonstrate the effects of the HS term on the quantities of the electrochemical potential $\bar{\mu}_k$, we present $\bar{\mu}_k$ associated with the models with and without the term (\ref{muLHS})
in Fig.~\ref{fig:smallL} and Fig.~\ref{fig:largeL}.

\section{Numerical Method: An adaptive moving mesh finite element method}
\label{SEC:numerical}
In this section we describe an adaptive moving mesh finite element method for the numerical solution of the boundary value problem (BVP) (\ref{dimensionlessPNP}) and (\ref{bound}). Under the effects of the permanent charge and the shape of the channel, the solution to the BVP is known to have discontinuous second order derivatives at $x=1/3$ and $x=2/3$ \cite{ZEL2019}. For better resolution and improvements of the accuracy of the numerical simulation, we adopt the so-called moving mesh PDE (MMPDE) method to dynamically relocate and concentrate more mesh nodes around $x= 1/3$ and $x=2/3$.
For notational simplicity, the method is described here only for the situation where the electrochemical potential contains only the ideal component, i.e., $\bar{\mu}_k = \bar{\mu}_k^{id} = z_k \phi + \ln c_k$. Other situations can be treated without major modification.

It should be pointed out that other moving mesh strategies or refinement-based mesh adaptation methods can also be used for the numerical solution of BVP (\ref{dimensionlessPNP}) and (\ref{bound}). The interested reader is referred to \cite{Bai94a,Baines-2011,BHR09,HR11,Tan05} for literature and references therein. A main reason for choosing the MMPDE method is that it can be used in two and three spatial dimensions and for both time-independent and time-dependent problems \cite{HR11}. The method has also been shown analytically and numerically to produce a nonsingular dynamically varying mesh \cite{HK2015}. These features are important since our long-term goal is to study time-dependent three-dimensional PNP models.


\subsection{Finite element discretization}
\label{subsec:fem}
We first describe the finite element approximation to (\ref{dimensionlessPNP}) for a given mesh.
Denote the mesh nodes by
\[
0 = x_1 < x_2 < ... < x_{N_v-1} < x_{N_v} = 1,
\] 
where $N_v$ is the number of mesh nodes. Define
\[
V_h = \text{span}\{\psi_1,...,\psi_{N_v}\}, \quad
V_h^0 = \text{span}\{\psi_2,...,\psi_{N_v - 1}\} ,
\]
where $\psi_j = \psi_j(x)$ denotes the piecewise linear basis function corresponding to the node $x_j$.
A linear finite element discretization of (\ref{dimensionlessPNP}) and (\ref{bound}) is  to find
$\phi_h$, $c_{1,h}$, $c_{2,h} \in V_h$, subject to (\ref{bound}), such that
\begin{equation}
\begin{split}
& \ds\int_0^1\epsilon^2 h(x) \frac{d \phi_h}{ dx} \frac{d v}{d x}\,d x- \ds\int_0^1 h(x)\left ( \sum_{i = 1}^2 z_i c_{i,h} + Q(x)\right ) v\,d x = 0, \quad \forall v \in V_h^0 \\
& \ds\int_0^1 D_k h(x)\left (z_kc_{k,h}\frac{d \phi_h}{d x} + \frac{d c_{k,h}}{d x}\right )\frac{d v}{d x} \,d x = 0, \quad k = 1, 2,\quad \quad \forall v \in V_h^0 .
\end{split}
\label{pnp-fem}
\end{equation}
Here, the subscript ``$h$'' used in $\phi_h$, $c_{1,h}$, and $c_{2,h}$ is different from the function $h(x)$.
It is used to distinguish these discrete functions from their continuous counterparts.
Express $\phi_h$, $c_{1,h}$, and $c_{2,h}$ as
\[
\phi_h = \sum_{j=1}^{N_v} \phi_h^{(j)} \psi_j (x), \quad
c_{k,h} = \sum_{j=1}^{N_v} c_{k,h}^{(j)} \psi_j (x), \quad k = 1, 2.
\]
Using these and the boundary conditions (\ref{bound}) and taking $v = \psi_2$, ..., $\psi_{N_v - 1}$
in (\ref{pnp-fem}) sequentially, we obtain a system of nonlinear algebraic equations for
the unknown variables $\phi_h^{(j)}$, $c_{1,h}^{(j)}$, and $c_{2,h}^{(j)}$, $j = 1, ..., N_v$.
%
This system is solved by the MATLAB\textsuperscript \textregistered\, function {\em fsolve}, a nonlinear system solver based on the trust-region-Dogleg algorithm.

The computation alternates between the mesh generation (see the next subsection) and the solution of BVP (\ref{dimensionlessPNP}). More specifically, at the $n$th iteration, we assume that the physical mesh
$\mathcal{T}^n_h$ and the solution thereon are known. Then, a new mesh 
$\mathcal{T}_h^{n+1}$ is generated based on the solution on $\mathcal{T}^n_h$.
The new solution on the new mesh are then obtained by solving (\ref{pnp-fem}). This procedure is repeated
until convergence is reached. In practice, we have found that the convergence is reached very quickly
and the solution changes very little after a few iterations. We have stopped the computations in five iterations. 

\begin{figure}[htbp]
\subfigure[$c_1$ and $c_2$]{\includegraphics[width = 5cm,height = 4.2cm]{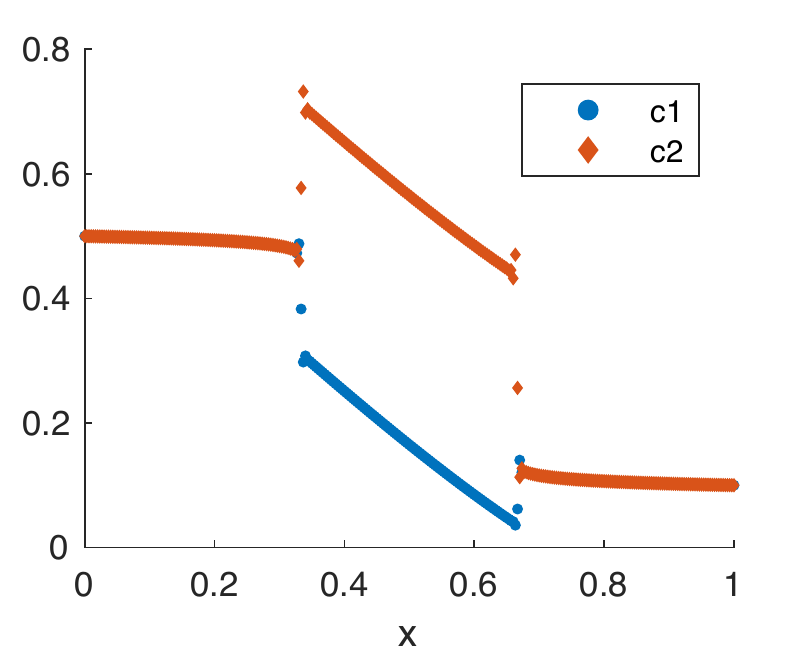} \label{fig:fixmesh-1}}
\subfigure[close view of (a)]{\includegraphics[width = 5cm,height = 4.2cm]{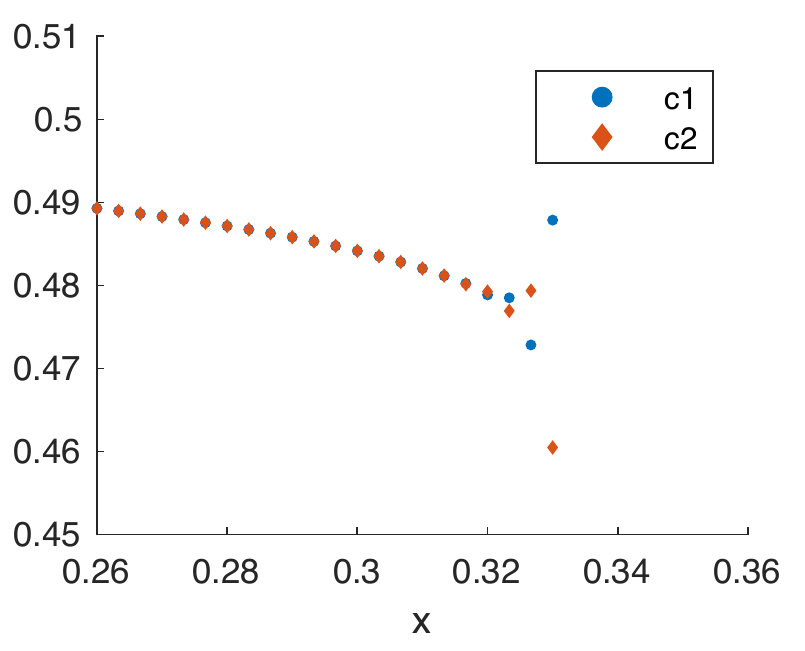} \label{fig:fixmesh-2}}
\subfigure[$J_1$ and $J_2$]{\includegraphics[width = 5cm,height = 4.2cm]{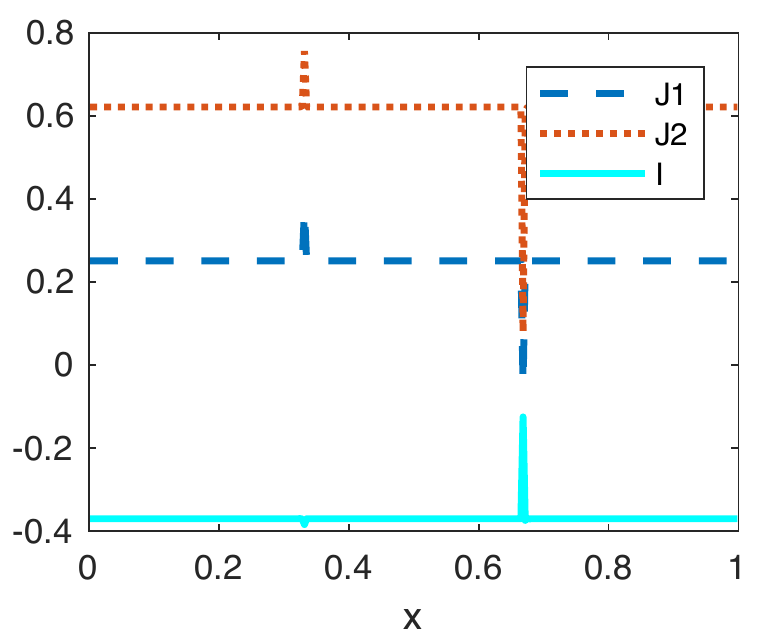} \label{fig:fixmesh-3}}
\subfigure[$c_1$ and $c_2$]{\includegraphics[width = 5cm,height = 4.2cm]{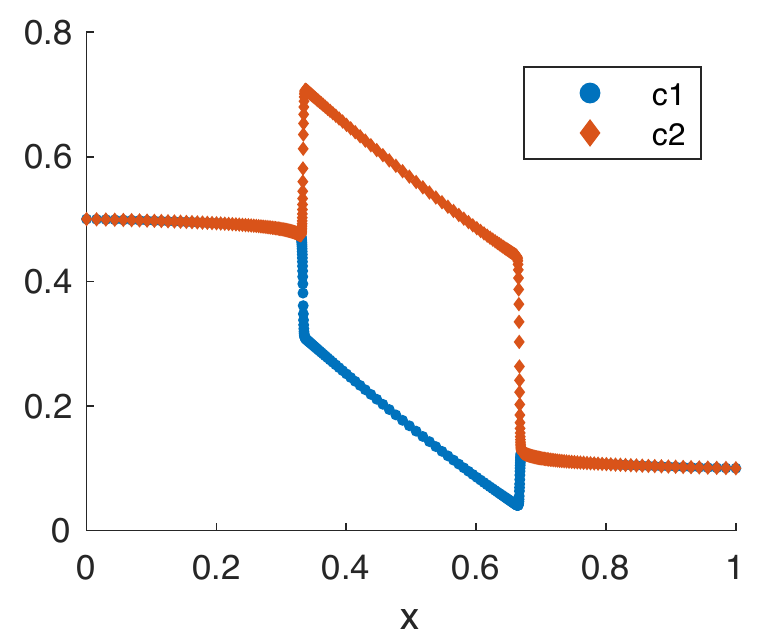} \label{fig:movingmesh-1}}
\hspace{1mm}
\subfigure[close view of (d)]{\includegraphics[width = 5cm,height = 4.2cm]{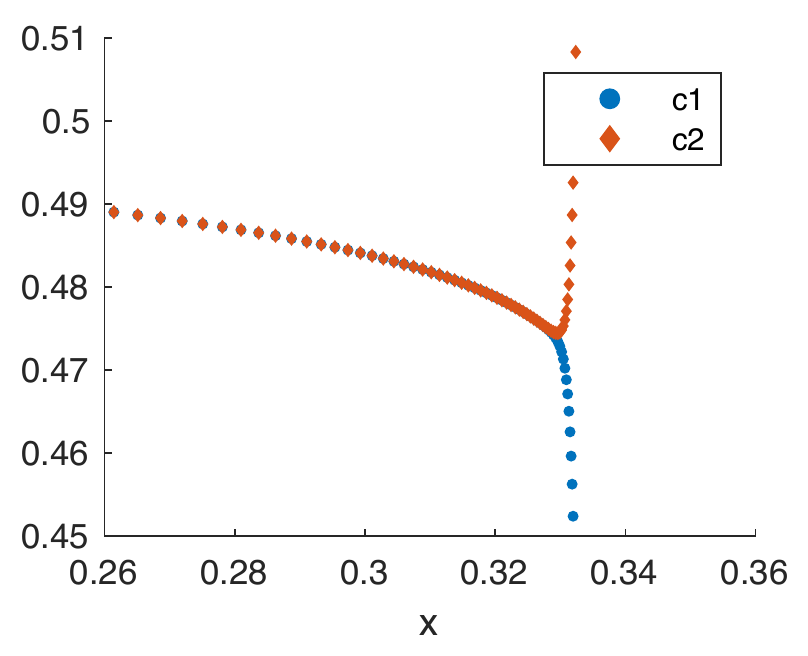} \label{fig:movingmesh-2}}
\hspace{1mm}
\subfigure[$J_1$ and $J_2$]{\includegraphics[width = 5cm,height = 4.2cm]{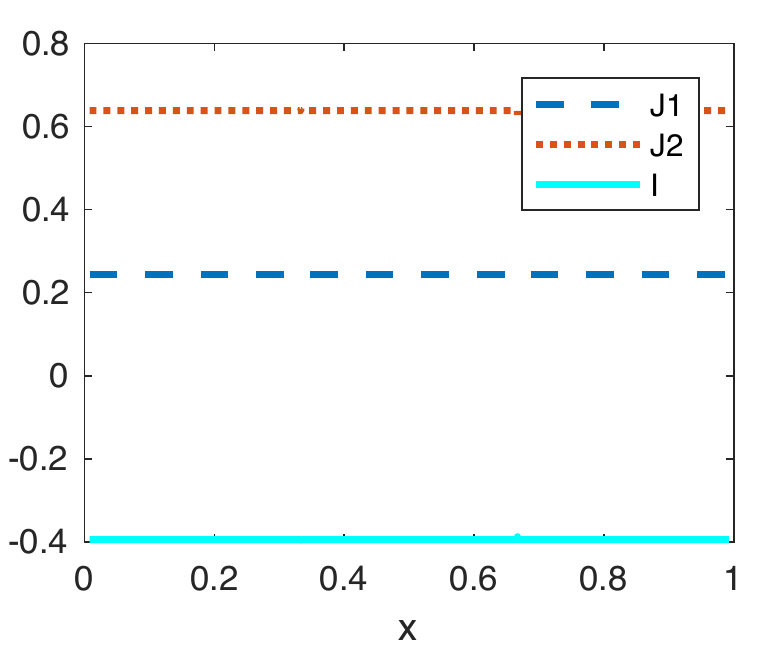} \label{fig:movingmesh-3}}
\caption{Solutions at the first iteration (without mesh adaptation) (a), (b), and (c) and the second iteration (with mesh adaptation), (d), (e), and (f). The second iteration is obtained with the adaptive mesh generated using the solutions from the first iteration.
The mesh size is $N_v = 301$. \label{fig:iteration}}
\end{figure}

Note that the analytical solution to BVP (\ref{dimensionlessPNP}) and (\ref{bound}) is not available for any
set of boundary values. Nevertheless, since the fluxes $J_1$ and $J_2$ should be constant throughout the domain,
we can check the accuracy of the computation by examining if the fluxes stay constant.
For example, Fig.~\ref{fig:fixmesh-1},~\ref{fig:fixmesh-2}, and \ref{fig:fixmesh-3} show the computational solution at the first iteration, which is obtained with the fixed mesh. Jumps in $J_1$ and $J_2$ are visible in the figures.
At the second iteration (with one round of mesh adaptation), as shown in Fig.~\ref{fig:movingmesh-1}, \ref{fig:movingmesh-2}, and \ref{fig:movingmesh-3}, the mesh is concentrated around $x = 1/3$ and $2/3$, as shown in the second panel, and the corresponding ionic flow $J_1$ and $J_2$ are almost constant.

\subsection{The MMPDE moving mesh method}
\label{SEC:MMPDE}
We now describe the MMPDE moving mesh method for the generation of the physical mesh $\mathcal{T}_h^{n+1}$ at the $(n + 1)$th iteration, given $\mathcal{T}_h^n$ and the solution to the system to (\ref{pnp-fem}) at the $n$th iteration. We introduce three different meshes: the physical mesh $\mathcal{T}_h = \{x_1,..., x_{N_v}\}$, the reference computational mesh $\hat{ \mathcal{T}}_{c,h} = \{ \hat{\xi}_1,..., \hat{\xi}_{N_v} \}$, and the computational mesh $\mathcal{T}_{c,h} = \{ \xi_1, ..., \xi_{N_v}\}$. All these three meshes have the same number of elements. $ \hat{ \mathcal{T}}_{c,h}$ is taken as a uniform mesh under Euclidean metric and is fixed during the computation, and $\mathcal{T}_{c,h}$ is introduced as an intermediate variable for computational purpose.

The MMPDE method employs a monitor function to control mesh concentration.
A choice for it is 
\begin{equation}
\rho(x) = (1 + |\phi''(x)|^2)^{\frac{1}{3}} ,
\label{meshdensity}
\end{equation}
where $\phi$ is the dimensionless electric potential.
The choice (\ref{meshdensity}) is known optimal for the $L^2$ norm of error for piecewise linear interpolation \cite{HR11}.
 This monitor function is expected to concentrate mesh elements in regions where the second-order
 derivative of $\phi$ is large.
 Although we expect that the second-order derivative of $\phi$ is large near the discontinuity points
 $x = 1/3$ and $2/3$, we modify the above monitor function for an extra explicit control of mesh concentration
 at these points as
\begin{equation}
\rho(x) = \sqrt{ (1 + |\phi''(x)|^2)^{\frac{2}{3}}
+ \frac{1}{e^{4(x - 1/3)^2} - 1 + c} + \frac{1}{e^{4(x - 2/3)^2} - 1 + c} },
\label{densityalt}
\end{equation}
where $c$ is a constant used to control the level of mesh concentration around the discontinuity points.
In our computation, $c$ has been taken as
\[
c = \frac{4}{\max\limits_{x}\, (1 + |\phi''(x)|^2)^{\frac{2}{3}}}.
\]
The exponential term in (\ref{meshdensity}) makes the level of concentration decay fast so that the mesh is not over-concentrated
near the discontinuity points.



The purpose of the MMPDE method is to generate a mesh in a way that the mesh density $\rho(x)$ is evenly distributed, i.e., 
\[
\int_{x_1}^{x_2} \rho(x) d x = \cdots = \int_{x_{N_v - 1} }^{x_{N_v}} \rho(x) d x,
\]
or
\begin{equation}
\int_{x_j}^{x_{j+1}} \rho(x) d x = \frac{\sigma}{N_v-1}, \qquad j = 1, ..., N_v - 1
\label{equidistribution-0}
\end{equation} 
where $\sigma = \int_0^1 \rho(x) dx$.
In practical computation, an approximation to the second derivative of $\phi$ at any mesh node is obtained
by first fitting a quadratic polynomial in a least-squares manner based on the nodal values of the computation
solution $\phi_h$ at the neighboring nodes and then differentiating the polynomial twice.
After that, a piecewise constant approximation to $\rho = \rho(x)$ is obtained. We denote
this approximation on $I_j = (x_{j-1}, x_{j})$ by $\rho_{I_j}$.
Then, a discrete version of (\ref{equidistribution-0}) is given by
\begin{equation}
(x_{j+1} - x_j) \rho_{I_j}  = \frac{\sigma_h}{N_v-1}, \qquad j = 1, ..., N_v - 1
\label{equidistribution}
\end{equation} 
where $\sigma_h = \sum_{j=1}^{N_v-1} \rho_{I_j} (x_{j+1} - x_j)$.
An energy function (or called the cost function) for (\ref{equidistribution}) is known \cite{HR11} to be
\begin{equation}
I_h(x_2, ..., x_{N_v-1}; \xi_2, ..., \xi_{N_v-1}) = \frac{1}{2} \sum_{n = 1}^{N_v-1} \frac{(\xi_{j+1} - \xi_{j})^2}{\rho_{I_j} (x_{j+1} - x_{j})^2} (x_{j+1} - x_{j}).
\label{energyformdis}
\end{equation}
There exist two approaches for minimizing this energy function. The first is a direct approach with which we take
$\mathcal{T}_{c,h}$ as $\hat{ \mathcal{T}}_{c,h}$ and then minimize $I_h$ with respect to $x_2, ..., x_{N_v-1}$.
In this case, the energy function becomes
\[
I_h(x_2, ..., x_{N_v-1}) = \frac{1}{2} \sum_{n = 1}^{N_v-1} \frac{(\hat{\xi}_{j+1} - \hat{\xi}_{j})^2}{\rho_{I_j} (x_{j+1} - x_{j})^2} (x_{j+1} - x_{j}).
\]
A disadvantage of this approach is that $\rho_{I_j}$ is highly nonlinear about the coordinates
of the physical nodes (i.e., $x_2, ..., x_{N_v-1}$),
which makes the minimization difficult. The other approach is an indirect one:
we take $\mathcal{T}_{h}$ as $\mathcal{T}_{h}^n$ and then minimize $I_h$ with respect to
$\xi_2, ..., \xi_{N_v-1}$. The energy function reduces to
\[
I_h(\xi_2, ..., \xi_{N_v-1}) = \frac{1}{2} \sum_{n = 1}^{N_v-1} \frac{(\xi_{j+1} - \xi_{j})^2}{\rho_{I_j} (x_{j+1}^n - x_{j}^n)^2} (x_{j+1}^n - x_{j}^n).
\]
Since $\rho_{I_j}$'s depend only on $x_1^n,\,x_2^n, ..., x_{N_v}^n$, this approach avoids the difficulty of $\rho$ being a highly nonlinear function of the physical coordinates.  But it
does not compute the physical mesh directly. To obtain the new physical mesh, we denote the new
computational mesh obtained through minimization by $\mathcal{T}_{c,h}^{n+1}$.
Then $\mathcal{T}_h^n$ and $\mathcal{T}_{c,h}^{n+1}$ form a correspondence which can be written
formally as $\mathcal{T}_h^n =  \Psi_h(\mathcal{T}_{c,h}^{n+1})$ or $x_j^n = \Psi_h(\xi_j^{n+1}),\; j = 1, ..., N_v$.
Then the new physical mesh is defined as $\mathcal{T}_h^{n+1} = \Psi_h(\hat{\mathcal{T}}_c)$
or $x_j^{n+1} = \Psi_h(\hat{\xi}_j),\; j = 1, ..., N_v$,
and can be computed readily via linear interpolation.

The minimization of $I_h$ with respect to $\xi_2, ..., \xi_{N_v-1}$ for the current one-dimensional situation
leads to a linear system which can readily be solved. On the other hand, the MMPDE method, which has been developed
for generating adaptive meshes in one and multiple dimensions, solves
the minimization problem by integrating the gradient system of $I_h$, i.e., 
\begin{equation}
\frac{\partial \xi_i}{\partial t} = - \frac{\partial I_h}{\partial \xi_j}, \quad j = 2,...,N_v - 1
\label{ximethod}
\end{equation}
where the partial derivative is
\begin{equation}
\frac{\partial I_h}{\partial \xi_j} = \frac{\xi_j - \xi_{j-1}}{\rho_n(x_j^n - x_{j-1}^n)} - \frac{\xi_{j+1} - \xi_j}{\rho_{j+1}(x_{j+1}^n - x_j^n)}. 
\label{Lagrange}
\end{equation}
The mesh equation (\ref{ximethod}) is integrated using the Matlab\textsuperscript \textregistered\, function {\em ode15s}, a Numerical Differentiation Formula based integrator, to a preset time $t = 10$ with the initial mesh $\hat{\mathcal{T}}_{c,h}$.

The MMPDE approach (\ref{ximethod}) performs very comparably with the direct minimization of $I_h$ with respect to $\xi_2, ..., \xi_{N_v-1}$ in one dimension. However, in multi-dimensions, $I_h$ has a more complex structure and its direct minimization with respect to the coordinates of the computational vertices will no longer result in a linear algebraic system. In this case, the MMPDE method has the advantages of being more stable in finding
approximate minimizers and avoiding the mesh from tangling and crossing over \cite{HK2015}.

\section{Conclusions}
In this work we have studied the effects of permanent charge along with boundary conditions on ionic flows via the numerical solution of a quasi-one-dimensional PNP model. The studies focus on the flux ratios $\lambda_k$'s 
defined in (\ref{lambdak}) that indicate if the fluxes $J_k$'s are enhanced by the permanent charge and boundary conditions. For an ionic mixture of two ion species with $z_1 = 1$ and $z_2 = -1$ and for fixed boundary ion concentrations and fixed shape of the ion channel, the behaviors of $\lambda_1$ and $\lambda_2$ as functions of $Q_0$ (for the permanent charge in the form (\ref{permanentcharge})) and the transmembrane electric potential $V$ (cf. (\ref{bound})) have been analyzed in \cite{JLZ2015, Liu18, ZEL2019,ZL2020} (also see \S\ref{SEC:recentresults}) for small and relatively large $Q_0$. In particular, for small $Q_0$, as $V$ increases,
$\lambda_1$ and $\lambda_2$ go through regions with $\lambda_1 < \lambda_2 < 1$
to $\lambda_1 < 1 < \lambda_2$ and to $1 < \lambda_1 < \lambda_2$, or vice versa.
On the other hand, for large $Q_0$, we have $\lambda_1 < 1$ for all $V$ but
$\lambda_2$ can be greater than one for a bounded interval of $V$ and less than one otherwise.
This difference means that there are transitions in the behaviors of $\lambda_1$ and $\lambda_2$
from small and relatively large $Q_0$ for which it is difficult, if not impossible, to study using the existing
analytical techniques. To study these transitions, we have used a numerical tool, an adaptive moving mesh
finite element method that can provide better resolution of the solution at the endpoints of the ``neck" of
an ion channel.

For fixed boundary concentrations, we have obtained a complete diagram, divided into different regions of $Q_0$-$V$ plane, that have verified and matched the existing analytical results for small and relatively large $Q_0$. We have observed saddle-node bifurcations, which might be useful reference in practice, and hopefully will bring up interesting subjects for future studies. Similar bifurcation diagrams have been obtained for different boundary conditions for the concentrations. Moreover, we have provided bifurcation diagrams corresponding to the quasi one-dimensional PNP model with a hard sphere model for the excess electrochemical potential. Similar bifurcation properties but with quantitative difference have been observed.  It remains an open question how these bifurcations have been generated, and further studies in mathematics as well as biology are highly demanded.
\bigskip

\noindent
{\bf Acknowledgement.} WL was partially supported by Simons Foundation 
Mathematics and Physical Sciences-Collaboration Grants for Mathematicians \#581822.

\end{document}